\documentclass[11pt,a4paper]{article}
\usepackage{lmodern}
\usepackage[T1]{fontenc}
\usepackage[utf8]{inputenc}
\usepackage[english]{babel}
\usepackage{fancyhdr}
\usepackage{setspace}
\usepackage{hyperref} 

\usepackage{amsmath}
\usepackage{amsfonts}
\usepackage{amssymb}
\usepackage{amsthm}
\usepackage{tikz}
\usepackage{subfig}

\newtheorem{theorem}{Theorem} 

\newtheorem{remarque}[theorem]{Remark}

\newcommand{\tr}{^{\mathsf{T}}}
\newcommand{\Real}{\mathbb{R}}
\newcommand{\unew}{^{\textrm{new}}}
\newcommand{\uold}{^{\textrm{old}}}
\newcommand{\veps}{\varepsilon}
\newcommand{\vesp}{\veps_p}
\newcommand{\vespc}{\veps_{p,c}}
\newcommand{\vespC}{\veps_{p,\mathcal{C}}}

\newcommand\ale[1]{%
  \protect\leavevmode
  \begingroup
        \color{red!65!black}%
    #1%
  \endgroup
}

\title{A variational discrete element method for quasi-static and dynamic elasto-plasticity}
\author{\begin{minipage}{\textwidth}\centering Frédéric
Marazzato$^{1,2,3}$, Alexandre Ern$^{1,3}$ and Laurent Monasse$^{1,4}$\\
   \small{$^{1}$Universit\'e Paris-Est, Cermics (ENPC), F-77455 Marne-la-Vall\'ee cedex 2, France}\\
   \small{email: \texttt{\{alexandre.ern, frederic.marazzato\}@enpc.fr}}\\
   \small{$^{2}$CEA, DAM, DIF, F-91297 Arpajon, France}\\
   \small{$^{3}$Inria Paris, EPC SERENA, F-75589 Paris, France}\\
   \small{$^4$ Universit\'e C\^ote d'Azur, Inria, CNRS, LJAD, EPC COFFEE, 06108 Nice, France}\\
   \small{email: \texttt{laurent.monasse@inria.fr}}\end{minipage}}
   
\begin{document}

\maketitle

\begin{abstract}
We propose a new discrete element method supporting general polyhedral meshes. The method can be understood as a lowest-order discontinuous Galerkin method parametrized by the continuous mechanical parameters (Young's modulus and Poisson's ratio). We consider quasi-static and dynamic elasto-plasticity, and in the latter situation, a pseudo-energy conserving time-integration method is employed.
The computational cost of the time-stepping method is moderate since it is explicit and used with a naturally diagonal mass matrix. 
Numerical examples are presented to illustrate the robustness and versatility of the method for quasi-static and dynamic elasto-plastic evolutions.
\end{abstract}

\section{Introduction}

Discrete element methods (DEM) constitute a large class of particle
methods which have originally been used for crystalline materials
\cite{hoover1974two} and geotechnical problems
\cite{cundall1979discrete} and have found applications in granular
materials, soil and rock mechanics. In their original formulation,
DEM consisted in representing a
domain by small spherical particles interacting by means of forces and
torques.  A wide range of models for the expression of these bonds
has been developed depending on the material
constitutive law. Computing the deformation of the domain then
consists in computing the evolution of the particle system. 
Advantages of DEM are their ability to deal with discontinuous
materials, such as fractured or porous materials, as well as the
possibility to take advantage of GPU computations
\cite{MR3606237}. Other, similar, particle methods have been derived in the
context of Smooth Particles Hydrodynamics (SPH) methods, which require an
interaction kernel \cite{MR1857694}. The main difficulty in DEM
consists in deriving a correct set of forces between elements to
discretize the continuous equations (in the present case, dynamic
elasto-plasticity). DEM originally used sphere
packing to discretize the domain \cite{MR2548422} and were forced to
fit parameters in order to obtain relevant values for the Young modulus $E$ or
the Poisson ratio $\nu$ \cite{jebahi2015discrete,MR3735741}.
Moreover simulating a material with a Poisson
ratio $\nu$ larger than $0.3$ met with difficulties \cite{andre2013using}. Note also the
possibility to use DEM only in a limited zone, where crack occurs for
instance, in order to mitigate these issues. For example, a modified
DEM (MDEM) has been coupled with a virtual element method (VEM)
for elasticity to discretize fracturing porous media
\cite{MR3736041}.

A discrete element method was developed in
\cite{LM_CM_2012} and was formally proved to be consistent with Cauchy
elasticity. A first attractive feature of this method was
that the discrete force parameters were directly
derived from the Young modulus and the Poisson ratio without the need for a fitting
process. Moreover the method exhibited robustness in the incompressibility limit
($\nu\to 0.5$). Similar ideas have been used to handle brittle
fracture \cite{andre2019novel}. However several limitations remain in
this approach. First the evaluation of the forces between particles 
hinges on the use of a
Voronoi mesh and does not adapt to general (not even tetrahedral)
meshes. This is due to a nearest-neighbour evaluation of the gradient on a
facet of the mesh (known in the Finite Volume community as the ``two-point flux
problem''). Secondly the expression of the forces for
a Cauchy continuum cannot be readily extended to more general behaviour
laws. Finally the convergence proof is mostly formal (on a Cartesian
grid) and no convergence proof is given on general (Voronoi) meshes,
apart from numerical evidence.

The main goal of the present contribution is to circumvent the above
issues by extending the discrete element
method of \cite{LM_CM_2012} to general polyhedral meshes and
elasto-plastic behaviour laws.
The present space-discretization scheme involves a diagonal mass matrix and shares a number of properties with finite volume \cite{eymard2004finite,eymard2009discretization} and lowest-order discontinuous
Galerkin (dG) methods \cite{di2012cell}. Specifically we use piecewise constant gradient 
reconstructions in each mesh cell evaluated from local displacement reconstructions at
the facets of the mesh. As a consequence, the present DEM is not a mesh-free method since it uses the mesh geometry to compute the above reconstructions. However, it can still be viewed as a particle method owing to the use of a diagonal mass matrix.
In addition to the displacement unknowns, volumetric unknowns representing plastic strains are also added. 
We devise the scheme for both quasi-static and dynamic elasto-plasticity, and in the latter
situation we perform the time discretization using the explicit pseudo-energy conserving
time-integration method developed in \cite{MARAZZATO2019906}.
There are two main differences with \cite{eymard2009discretization}
regarding the discretization of the elliptic (here linear elastic) part. On the one hand
\cite{eymard2009discretization} uses only cell dofs, whereas the present
DEM uses also boundary vertex dofs for the displacement. As we shall see below,
the introduction of boundary vertex dofs can have a sizable impact on tempering
the CFL restriction on the time step in the context of explicit time-marching methods.
On the other hand the mass distribution of the present DEM is different. 
Our choice is motivated by the fact that, to use an explicit integration with a diagonal mass matrix, every dof needs a mass to compute its velocity, even dofs on boundary vertices. 
Moreover numerical results
are presented to illustrate the robustness and versatility of the proposed method
in two and three space dimensions. Finally, we mention that the convergence of the scheme
can be studied using the framework of gradient discretization methods
(GDM) \cite{droniou2018gradient}. GDM lead to a unified and powerful framework
allowing one to prove convergence and error estimates for a wide range of
numerical schemes. 

This paper is organised as follows. Section~\ref{sec:governing} briefly 
recalls the equations of dynamic
elasto-plasticity in a Cauchy continuum. Section~\ref{sec:dofs} 
introduces the proposed DEM and presents
the space discretization of the governing equations.
Some numerical tests to verify the convergence of the space
discretization in a steady setting are reported. 
Section~\ref{sec:quasi-static} deals with the DEM discretization
for quasi-static elasto-plasticity and presents
test cases in two and three space dimensions.  Section~\ref{sec:MEMM}
addresses the time discretization of the dynamic elasto-plasticity problem
using the explicit pseudo-energy conserving time-integrator developed in
\cite{MARAZZATO2019906}. This section also assesses the coupled
DEM and time discretization on test cases in three space dimensions. 
Finally Section~\ref{sec:conclusions} draws some 
conclusions.

\section{Governing equations for dynamic elasto-plasticity}
\label{sec:governing}

We consider an elasto-plastic material occupying the domain
$\Omega \subset \mathbb{R}^d$, $d\in\{2,3\}$, in the reference configuration
and evolving dynamically on the finite time interval $(0,T)$, $T > 0$, under the 
action of volumetric forces and boundary conditions.
The strain regime is restricted to small strains so that the linearized strain tensor
is $\veps(u):=\frac12(\nabla u+(\nabla u)\tr)$, where $u$ is the displacement field.
The plastic constitutive law hinges on a von Mises criterion with nonlinear isotropic hardening. 
The material is supposed to be homogeneous, isotropic and rate-independent. 
The present formalism can be extended to the case of anisotropic, inhomogeneous, 
rate-dependent, anisothermal materials as well as finite strains.
The stress tensor $\sigma \in \mathbb{R}_\mathrm{sym}^{d\times d}$ is such that
\begin{equation}
\sigma := \mathbb{C}: (\veps(u) - \vesp), 
\end{equation}
where $\mathbb{C}$ is the fourth-order stiffness tensor, the subscript $\mathrm{sym}$ stands for symmetric tensors and $\vesp \in \mathbb{R}_\mathrm{sym}^{d\times d}$ is 
the (trace-free) tensor of remanent plastic strain. The 
von Mises yield function $\varphi$ is given by
\begin{equation}
\label{eq:yield}
\varphi(\sigma, p) := \sqrt{\frac{3}{2}} \vert \mathrm{dev}(\sigma) \vert - (\sigma_0+R(p)),
\end{equation} 
where $\mathrm{dev}(\sigma)$ is the deviatoric part of $\sigma$
and $|\tau|=(\sum_{i,j=1}^d \tau_{ij}^2)^{\frac12}$ for a second-order tensor $\tau$, 
$p$ is the scalar cumulated plastic deformation, 
$R(p) := \frac{d \omega_p}{dp}$ where the function $\omega_p$ is the part of the Helmholtz free 
energy related to isotropic hardening, and $\sigma_0$ is the initial yield stress, so that
the actual yield stress is $\sigma_0 + R(p)$. Admissible states are characterized by the inequality
$\varphi(\sigma,p) \leq 0$, the material is in the elastic domain if $\varphi(\sigma,p) < 0$
and in the plastic domain if $\varphi(\sigma,p) = 0$.

In strong form, the dynamic elasto-plasticity equations consist in searching 
for the displacement field $u:(0,T) \times \Omega\rightarrow \mathbb{R}^d$, 
the remanent plastic strain tensor $\varepsilon_p:(0,T) \times 
\Omega\rightarrow \mathbb{R}^{d\times d}_{\mathrm{sym}}$, and the scalar cumulated plastic deformation
$p:(0,T) \times \Omega\rightarrow \mathbb{R}$ such that the following equations hold in $\Omega$
for all $t\in(0,T)$:
\begin{equation}
\label{eq:elasto-plasticity}
\left\{ \begin{aligned}
&\rho \ddot{u} - \mathop{\rm div}(\sigma) = f, \\
&\lambda \ge 0, \quad \varphi(\sigma,p) \leq 0, \quad \lambda \varphi(\sigma, p) = 0\\
&\dot p = \lambda, \quad
\dot{\varepsilon}_p = \lambda \frac{\partial \varphi}{\partial \sigma}(\sigma),
\end{aligned} \right.
\end{equation}
where $\rho > 0$ is the density of the material,
dots indicate time derivatives, 
$f$ is the imposed volumetric force, and $\lambda$ is the Lagrange
multiplier associated with the constraint $\varphi(\sigma,p) \leq 0$. 
Note that owing to~\eqref{eq:yield}, we have $\dot{\varepsilon}_p=\lambda
\sqrt{\frac32} \frac{\mathrm{dev}(\sigma)}{|\mathrm{dev}(\sigma)|}$, so that
$\dot{p} = \lambda = \sqrt{\frac{2}{3}} |\dot{\varepsilon}_p|$.

Let $\partial \Omega = \partial \Omega_N \cup \partial \Omega_D$ be a partition of 
the boundary of $\Omega$. By convention $\partial\Omega_D$ is a closed set and 
$\partial\Omega_N$ is a relatively open set in $\partial\Omega$. 
The boundary $\partial \Omega_D$ has an imposed displacement $u_D$, 
whereas a normal stress $g$ is imposed on $\partial \Omega_N$, i.e. we enforce
\begin{equation} \label{eq:BC}
u = u_D \ \text{ on } (0,T) \times \partial \Omega_D, \qquad
\sigma \cdot n = g_N\  \text{ on } (0,T) \times \partial \Omega_N.
\end{equation}
Note that $u_D$ and $g_N$ can  be time-dependent.
Finally the initial conditions prescribe that $u(0) = u_0$, $\dot{u}(0) = v_0$
and $p(0)=0$ in $\Omega$. 

To formulate the governing equations (\ref{eq:elasto-plasticity}) together with the
Neumann boundary condition from \eqref{eq:BC} in weak form,
we consider time-dependent functions with values in space-dependent functional
spaces. Let us set 
\begin{equation}
V_D := \left\{ v \in H^1(\Omega;\mathbb{R}^d) \ | \ v_{|\partial \Omega_D} = u_D \right\},
\qquad
V_0 := \left\{ v \in H^1(\Omega;\mathbb{R}^d) \ | \ v_{|\partial \Omega_D} = 0 \right\}.
\end{equation}
(Note that the space $V_D$ is actually time-dependent if 
the Dirichlet data is time-dependent.) We also set
\begin{equation} 
Q:= L^2(\Omega;\mathbb{R}^{d \times d}_{\mathrm{sym}}), \qquad Q_0 := \left\{ \eta_p \in Q \ | \ 
\mathrm{tr}(\eta_p) = 0 \right\},
\end{equation}
where $P:=L^2(\Omega)$. 
Here, for any vector space S, $L^2(\Omega;S)$ is the Hilbert space composed of
$S$-valued square-integrable functions in $\Omega$, and 
$H^1(\Omega;S)$ is the subspace of $L^2(\Omega;S)$ composed
of those functions whose weak gradient is also square-integrable. All of the above
functional spaces are equipped with their natural inner product.
Then the weak solution is searched as a triple 
$(u,\varepsilon_p,p):(0,T) \rightarrow V_D\times Q_0\times P$.
To alleviate the mathematical
formalism, we do not specify here the regularity in time (see 
\cite{han2012plasticity} and \cite{carstensen1999numerical}). We introduce the
mass bilinear form such that
\begin{equation}
\label{eq:kinetic power}
m(a,\tilde{v}) := \left< \rho a, \tilde{v} \right>_{V_0',V_0},
\quad \forall (a,\tilde{v}) \in V_0'\times V_0,
\end{equation}
where $V_0'$ denotes the dual space of $V_0$ and $\left< \cdot,\cdot \right>_{V_0',V_0}$
the duality pairing,
and the stiffness bilinear form parameterized by a member $\eta_p\in Q_0$ such that
\begin{equation}
\label{eq:bilinear form}
a(\eta_p;v,\tilde{v}) := \int_{\Omega} \big( \mathbb{C} : (\varepsilon(v) - \eta_p) \big) : \varepsilon(\tilde{v}), 
\quad \forall (v,\tilde{v}) \in V_D\times V_0,
\end{equation}
The governing equations (\ref{eq:elasto-plasticity}) are rewritten as follows:
Find $(u,\varepsilon_p,p):(0,T) \rightarrow V_D\times Q_0\times P$ such that, for all $t\in (0,T)$,
\begin{equation}
\label{eq:elasto-plasticity-weak}
\left\{ \begin{alignedat}{2}
&m(\ddot{u}(t),\tilde v) +a(\varepsilon_p(t);u(t),\tilde{v}) = l(t;\tilde{v}),
&\quad&\forall \tilde{v}\in V_0,\\
&\lambda\ge 0, \quad \varphi(\sigma,p) \leq 0, \quad 
\lambda \varphi(\sigma, p) = 0,&\quad&\text{in $\Omega$},\\
&\dot p = \lambda, \quad 
\dot{\varepsilon}_p = \lambda \frac{\partial \varphi}{\partial \sigma}(\sigma),&\quad&\text{in $\Omega$},
\end{alignedat} \right.
\end{equation}
where the time-dependency is left implicit in the second and third equations, and
with the linear form $l(t;\cdot)$ acting on $V_0$ as follows:
\begin{equation}
\label{eq:linear form}
l(t;\tilde{v}) := \int_{\Omega} f(t) \cdot \tilde{v} + \int_{\partial \Omega_N} g_N(t) \cdot \tilde{v}. 
\end{equation}
Note that the Dirichlet condition is enforced strongly, whereas the Neumann condition
is enforced weakly. 
Define the elastic energy $E_{\mathrm{elas}}(t):=\frac12 \int_\Omega \sigma(t) : 
\mathbb{C}^{-1} : \sigma(t)$ with $\sigma(t):=\mathbb{C}:(\veps(u(t))-\vesp(t))$, the kinetic energy $E_{\mathrm{kin}}(t):=\frac12 m(\dot u(t),\dot u(t))$, the plastic dissipation
$E_{\mathrm{plas}}(t) := \int_\Omega \sigma_0p(t)+\omega_p(p(t))$, and the work of external loads
$E_{\mathrm{ext}}(t) := \int_0^t l(s;\dot u(s))ds$. Then assuming for simplicity a homogeneous Dirichlet condition, and recalling the assumption $p(0) = 0$, we have the following energy equation:
\begin{equation} \label{eq:energy_balance} E_{\mathrm{elas}}(t) + E_{\mathrm{kin}}(t) + E_{\mathrm{plas}}(t)
= E_{\mathrm{ext}}(t) + E_{\mathrm{elas}}(0) + E_{\mathrm{kin}}(0),
\end{equation}
showing that the total energy at time $t$ is balanced with the work of external loads up to time $t$ (notice that $E_{\mathrm{plas}}(0)=0$ in our setting).

\section{Space semi-discretization}
\label{sec:dofs}

In this section we present the DEM space semi-discretization 
of the weak formulation~\eqref{eq:elasto-plasticity-weak}, and we present a few
verification test cases for static linear elasticity.

\subsection{Degrees of freedom}

The domain $\Omega$ is discretized with a
mesh $\mathcal{T}_h$ of size $h$ made of polyhedra with planar facets
in three space dimensions or polygons with straight edges in two space dimensions. 
We assume that $\Omega$ is itself a
polyhedron or a polygon so that the mesh covers $\Omega$ exactly, and 
we also assume that the mesh
is compatible with the partition of the boundary $\partial\Omega$ 
into the Dirichlet and Neumann parts. 

\begin{figure} [!htp] 
\begin{center}
\begin{tikzpicture}[scale=1.3]
\coordinate (a) at (0,-1);
\coordinate (b) at (0,1);
\coordinate (c) at (1.3,1);
\coordinate (d) at (-1.7,-0.8);
\coordinate (e) at (-2,-3);
\coordinate (f) at (0.5,-2.5);
\coordinate (g) at (1.7,0);
\coordinate (h) at (2.8,0.1);
\coordinate (i) at (2.1,-1.1);
\coordinate (j) at (3.5,-1.5);
\coordinate (k) at (2.2,-3.5);
\coordinate (l) at (0.8,-4);
\coordinate (m) at (-1.7,-4);

\draw[fill=gray,opacity=0.1] (b) --(c) -- (h) -- (j) -- (k) -- (l) --
(m) -- (e) -- (d) -- cycle;
\node[right] at (barycentric cs:f=0.8,a=0.25,g=0.25,i=0.25) {\Large
  $\Omega$};
\node[right] at (barycentric cs:k=0.5,j=0.5) {\Large $\partial\Omega$};
\path[draw] (a)-- (b)-- (c) -- (g) -- cycle;
\path[draw] (a)-- (b)-- (d) -- cycle;
\path[draw] (a) -- (e) -- (d) -- cycle;
\path[draw] (a) -- (f) -- (e) -- cycle;
\path[draw] (a) -- (g) -- (i) -- (f) -- cycle;
\path[draw] (c) -- (h) -- (g) -- cycle;
\path[draw] (g) -- (h) -- (j) -- (i) -- cycle;
\path[draw] (i) -- (j) -- (k) -- (l) -- (f) -- cycle;
\path[draw] (f) -- (l) -- (m) -- (e) -- cycle;

\fill[red] (barycentric cs:a=0.3,b=0.3,d=0.3) circle (3pt);
\fill[red] (barycentric cs:a=0.3,e=0.3,d=0.3) circle (3pt);
\fill[red] (barycentric cs:f=0.3,e=0.3,a=0.3) circle (3pt);
\fill[red] (barycentric cs:a=0.25,g=0.25,c=0.25,b=0.25) circle (3pt);
\fill[red] (barycentric cs:g=0.3,h=0.3,c=0.3) circle (3pt);
\fill[red] (barycentric cs:a=0.25,g=0.25,i=0.25,f=0.25) circle (3pt);
\fill[red] (barycentric cs:h=0.25,g=0.25,i=0.25,j=0.25) circle (3pt);
\fill[red] (barycentric cs:i=0.2,j=0.2,k=0.2,l=0.2,f=0.2) circle
(3pt);
\fill[red] (barycentric cs:f=0.25,l=0.25,m=0.25,e=0.25) circle (3pt);

\node at (b) [fill=blue,minimum width=3pt,minimum height=3pt] {};
\node at (c) [fill=blue,minimum width=3pt,minimum height=3pt] {};
\node at (d) [fill=blue,minimum width=3pt,minimum height=3pt] {};
\node at (e) [fill=blue,minimum width=3pt,minimum height=3pt] {};
\node at (h) [fill=blue,minimum width=3pt,minimum height=3pt] {};
\node at (j) [fill=blue,minimum width=3pt,minimum height=3pt] {};
\node at (k) [fill=blue,minimum width=3pt,minimum height=3pt] {};
\node at (l) [fill=blue,minimum width=3pt,minimum height=3pt] {};
\node at (m) [fill=blue,minimum width=3pt,minimum height=3pt] {};

\fill[red] (4.6,-1) circle (3pt);
\node[right] at (4.7,-1) {$(u_c)_{c\in\mathcal{C}}$};
\node at (4.6,-2) [fill=blue,minimum width=3pt,minimum height=3pt]{};
\node[right] at (4.7,-2) {$(u_z)_{z\in\mathcal{Z}^\partial}$};
\end{tikzpicture}
\caption{Continuum $\Omega$ covered by a polyhedral mesh and vector-valued degrees of freedom for the displacement.}
\label{fig:dofs}
\end{center}
\end{figure}
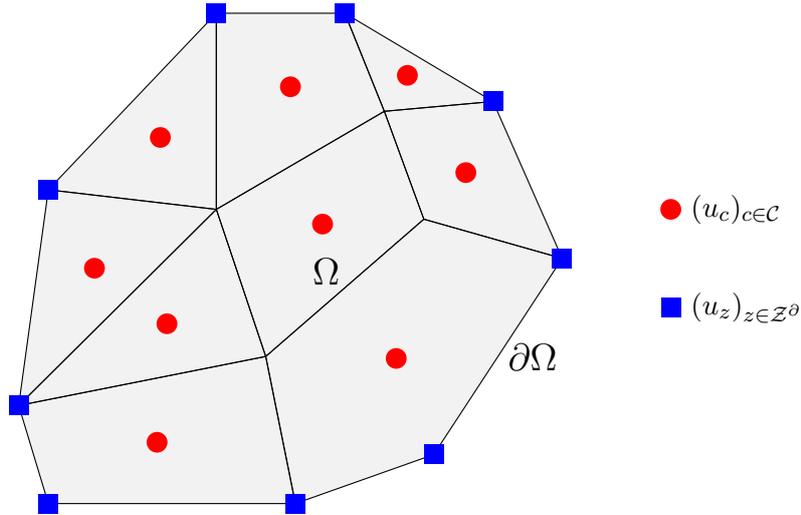

Let $\mathcal{C}$ denote the
set of mesh cells and
$\mathcal{Z}^\partial$ the set of mesh vertices sitting on the boundary of $\Omega$.  
Vector-valued volumetric degrees of freedom (dofs) for a generic 
displacement field $v_{\mathcal{C}}:=(v_c)_{c\in\mathcal{C}}\in\mathbb{R}^{d\#(\mathcal{C})}$ 
are placed at the barycentre of
every mesh cell $c\in\mathcal{C}$, where $\#(S)$ denotes the cardinality of any set $S$. 
Additional vector-valued boundary degrees of freedom
$v_{\mathcal{Z}^\partial}:=(v_z)_{z\in\mathcal{Z}^\partial}\in\mathbb{R}^{d\#(\mathcal{Z}^\partial)}$ 
for the displacement are added at
every boundary vertex $z\in\mathcal{Z}^\partial$. The reason why we introduce boundary vertex
dofs is motivated in Section~\ref{sec:reconstruction}. These 
dofs are also used to enforce the Dirichlet condition on $\partial\Omega_D$.
We use the compact notation $v_h:=(v_{\mathcal{C}},v_{\mathcal{Z}^\partial})$ for
the collection of all the cell dofs and all the boundary vertex
dofs. Figure \ref{fig:dofs} illustrates the position of the displacement
dofs. In
addition a (trace-free) symmetric tensor-valued dof representing the internal plasticity variable
$\eta_{p,c}\in\mathbb{R}^{d\times d}_{\mathrm{sym}}$ is attached to every mesh cell $c\in\mathcal{C}$, as well
as a scalar dof $p_c$ representing the cumulated plastic deformation. 
We write $\eta_{p,\mathcal{C}}:=(\eta_{p,c})_{c\in\mathcal{C}}\in Q_h := 
\left( \mathbb{R}^{d\times d}_\mathrm{sym} \right)^{\#(\mathcal{C})}$ and $p_{\mathcal{C}}:=(p_c)_{c\in\mathcal{C}}\in P_h := 
\mathbb{R}^{\#(\mathcal{C})}$.

Let $\mathcal{F}$ denote the set of mesh facets. We partition 
this set as $\mathcal{F}=\mathcal{F}^i\cup \mathcal{F}^b$, where 
$\mathcal{F}^i$ is the collection of the internal facets shared
by two mesh cells and $\mathcal{F}^b$ is the collection of the boundary facets 
sitting on the boundary $\partial\Omega$ (such facets belong to the boundary of only
one mesh cell). 
Using the cell and boundary-vertex dofs introduced above, we reconstruct a collection of
displacements $v_{\mathcal{F}}:=(v_F)_{F\in\mathcal{F}}\in
\mathbb{R}^{d\#(\mathcal{F})}$ on the mesh facets. 
The facet reconstruction operator is
denoted $\mathcal{R}$ and we write
\begin{equation}
v_{\mathcal{F}} := \mathcal{R}(v_h).
\end{equation}
The precise definition of the facet reconstruction operator is given in
Section~\ref{sec:reconstruction}. 
Using the reconstructed facet displacements and a discrete Stokes formula, it is possible to devise a discrete $\mathbb{R}^{d\times d}$-valued
piecewise-constant gradient field for the displacement
that we write $G_{\mathcal{C}}(v_{\mathcal{F}}) := (G_{c}(v_{\mathcal{F}}))_{c\in\mathcal{C}}
\in \mathbb{R}^{d^2\#(\mathcal{C})}$. Specifically we set in every
mesh cell $c\in\mathcal{C}$,
\begin{equation}
\label{eq:gradient reconstruction}    
G_c(v_{\mathcal{F}}) := \sum_{F \in \partial c} \frac{|F|}{|c|} v_F \otimes n_{F,c},
\qquad \forall v_{\mathcal{F}} \in\mathbb{R}^{d\#(\mathcal{F})},
\end{equation}
where the summation is over the facets $F$ of $c$ and $n_{F,c}$ is the outward 
normal to $c$ on $F$. 
Note that for all $v_h\in V_h$, we have
\begin{equation}
G_c(\mathcal{R}(v_h)) = \sum_{F \in \partial c} \frac{|F|}{|c|} (\mathcal{R}(v_h)_F -
v_c) \otimes n_{F,c},
\end{equation} 
since $\sum_{F \in \partial c}|F|n_{F,c}=0$.
We define a constant
linearized strain tensor in every mesh cell $c\in\mathcal{C}$ such that 
\begin{equation}
\varepsilon_{c}(v_{\mathcal{F}}) := \frac{1}{2}(G_{c}(v_{\mathcal{F}})+
G_{c}(v_{\mathcal{F}})\tr)\in \mathbb{R}^{d\times d}_{\mathrm{sym}}.
\end{equation}
Finally, we define two additional reconstructions. The first is a cellwise nonconforming $P^1$ reconstruction $\mathfrak{R}$ defined for all $c \in \mathcal{C}$ by
\begin{equation}
\label{eq:DG reconstruction}
\mathfrak{R}(v_h)_c(\mathbf{x}) := v_c + G_c(\mathcal{R}(v_h)) \cdot (\mathbf{x} - \mathbf{x}_c),
\end{equation}
where $\mathbf{x} \in c$ and $\mathbf{x}_c$ is the barycentre of the cell $c$. The second is a $P^1$ reconstruction on boundary facets, written $\mathfrak{R}^\partial(v_{\mathcal{Z}^\partial})$, and computed using the vertex dofs of the boundary facets. In case of simplicial facets, it reduces to classical $P^1$ Lagrange interpolation. For non-simplicial facets, generalised barycentric coordibnates need to be used, see \cite{budninskiy2016power}, for instance.

\subsection{Discrete problem}

Let us set $V_h:=\mathbb{R}^{d\#(\mathcal{C})}\times \mathbb{R}^{d\#(\mathcal{Z}^\partial)}$
and (recall that $\partial\Omega_D$ is a closed set)
\begin{equation} \left\{ \begin{aligned}
V_{hD} &:=\{v_h\in V_h
\ | \ v_z=u_D(z) \ \forall z\in \mathcal{Z}^\partial\cap\partial\Omega_D\},
\\
V_{h0} &:=\{v_h\in V_h
\ | \ v_z=0 \ \forall z\in \mathcal{Z}^\partial\cap\partial\Omega_D\}, \\
W_{hD} &:= \{v_h\in V_h
\ | \ v_z=\dot{u}_D(z) \ \forall z\in \mathcal{Z}^\partial\cap\partial\Omega_D\}.
\end{aligned}\right. \end{equation}
(Note that the spaces $V_{hD}$ and $W_{hD}$ are actually time-dependent if 
the Dirichlet data is time-dependent.) 
The discrete stiffness bilinear form is parameterized by a member
$\eta_{p,\mathcal{C}}\in Q_h$ and is such that,
for all $(v_h,\tilde{v}_h)\in V_{hD}\times V_{h0}$
(compare with (\ref{eq:bilinear form})),
\begin{equation}
\label{eq: discrete bilinear form}
a_h(\eta_{p,\mathcal{C}};v_h,\tilde{v}_h) :=
\sum_{c\in\mathcal{C}}|c|\big(\mathbb{C} : (\varepsilon_{c}(\mathcal{R}(v_h))-\eta_{p,c})\big)
: \varepsilon_{c}(\mathcal{R}(\tilde{v}_h)) + s_h(v_h,\tilde{v}_h).
\end{equation}
Here $s_h$  is a weakly consistent stabilization bilinear form intended to
render $a_h$ coercive and which is defined on $V_h\times V_h$ as follows:
\begin{equation}
\label{eq:penalty term}
s_h(v_h,\tilde{v}_h) = \sum_{F \in \mathcal{F}}  \frac{\eta}{h_F} |F| [\mathfrak{R}(v_h)]_F \cdot [ \mathfrak{R}(\tilde{v}_h)]_F 
\end{equation} 
where $h_F$ is the diameter of the facet $F\in\mathcal{F}$.
For an interior facet $F\in\mathcal{F}^i$, writing $c_-$
and $c_+$ the two mesh cells sharing $F$, i.e., $F=\partial c_-\cap \partial c_+$, 
and orienting $F$ by the unit normal
vector $n_F$ pointing from $c_-$ to $c_+$, one has 
\begin{equation}
[\mathfrak{R}(v_h)]_F :=
\mathfrak{R}(v_h)_{c_-}(\mathbf{x}_F) - \mathfrak{R}(v_h)_{c_+}(\mathbf{x}_F).
\end{equation}
The sign of the jump is irrelevant in what follows. The role of the 
summation over the interior facets in~\eqref{eq:penalty term}
is to penalize the jumps of the cell reconstruction $\mathfrak{R}$ across the interior facets. 
For a boundary facet $F\in\mathcal{F}^b$, we denote $c_-$ the unique mesh cell 
containing $F$, we orient $F$ by the unit normal vector $n_F:=n_{c_-}$ which
points outward $\Omega$, and we define 
\begin{equation}
[\mathfrak{R}(v_h)]_F := \mathfrak{R}^\partial(v_{\mathcal{Z}^\partial})_F(\mathbf{x}_F)-\mathfrak{R}(v_h)_{c_-}(\mathbf{x}_F).
\end{equation}
The role of the summation over the boundary facets 
in~\eqref{eq:penalty term} is to
penalize the jumps between the cell reconstruction $\mathfrak{R}$ and the boundary facet reconstruction $\mathfrak{R}^\partial$.
The parameter $\eta > 0$ in~\eqref{eq:penalty term} is user-defined with the only
requirement that $\eta > 0$.
The bilinear form $s_h$ is classical in the context of discontinuous Galerkin methods (see \cite{arnold1982interior,ern_discontinuous} for instance, see also \cite{di2012cell} for cell-centred Galerkin methods).
In practice, the penalty parameter $\eta$ scales as $\eta = \beta \mu$ where $\mu$ is the second Lamé coefficient of the material and $\beta$ is a dimensionless factor that remains user-dependent. Notice that this choice is robust with respect to $\nu \to 0.5$.
We present a numerical test in Section \ref{sec:choice penalty} illustrating the moderate impact of $\beta$ on the numerical computations.

We can next define a discrete mass bilinear form $m_h$ 
similar to (\ref{eq:kinetic power}) and a discrete
load linear form $l_h(t)$ similar to (\ref{eq:linear form}); details are given below.
Then the space semi-discrete version of the evolution 
problem~\eqref{eq:elasto-plasticity-weak} amounts to seeking
$(u_h,\varepsilon_{p,\mathcal{C}},p_{\mathcal{C}}):(0,T)\rightarrow V_{hD}\times Q_h\times P_h$ such that,
for all $t\in (0,T)$, 
\begin{equation} \label{eq:discretization} \left\{ 
\begin{alignedat}{2} 
& m_h(\ddot{u}_h(t),\tilde{v}_h) +
a_h(\varepsilon_{p,\mathcal{C}}(t);u_h(t),\tilde{v}_h) =
l_h(t,\tilde{v}_h),&\quad&\forall \tilde{v}_h\in V_{h0}, \\
&\lambda_{c}\ge 0, \quad \varphi(\Sigma_{c}(u_h),p_c)\leq 0,
\quad  \lambda_c \varphi(\Sigma_{c}(u_h),p_c)= 0,
&\quad&\forall c\in\mathcal{C},\\
&\dot{p}_c=\lambda_{c}, \quad \dot{\varepsilon}_{p,c} =
\lambda_c\frac{\partial\varphi}{\partial\sigma}(\Sigma_{c}(u_h)),
&\quad&\forall c\in\mathcal{C},\\
\end{alignedat}\right.\end{equation}
where the time-dependency is left implicit in the second and third equations
and where we introduced in every mesh cell $c\in\mathcal{C}$ the local
stress tensor 
\begin{equation} \label{eq:stress_uh}
\Sigma_{c}(u_h) := \mathbb{C}:(\varepsilon_{c}(\mathcal{R}(u_h)) - \vespc) 
\in \mathbb{R}^{d\times d}_{\mathrm{sym}}.
\end{equation}
Note that the plasticity relations in~\eqref{eq:discretization}
are enforced cellwise, i.e., a mesh cell $c\in\mathcal{C}$ is either
in the elastic state or in the plastic state depending on the value of 
$\varphi(\Sigma_{c}(u_h),p_c)$. This is due to the fact that stresses are cell-wise constant and thus the second relation of \eqref{eq:discretization} can only be enforced cell-wise.

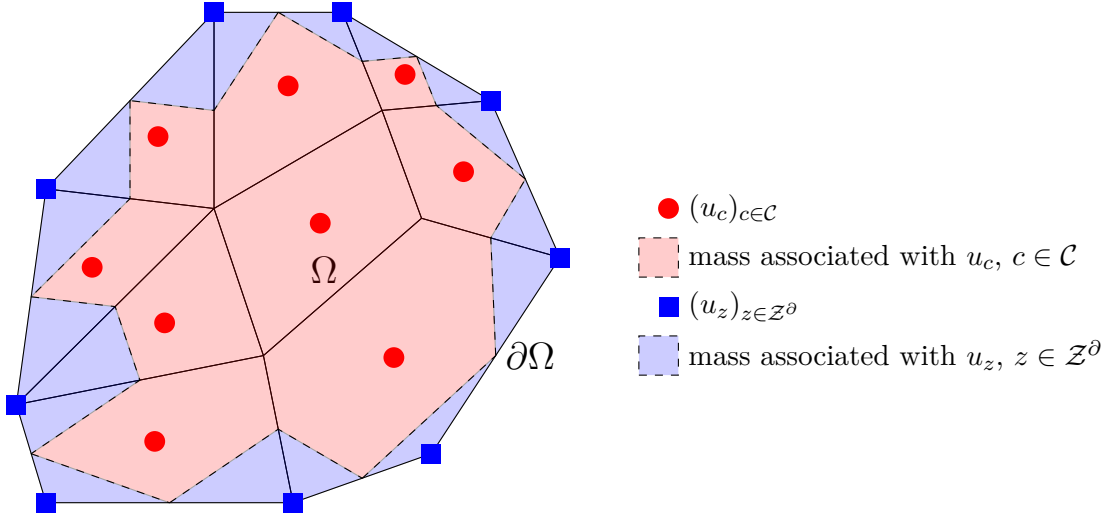
\begin{figure} [!htp] 
\begin{center}
\begin{tikzpicture}[scale=1.3]
\coordinate (a) at (0,-1);
\coordinate (b) at (0,1); 
\coordinate (c) at (1.3,1); 
\coordinate (d) at (-1.7,-0.8); 
\coordinate (e) at (-2,-3); 
\coordinate (f) at (0.5,-2.5); 
\coordinate (g) at (1.7,0);
\coordinate (h) at (2.8,0.1); 
\coordinate (i) at (2.1,-1.1);
\coordinate (j) at (3.5,-1.5); 
\coordinate (k) at (2.2,-3.5); 
\coordinate (l) at (0.8,-4); 
\coordinate (m) at (-1.7,-4); 

\node[right] at (barycentric cs:f=0.8,a=0.25,g=0.25,i=0.25) {\Large
  $\Omega$};
\node[right] at (barycentric cs:k=0.5,j=0.5) {\Large $\partial\Omega$};
\path[draw] (a)-- (b)-- (c) -- (g) -- cycle;
\path[draw] (a)-- (b)-- (d) -- cycle;
\path[draw] (a) -- (e) -- (d) -- cycle;
\path[draw] (a) -- (f) -- (e) -- cycle;
\path[draw] (a) -- (g) -- (i) -- (f) -- cycle;
\path[draw] (c) -- (h) -- (g) -- cycle;
\path[draw] (g) -- (h) -- (j) -- (i) -- cycle;
\path[draw] (i) -- (j) -- (k) -- (l) -- (f) -- cycle;
\path[draw] (f) -- (l) -- (m) -- (e) -- cycle;

\coordinate (bc) at (barycentric cs:b=0.5,c=0.5);
\coordinate (ab) at (barycentric cs:a=0.5,b=0.5);
\coordinate (bd) at (barycentric cs:b=0.5,d=0.5);
\draw[fill=blue,opacity=0.2] (b) -- (bc) -- (ab) -- (bd) --cycle;
\path[draw,dashed] (bc) -- (ab) -- (bd);
\coordinate (ad) at (barycentric cs:a=0.5,d=0.5);
\coordinate (de) at (barycentric cs:d=0.5,e=0.5);
\draw[fill=blue,opacity=0.2] (d) -- (bd) -- (ad) -- (de) --cycle;
\path[draw,dashed] (bd) -- (ad) -- (de);
\coordinate (ae) at (barycentric cs:a=0.5,e=0.5);
\coordinate (ef) at (barycentric cs:e=0.5,f=0.5);
\coordinate (em) at (barycentric cs:e=0.5,m=0.5);
\draw[fill=blue,opacity=0.2] (e) -- (de) -- (ae) -- (ef) -- (em) --cycle;
\path[draw,dashed] (de) -- (ae) -- (ef) -- (em);
\coordinate (lm) at (barycentric cs:l=0.5,m=0.5);
\draw[fill=blue,opacity=0.2] (m) -- (em) -- (lm) --cycle;
\path[draw,dashed] (em) -- (lm);
\coordinate (fl) at (barycentric cs:f=0.5,l=0.5);
\coordinate (kl) at (barycentric cs:k=0.5,l=0.5);
\draw[fill=blue,opacity=0.2] (l) -- (lm) -- (fl) -- (kl) --cycle;
\path[draw,dashed] (lm) -- (fl) -- (kl);
\coordinate (jk) at (barycentric cs:j=0.5,k=0.5);
\draw[fill=blue,opacity=0.2] (k) -- (kl) -- (jk) --cycle;
\path[draw,dashed] (kl) -- (jk);
\coordinate (ij) at (barycentric cs:i=0.5,j=0.5);
\coordinate (hj) at (barycentric cs:h=0.5,j=0.5);
\draw[fill=blue,opacity=0.2] (j) -- (jk) -- (ij) -- (hj) --cycle;
\path[draw,dashed] (jk) -- (ij) -- (hj);
\coordinate (gh) at (barycentric cs:g=0.5,h=0.5);
\coordinate (ch) at (barycentric cs:c=0.5,h=0.5);
\draw[fill=blue,opacity=0.2] (h) -- (hj) -- (gh) -- (ch) --cycle;
\path[draw,dashed] (hj) -- (gh) -- (ch);
\coordinate (cg) at (barycentric cs:c=0.5,g=0.5);
\draw[fill=blue,opacity=0.2] (c) -- (ch) -- (cg) -- (bc) --cycle;
\path[draw,dashed] (ch) -- (cg) -- (bc);
\draw[fill=red,opacity=0.2] (bc) -- (ab) -- (bd) -- (ad) -- (de) --
(ae) -- (ef) -- (em) -- (lm) -- (fl) -- (kl) -- (jk) -- (ij) -- (hj)  -- (gh) -- (ch) -- (cg) -- cycle;

\fill[red] (barycentric cs:a=0.3,b=0.3,d=0.3) circle (3pt);
\fill[red] (barycentric cs:a=0.3,e=0.3,d=0.3) circle (3pt);
\fill[red] (barycentric cs:f=0.3,e=0.3,a=0.3) circle (3pt);
\fill[red] (barycentric cs:a=0.25,g=0.25,c=0.25,b=0.25) circle (3pt);
\fill[red] (barycentric cs:g=0.3,h=0.3,c=0.3) circle (3pt);
\fill[red] (barycentric cs:a=0.25,g=0.25,i=0.25,f=0.25) circle (3pt);
\fill[red] (barycentric cs:h=0.25,g=0.25,i=0.25,j=0.25) circle (3pt);
\fill[red] (barycentric cs:i=0.2,j=0.2,k=0.2,l=0.2,f=0.2) circle
(3pt);
\fill[red] (barycentric cs:f=0.25,l=0.25,m=0.25,e=0.25) circle (3pt);

\node at (b) [fill=blue,minimum width=3pt,minimum height=3pt] {};
\node at (c) [fill=blue,minimum width=3pt,minimum height=3pt] {};
\node at (d) [fill=blue,minimum width=3pt,minimum height=3pt] {};
\node at (e) [fill=blue,minimum width=3pt,minimum height=3pt] {};
\node at (h) [fill=blue,minimum width=3pt,minimum height=3pt] {};
\node at (j) [fill=blue,minimum width=3pt,minimum height=3pt] {};
\node at (k) [fill=blue,minimum width=3pt,minimum height=3pt] {};
\node at (l) [fill=blue,minimum width=3pt,minimum height=3pt] {};
\node at (m) [fill=blue,minimum width=3pt,minimum height=3pt] {};

\fill[red] (4.6,-1) circle (3pt);
\node[right] at (4.7,-1) {$(u_c)_{c\in\mathcal{C}}$};
\draw[fill=red,fill opacity=0.2,dashed] (4.3,-1.7) -- (4.7,-1.7) --
(4.7,-1.3) -- (4.3,-1.3) -- cycle;
\node[right] at (4.7,-1.5) {mass associated with $u_c$, $c\in\mathcal{C}$};
\node at (4.6,-2) [fill=blue,minimum width=3pt,minimum height=3pt]{};
\node[right] at (4.7,-2) {$(u_z)_{z\in\mathcal{Z}^\partial}$};
\draw[fill=blue,fill opacity=0.2,dashed] (4.3,-2.7) -- (4.7,-2.7) --
(4.7,-2.3) -- (4.3,-2.3) -- cycle;
\node[right] at (4.7,-2.5) {mass associated with $u_z$, $z\in\mathcal{Z}^\partial$};
\end{tikzpicture}
\caption{Integration domains to determine the mass associated with the displacement dofs.}
\label{fig:mass}
\end{center}
\end{figure}

The definition of the discrete mass bilinear form $m_h$ 
hinges on subdomains to condense the mass 
associated with the dofs. Figure \ref{fig:mass} represents our
choice for the subdomains. For all the interior cells, the subdomain
$\omega_c$ is chosen as the whole cell, i.e., $\omega_c=c$. 
For the boundary vertices and for the
cells having a boundary face, a dual
barycentric subdomain is constructed, leading to subdomains denoted by
$\omega_z$ and $\omega_c$, respectively (see Figure \ref{fig:mass}).
For the discrete load linear form, we compute averages  of the external loads $f$ and $g_N$
in the mesh cells and on the Neumann boundary facets, respectively. 
As a consequence, $m_h(\cdot,\cdot)$ and $l_h(t;\cdot)$ can be written as follows for
all $(v_h,\tilde{v}_h)\in V_h\times V_h$ (compare with 
(\ref{eq:kinetic power}) and (\ref{eq:linear form})):
\begin{align}
  m_h(v_h,\tilde{v}_h) &:=
  \sum_{z\in\mathcal{Z}^\partial} m_z v_z\cdot \tilde{v}_z +
  \sum_{c\in\mathcal{C}} m_c v_c\cdot \tilde{v}_c,\\
  l_h(t,\tilde{w}_h) &:= \sum_{c\in\mathcal{C}}f_c(t) \cdot \tilde{w}_c + 
  \sum_{F\in\mathcal{F}^b\cap\partial\Omega_N} g_F(t) \cdot \mathcal{R}(\tilde{v}_h)_F,
\end{align}
with $m_z:=\int_{\omega_z}\rho$, $m_c:=\int_{\omega_c}\rho$,
$f_c(t):=\int_{c}f(t)$ and $g_F(t):=\int_F g_N(t)$.

\subsection{Reconstruction operator on facets}
\label{sec:reconstruction}

The reconstruction operator $\mathcal{R}$ is constructed in the same way as
in the finite volume methods studied in
\cite[Sec. 2.2]{eymard2009discretization} and in the cell-centered Galerkin methods from \cite{di2012cell}. For a
given facet $F\in\mathcal{F}$, we select neighbouring boundary vertices
collected in a subset denoted $\mathcal{Z}_F^\partial$ and neighbouring cells collected
in a subset denoted $\mathcal{C}_F$, as well as coefficients
$(\alpha_F^z)_{z\in\mathcal{Z}_F^\partial}$ and
$(\alpha_F^c)_{c\in\mathcal{C}_F}$, and we set
\begin{equation}\label{eq:barycentric}
  \mathcal{R}(v_h)_F :=
  \sum_{z\in\mathcal{Z}_F^\partial}{\alpha_F^zv_z}+
  \sum_{c\in\mathcal{C}_F}{\alpha_F^cv_c},
\qquad \forall v_h\in V_h.
\end{equation}
The neighbouring dofs should stay $\mathcal{O}(h)$ close
to the facet $F$. An algorithm is detailed thereafter to explain the selection of the neighbouring dofs in $\mathcal{Z}^{\partial}_F$ and $\mathcal{C}_F$.
The coefficients $\alpha_F^z$ and $\alpha_F^c$ are
chosen as the barycentric coordinates of the facet barycentre $\mathbf{x}_F$ in
terms of the location of the boundary vertices in $\mathcal{Z}_F^\partial$ and the barycentres
of the cells in $\mathcal{C}_F$.
For any interior facet $F \in \mathcal{F}^i$, the set $\mathcal{Z}_F^{\partial} \cup \mathcal{C}_F$ is constructed so as to contain exactly $(d+1)$ points forming the vertices of a non-degenerate simplex. Thus, the barycentric coefficients $\alpha^\mathbf{z}_F$ and $\alpha^c_F$ are computed by solving the linear system:
\begin{equation}
\label{eq:system barycentric coordinates}
\left\{
\begin{alignedat}{2}
&\sum_{\mathbf{z}\in\mathcal{Z}_F^\partial}{\alpha_F^\mathbf{z}}+
\sum_{c\in\mathcal{C}_F}{\alpha_F^c} = 1,\qquad &\forall F\in\mathcal{F}, \\
&\sum_{\mathbf{z}\in\mathcal{Z}_F^\partial}{\alpha_F^\mathbf{z} \mathbf{z}}+
\sum_{c\in\mathcal{C}_F}{\alpha_F^c \mathbf{x}_c} = \mathbf{x}_F,\qquad &\forall F\in\mathcal{F}, \\
\end{alignedat}
\right.
\end{equation}
where the vertex and its position are written $\mathbf{z}$ and $\mathbf{x}_c$ is the position of the barycentre of the cell $c$.

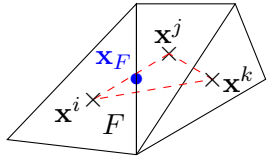
\begin{figure} [!htp]
\begin{center}
\begin{tikzpicture}
\path[draw] (0,-1)-- (0,1)-- (1.3,1) -- (0, -1);
\draw (0.43,0.34) node {$\times$};
\draw (0.43,0.34) node[above]{$\mathbf{x}^j$};

\draw (-0.3, -0.6) node{$F$};

\draw (0,0) node[blue] {$\bullet$};
\draw (-0.3,0.25) node[blue] {$\mathbf{x}_F$};
\path[draw, red, dashed] (0.43,0.34) -- (-0.57,-0.27)-- (1,0) -- (0.43,0.34);

\path[draw] (0,-1)-- (0,1)-- (-1.7,-0.8) -- (0, -1);
\draw (-0.57,-0.27) node {$\times$};
\draw (-0.9,-0.1) node[below]{$\mathbf{x}^i$};


\path[draw] (0,-1)-- (1.7,0)-- (1.3,1);
\draw (1,0) node {$\times$};
\draw (1,0) node[right]{$\mathbf{x}^k$};
\end{tikzpicture}
\caption{Dofs associated with the interior facet $F$ used for the reconstruction.}
\label{fig:tetra associe a une face}
\end{center}
\end{figure}

The main rationale for choosing the neighboring dofs is to ensure as much as possible
that all the coefficients
$\alpha_F^z$ or $\alpha_F^c$ lie in the interval $(0,1)$, so that 
the definition of $\mathcal{R}(v_h)_F$ in \eqref{eq:barycentric}
is based on an interpolation formula (rather than an extrapolation formula
if some coefficients lie outside the interval $(0,1)$.)
For most internal facets $F\in\mathcal{F}^i$, 
far from the boundary $\partial\Omega$, it is
possible to choose an interpolation-based 
reconstruction operator using only cell dofs, i.e., we usually have 
$\mathcal{Z}_F^\partial:=\emptyset$. 
Figure \ref{fig:tetra associe a une face} presents an example for 
an interior facet $F$ using three
neighbouring cell dofs located at the cell barycentres $\mathbf{x}^i$, $\mathbf{x}^j$ and $\mathbf{x}^k$. Close to
the boundary $\partial\Omega$, the use of boundary vertex dofs 
helps to prevent extrapolation. In all the cases we
considered, interpolation was always possible using the algorithm described below.

On the boundary facets, the reconstruction operator only uses interpolation
from the boundary vertices of the facet, i.e., we always set 
$\mathcal{C}_F:=\emptyset$ for all $F\in\mathcal{F}^b$. 
In three space dimensions, the facet can be
polygonal and the barycentric coordinates are generalized barycentric
coordinates. This is achieved using \cite{budninskiy2016power} and
the package 2D Triangulation of the geometric library CGAL. In
the case of simplicial facets (triangles in three space dimensions 
and segments in two space dimensions), we recover the classical barycentric interpolation as described above.

The advantage of using interpolation rather than extrapolation
is relevant in the context of explicit time-marching schemes where the 
time step is restricted by a CFL condition depending on the largest
eigenvalue $\lambda_{\max{}}$ of the stiffness matrix associated with the
discrete bilinear form $a_h(0;\cdot,\cdot)$ (see, e.g., \eqref{eq:CFL stability MEMM}). 
It turns out that using 
extrapolation can have an adverse effect on the maximal eigenvalue of the
stiffness matrix, thereby placing a severe restriction on the time step,
and this restriction is significantly alleviated if enough neighboring dofs are
used in~\eqref{eq:barycentric} to ensure that interpolation is being used. 
We refer the reader to Table~\ref{tab:interpolation VS extrapolation} below 
for an illustration.

Let us now briefly outline the algorithm used for selecting the 
reconstruction dofs associated with a given mesh inner facet $F\in\mathcal{F}^{i}$. 
This algorithm has to be viewed more as a proof-of-concept than as an optimized 
algorithm, and we observe that this algorithm is only used in a 
pre-processing stage of the computations. Fix an integer $I\ge d+1$. 
\begin{enumerate}
\item Compute a list of points $(\mathbf{x}_i)_{1\le i\le I}$ ordered by increasing distance to $\mathbf{x}_F$; each point can be either the barycentre of a mesh cell or a boundary vertex. To this purpose, 
we use the
KDTree structure of the \texttt{scipy.spatial} module of Python.
\item 
Check if $\mathbf{x}_F$ lies in the convex hull of the set $(\mathbf{x}_i)_{1\le i\le I}$. To this purpose we use the ConvexHull structure of \texttt{scipy.spatial}.
If that is not true, then extrapolation must be used. Otherwise find a subset of $(\mathbf{x}_i)_{1\le i\le I}$ containing $(d+1)$ points and use the barycentric coordinates of the resulting simplex to evaluate the interpolation coefficients to be used in~\eqref{eq:barycentric}. 

\end{enumerate}
Note that $I$ must be large enough to enable the recovery of at least one simplex per inner facet that is not too flat, independently of the use of extrapolation or interpolation. In our computations, we generally took $I=10$ if $d=2$ and $I=25$ if $d=3$.

To illustrate the performance of the above algorithm in alleviating time step restrictions
based on a CFL condition, we report in Table~\ref{tab:interpolation VS extrapolation} 
the largest eigenvalue $\lambda_{\max{}}$ of the stiffness matrix and the percentage of the 
mesh facets where extrapolation has to be used as a function of the integer parameter $I$ from the above algorithm. The results are obtained on the two
three-dimensional meshes (called "coarse" and "fine") considered in Section \ref{sec:dynamic flexion} together with the DEM discretization for the dynamic elasto-plastic evolution of a beam undergoing flexion. Note that the minimal value is $I=7$ for the coarse mesh and $I=9$ for the fine mesh.

\begin{table}[!htp]
\begin{center}
   \begin{tabular}{ | c | c | c | c | c | c | c | c | c | c | c |}
     \hline 
mesh & \multicolumn{2}{|c|}{$I=7$} & \multicolumn{2}{|c|}{$I=9$} & \multicolumn{2}{|c|}{$I=12$} & \multicolumn{2}{|c|}{$I=15$} & \multicolumn{2}{|c|}{$I=18$} \\ \hline
      coarse & $4{\cdot}10^{10}$ & $15\%$ & $6{\cdot}10^{09}$ & $4\%$ & $7{\cdot}10^{07}$ & $0.8\%$ & $2{\cdot}10^{05}$ & $0\%$ & $2{\cdot}10^{05}$ & $0\%$ \\ \hline
      fine & - & - & $7{\cdot}10^{09}$ & $1.6\%$ & $3{\cdot}10^{07}$ & $0.1\%$ & $1{\cdot}10^{07}$ & $0.01\%$ & $8{\cdot}10^{05}$ & $0\%$ \\ \hline
   \end{tabular}
   \caption{Largest eigenvalue of the stiffness matrix and percentage of inner facets with extrapolation for various values of the parameter $I$ on the coarse and fine meshes used in the DEM discretizations reported in Section \ref{sec:dynamic flexion}.}
   \label{tab:interpolation VS extrapolation}
\end{center}
\end{table}

\subsection{Interpretation as a Discrete Element Method}

In this section we rewrite the first equation in~\eqref{eq:discretization}
as a particle method by introducing the dofs of the discrete displacement $u_h(t)
\in V_{hD}$ attached to the mesh cells and to the boundary vertices lying on the Neumann
boundary, which we write 
$U_{\mathrm{DEM}}:=(U_p(t))_{p\in\mathcal{P}}$ 
with $\mathcal{P}:=\mathcal{C} \cup \mathcal{Z}_N^\partial$ and
$\mathcal{Z}_N^\partial:=\{z\in \mathcal{Z}^\partial \ | \ z\in\partial\Omega_N\}$.
Here $\mathcal{P}$ can be viewed as the indexing set for the set of particles. Note that $U_\mathrm{DEM}$ is a collection of dof values forming a column vector, whereas $u_h$ is a piecewise-constant function.
Recalling the definition of the discrete mass bilinear form, we set 
$m_p:= \int_{\omega_c}\rho$ if $p=c\in \mathcal{C}$ and $m_p:=\int_{\omega_\mathbf{z}}\rho$
if $p=\mathbf{z}\in \mathcal{Z}_N^\partial$. Concerning the external loads, we set
$F_{\mathrm{DEM}}(t):=(F_p(t))_{p\in\mathcal{P}}$ with $F_p(t):=f_c(t)=\int_c f(t)$ if
$p=c$ and $F_p(t):=\sum_{F\in\mathcal{F}_\mathbf{z}} \alpha_F^\mathbf{z} g_F(t) = \sum_{F\in\mathcal{F}_\mathbf{z}} \alpha_F^\mathbf{z}
\int_F g_N(t)$ if $p=\mathbf{z}$, where $\mathcal{F}_\mathbf{z}\subset \mathcal{F}^b$ is the collection of
the boundary facets to which $\mathbf{z}$ belongs and the coefficients $\alpha_F^\mathbf{z}$ are those used
in~\eqref{eq:barycentric} for the facet reconstruction. Since the Neumann boundary
$\partial\Omega_N$ is relatively open in $\partial\Omega$, all the facets in $\mathcal{F}_\mathbf{z}$
belong to $\partial\Omega_N$ if $\mathbf{z}\in \mathcal{Z}_N^\partial$.

Recall that $\vespC:(0,T)\to Q_h$ is the discrete tensor of remanent plastic strain.
Let us use the shorthand notation $\Sigma_c(t):= \Sigma_c(u_h(t))$ as defined in~\eqref{eq:stress_uh}, as well as $\Sigma_{\mathcal{C}}(t)
:= (\Sigma_c(t))_{c\in\mathcal{C}}$. For a piecewise-constant function defined on the
mesh cells, say $w_{\mathcal{C}}=(w_c)_{c\in\mathcal{C}}$, we define the mean-value
$\{w_{\mathcal{C}}\}_F :=\frac12(w_{c_-}+w_{c_+})$ for all $F=\partial c_-\cap \partial c_+ \in
\mathcal{F}^i$. Recall that the interior facet $F$ is oriented by the unit normal
vector $n_F$ pointing from $c_-$ to $c_+$ and that the jump across $F\in\mathcal{F}^i$ 
is defined such that
$[w_{\mathcal{C}}]_F :=w_{c_-}-w_{c_+}$. Recall also that for a boundary facet 
$F\in\mathcal{F}^b$, $c_-$ denotes the mesh cell to which $F$ belongs and that $n_F$
is the unit normal vector to $F$ pointing outward $\Omega$. 
Then a direct calculation shows that for all 
$\tilde{v}_h\in V_{h0}$, 
\begin{equation}\begin{aligned}
-a_h(\vespC(t);u_h(t),\tilde{v}_h) = {}& 
\sum_{F\in\mathcal{F}^i} |F| (\{\Sigma_{\mathcal{C}}(t)\}_F{\cdot} n_F)\cdot [\tilde{v}_{\mathcal{C}}]_F \\
&+ \sum_{F\in\mathcal{F}^i} |F| ([\Sigma_{\mathcal{C}}(t)]_F{\cdot} n_F)\cdot (\{\tilde{v}_{\mathcal{C}}\}_F-\mathcal{R}(\tilde{v}_h)_F) \\
&+ \sum_{F\in\mathcal{F}^b} |F| (\Sigma_{c_-}(t){\cdot} n_F)\cdot (\tilde{v}_{c_-}-\mathcal{R}(\tilde{v}_h)_F) \\ 
&-\sum_{F \in \mathcal{F}} \frac{\eta}{h_F} |F| [\mathfrak{R}(u_h(t))]_F \cdot [\mathfrak{R}(\tilde{v}_h)]_F. \\ 
\end{aligned}\end{equation}
To simplify some expressions, we are going to neglect the second term on the above 
right-hand side since this term is of higher-order (it is essentially the product of two
jumps).
Recall that, by definition, the reconstruction
operator $\mathcal{R}$ on a boundary facet $F\in\mathcal{F}^b$ 
makes use only of the vertex dofs of that facet. Then, letting
$(\tilde{V}_p)_{p\in\mathcal{P}}$ be the collection of the dofs of the discrete 
test function $\tilde{v}_h$, we infer that 
\begin{equation}
-a_h(\vespC(t);u_h(t),\tilde{v}_h) \simeq \sum_{p\in\mathcal{P}} \Phi_p^{\mathrm{ep}}(t)\cdot \tilde{V}_p + \Phi_p^{\mathrm{pen}}(t)\cdot \tilde{V}_p,
\end{equation}
where $\Phi_p^{\mathrm{ep}}(t)$ is the elasto-plastic force acting 
on the particle $p\in\mathcal{P}$ and $\Phi_p^{\mathrm{pen}}(t)$ is the force induced by the penalty and acting on the same particle.
For all $c\in\mathcal{C}$, we have $\Phi_c^{\mathrm{ep}}(t) := 
\sum_{F\in\mathcal{F}_c^{i}\cup \mathcal{F}_c^{N}}
\Phi_{c,F}^{\mathrm{ep}}(t)$ with $\mathcal{F}_c^{i} := \{F\in \mathcal{F}^i \ | \ 
F\subset \partial c\}$,
$\mathcal{F}_c^{N} := \{F\in \mathcal{F}^b \ | \ 
F\subset \partial c \cap \partial \Omega_N\}$, and
\begin{equation} \label{eq:elastic_fluxes}
\Phi_{c,F}^{\mathrm{ep}}(t) := \begin{cases}
\iota_{c,F} |F| \{\Sigma_{\mathcal{C}}(t)\}_F{\cdot} n_F 
&\text{if $F\in \mathcal{F}_c^{i}$}, \\
|F| \Sigma_{c_-}(t){\cdot} n_F 
&\text{if $F\in \mathcal{F}_c^{N}$},
\end{cases}
\end{equation}
with $\iota_{c,F}:=n_c\cdot n_F$,
and for all $z\in \mathcal{Z}_N^\partial$, we have
\begin{equation}
\Phi_z^{\mathrm{ep}}(t) := -\sum_{F\in\mathcal{F}_z} \alpha_F^z \Phi_{c_-,F}^{\mathrm{ep}}(t),
\end{equation}
with $\Phi^{\mathrm{ep}}_{c_-,F}$ defined in~\eqref{eq:elastic_fluxes} (recall that $c_-$ denotes the unique mesh cell 
containing the boundary facet $F\in \mathcal{F}^b$).
Note that the principle of action and reaction is encoded in the fact that
$\iota_{c_-,F}+\iota_{c_+,F}=0$ for all $F=\partial c_-\cap \partial c_+\in\mathcal{F}^i$.

\begin{remarque} [Physical forces]
The quantities of \eqref{eq:elastic_fluxes} are remarkable in the sense that they correspond to the physical force that one expects between two particles: the average of the stress in each particle multiplied by the surface shared by the two particles and contracted with the normal of the corresponding facet. Notice that this force only depends on the macroscopic material parameters and does not depend on any added microscopic parameter.
\end{remarque}

The penalty force is composed of two terms, that is, for all $p \in \mathcal{P}$, $\Phi_p^{\mathrm{pen}}(t) := \Phi_p^{\mathrm{pen,1}}(t) + \Phi_p^{\mathrm{pen,2}}(t) $. We define the first term for every cell $c\in\mathcal{C}$ and every facet $F \in \mathcal{F}$ such that $F \subset \partial c$ as
\begin{equation} \label{eq:penalty_fluxes 1}
\Phi_{c,F}^{\mathrm{pen,1}}(t) := - \iota_{c,F} \frac{\eta}{h_F} |F| [\mathfrak{R}(u_h(t))]_F,
\end{equation}
and we define it for every Neumann boundary vertex $z\in \mathcal{Z}_N^\partial$ as
\begin{equation}
\Phi_z^{\mathrm{pen,1}}(t) := -\sum_{F\in\mathcal{F}_z} \alpha_F^z \Phi_{c_-,F}^{\mathrm{pen,1}}(t).
\end{equation}
The second term is more intricate since it is a byproduct of the extended stencil of the method. 
The set of facets using the dof of the particle $p \in \mathcal{P}$ (whether cell or boundary vertex) is defined as $\mathfrak{F}_p$. This means that if $F \in \mathfrak{F}_p$, then one has either $p \in \mathcal{C}_F$ or $p \in \mathcal{Z}_F^\partial$.
Let us recall that for a facet $F \in \mathcal{F}$, the cell dofs used for the reconstruction $\mathcal{R}(u_h(t))_F$ are collected in $\mathcal{C}_F$ and the boundary vertex dofs in $\mathcal{Z}_F^\partial$. Then the second penalty term writes for all $p \in \mathcal{P}$ and all $F \in \mathfrak{F}_p$ as
\begin{equation} \label{eq:penalty_fluxes 2}
\Phi_{p,F}^{\mathrm{pen,2}}(t) := - \sum_{c\ale{,} F \supset \partial c} \sum_{F' \subset \partial c} \iota_{c,F'}\frac{\eta}{h_{F'}} |F'| [\mathfrak{R}(u_h(t))]_{F'} \frac{|F|}{|c|} n_{F,c} \cdot (\mathbf{x}_{F'} - \mathbf{x}_c) \alpha^p_{F},
\end{equation}
where $\alpha_F^p$ is the barycentric coordinate associated with the particle $p \in \mathcal{P}$ in the reconstruction over the facet $F$.
As a consequence, the total second penalty term writes for all $p \in \mathcal{P}$ as
\begin{equation}
\Phi_{p}^{\mathrm{pen,2}}(t) := \sum_{F \in \mathfrak{F}_p} \Phi_{p,F}^{\mathrm{pen,2}}(t)
\end{equation}
Finally, putting everything together, we infer that the
the first equation in~\eqref{eq:discretization} becomes
\begin{equation}
m_p \ddot U_p(t) \simeq \Phi_p^{\mathrm{ep}}(t) + \Phi_p^{\mathrm{pen}}(t) + F_p(t), \qquad \forall p\in\mathcal{P}.
\end{equation}

\begin{remarque}[Matrix formulation]
Let us briefly describe the matrix formulation of the space semi-discrete 
problem~\eqref{eq:discretization}  in the case of elastodynamics, i.e., 
without plasticity. For simplicity we focus on the pure Neumann problem.
A matrix $\mathbf{R}\in \mathbb{R}^{d\#(\mathcal{F}) \times d\#(\mathcal{P})}$ 
corresponding to the reconstruction operator
$\mathcal{R}$ is first constructed. Its entries are the barycentric
coefficients of the dofs used for the reconstruction on
the face associated with the given line of $\mathbf{R}$. 
The lines of $\mathbf{R}$ associated with boundary facets have, by construction,
non-zero entries only for boundary vertex dofs.
The linearized strain matrix
$\mathbf{E}\in \left(\mathbb{R}_\mathrm{sym}^{d \times d} \right)^{\#(\mathcal{C})} \times \mathbb{R}^{d\#(\mathcal{F})}$ 
is composed of the tensorial coefficients
$\frac{1}{2}\frac{|F|}{|c|}(\otimes n_F+ n_F \otimes)$ on the lines associated with the mesh cell $c\in\mathcal{C}$ and the columns associated with the facets $F\subset\partial c$. 
The linear
elasticity matrix $\mathbf{C}\in \left(\mathbb{R}_\mathrm{sym}^{d \times d} \right)^{\#(\mathcal{C}) \times \#(\mathcal{C})}$ 
can be written as the block-diagonal
matrix where each block corresponds to the double contraction with the fourth-order 
elastic tensor $\mathbb{C}$ and multiplication by $|c|$.
The jump matrix $\mathbf{J}\in \mathbb{R}^{d\#(\mathcal{F})\times d\#(\mathcal{P})}$ discretises the jumps $[\mathfrak{R}(u_h)]_F$ on a facet F. Its assembly is not detailed for simplicity. However, it is assembled using the connectivity matrix of adjacent cells, the gradient matrix (similar to $\mathbf{E}$ but non-symetrized and composed of the tensorial coefficients
$ \frac{|F|}{|c|} \otimes n_F$) and $\mathbf{R}$.
Denoting $\mathbf{D}\in \mathbb{R}^{d\#(\mathcal{F})\times d\#(\mathcal{F})}$ 
the diagonal matrix with entry $\frac{\eta}{h_F}|F|$ 
for the facet $F$, the stabilization matrix $\mathbf{S}$
corresponding to the bilinear form $s_h$ can be written
$\mathbf{S} := \mathbf{J}\tr\mathbf{D}\mathbf{J}\in
\mathbb{R}^{d\#(\mathcal{P})\times d\#(\mathcal{P})}$.
Finally, denoting
$\mathbf{K}:=\mathbf{R}\tr\mathbf{E}\tr\mathbf{C}\mathbf{E}\mathbf{R}+\mathbf{S}
\in \mathbb{R}^{d\#(\mathcal{P})\times d\#(\mathcal{P})}$
the stiffness matrix, the space semi-discrete system
(\ref{eq:discretization}) in the case of elastodynamics reduces to
$\mathbf{M} \ddot{U}_{\mathrm{DEM}}(t) + \mathbf{K} U_{\mathrm{DEM}}(t)= F_{\mathrm{DEM}}(t)$.
\end{remarque}

\subsection{Convergence tests for linear elasticity}
\label{sec:convergence tests}

The goal of this section is to briefly verify the correct implementation of the method 
in the case of static linear elasticity by comparing the numerical predictions using DEM
with some analytical solutions and reporting the orders of convergence on sequences of
uniformly refined meshes. The model problem thus consists of finding $u\in V_D$ such that
\begin{equation}
\int_{\Omega} \varepsilon(\nabla u) : \mathbb{C} : \varepsilon(\nabla \tilde{u}) = \int_{\Omega} f \cdot \tilde{u}, \qquad \forall \tilde{u}\in V_0.  
\end{equation}
The $L^2$-error is computed as $\Vert u - \mathfrak{R}(u_h) \Vert_{L^2(\Omega)}$.
The energy error is based on the reconstructed linearized strain of the discrete solution in each mesh cell and is computed as $\Vert u - u_h \Vert_{\mathrm{en}} := \frac12 a_h(0; u-u_h,u-u_h)$ where $a_h$ is defined in \eqref{eq: discrete bilinear form}. Note that this last term also contains the energy associated with the penalty bilinear form $s_h$.
The convergence rates are approximated as
\begin{equation}
\text{order} = d \log\left(\frac{e_1}{e_2}\right)  \left(\log\left( \frac{n_2}{n_1} \right) \right)^{-1},
\end{equation}
where $e_1,e_2$ denote the errors on the computations with mesh sizes $h_1,h_2$ and the number of dofs $n_1,n_2$.

\subsubsection{Choice of penalty factor}
\label{sec:choice penalty}
\begin{figure}[htp]
\centering
\subfloat{
\includegraphics[width=0.5\textwidth]{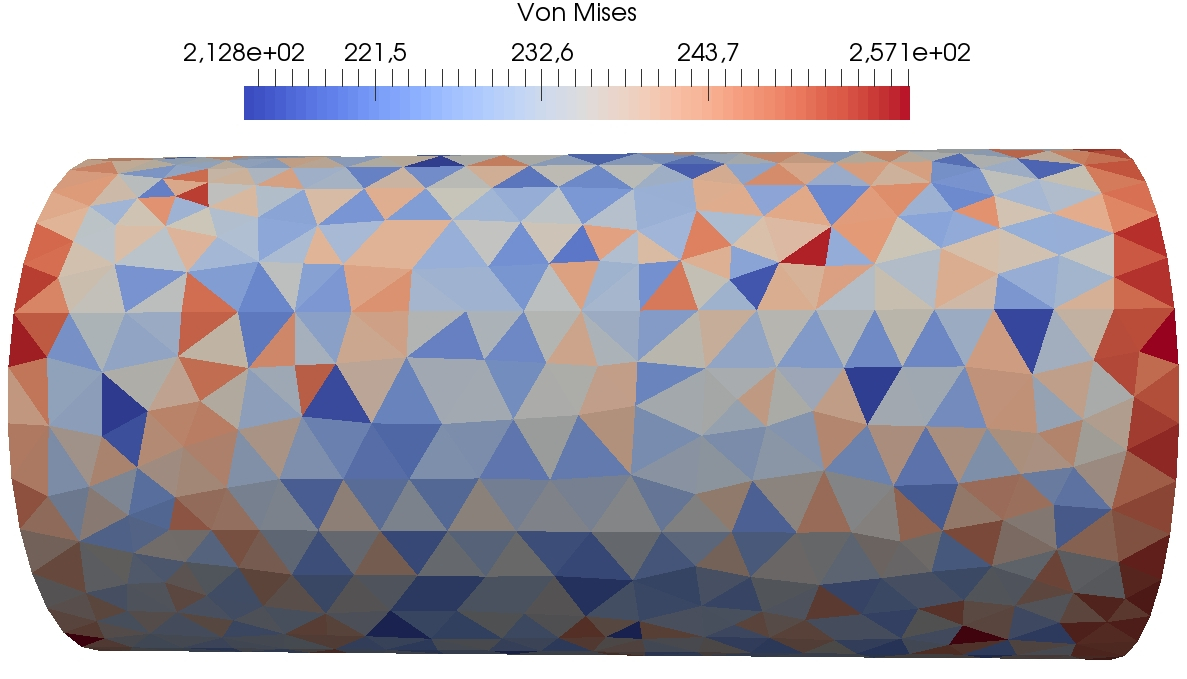}
}
\subfloat{
\includegraphics[width=0.5\textwidth]{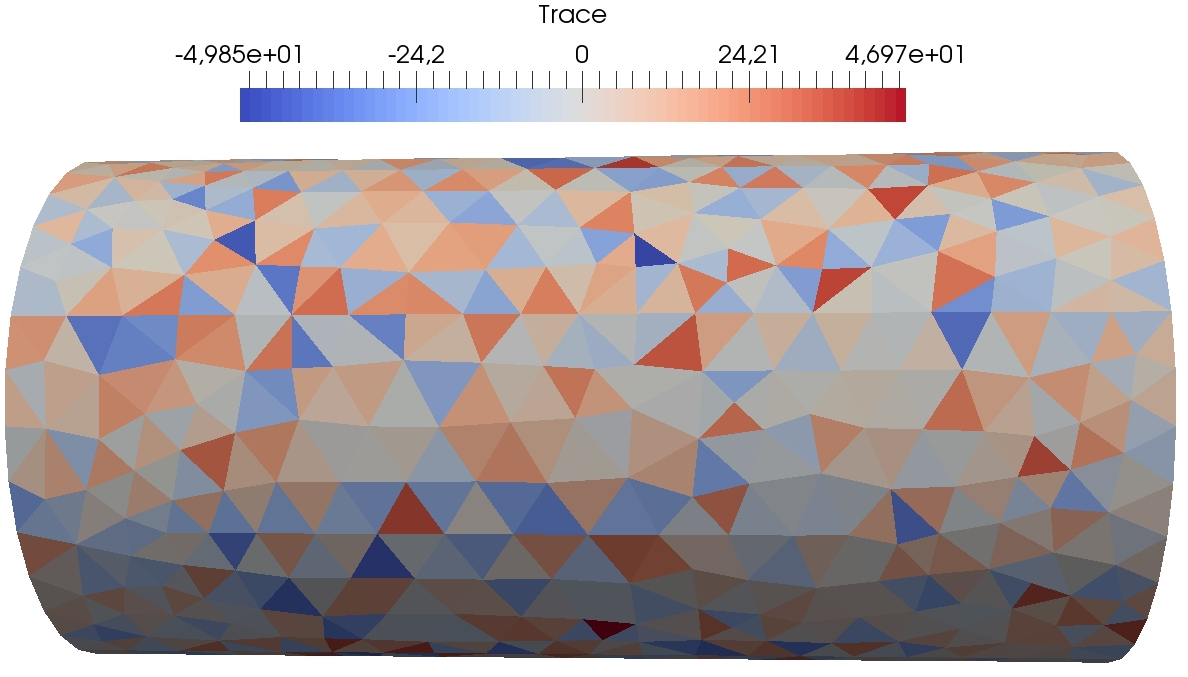}
} \\
\subfloat{
\includegraphics[width=0.5\textwidth]{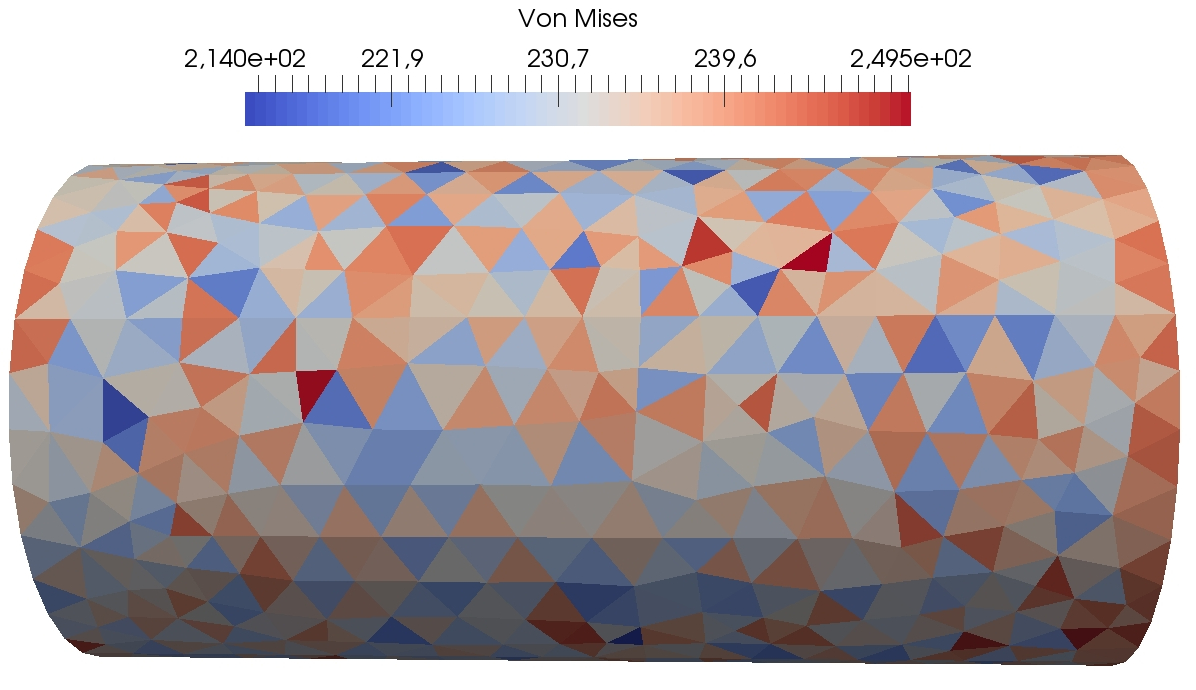}
}
\subfloat{
\includegraphics[width=0.5\textwidth]{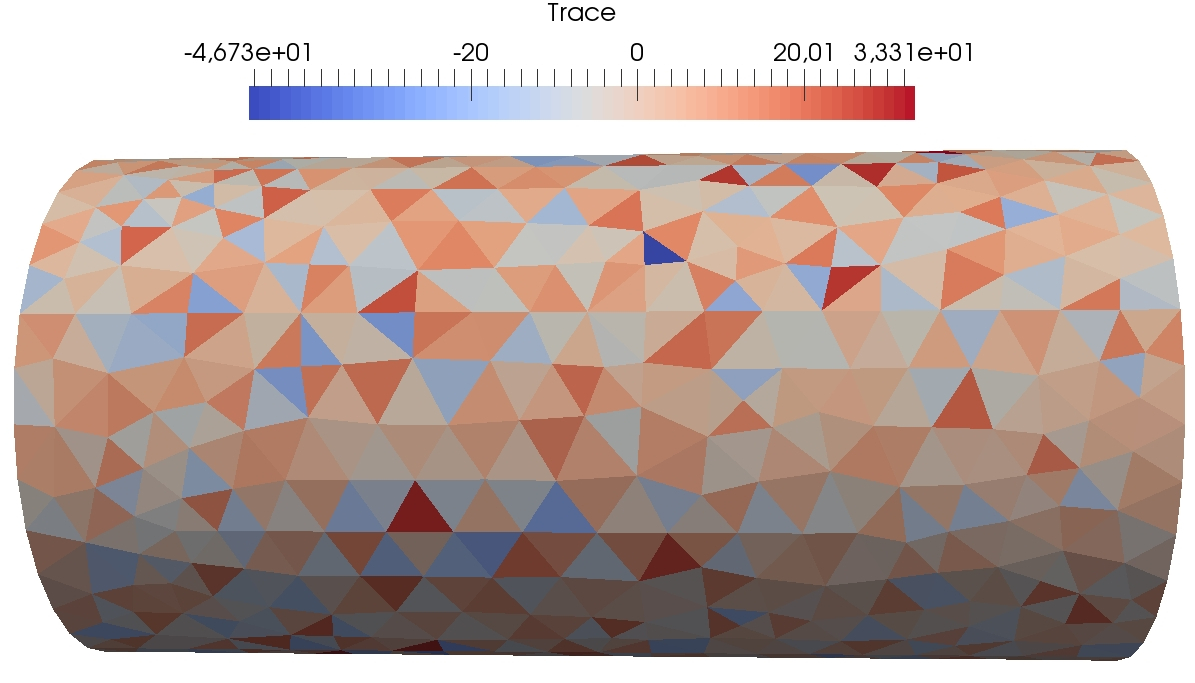}
} \\
\subfloat{
\includegraphics[width=0.5\textwidth]{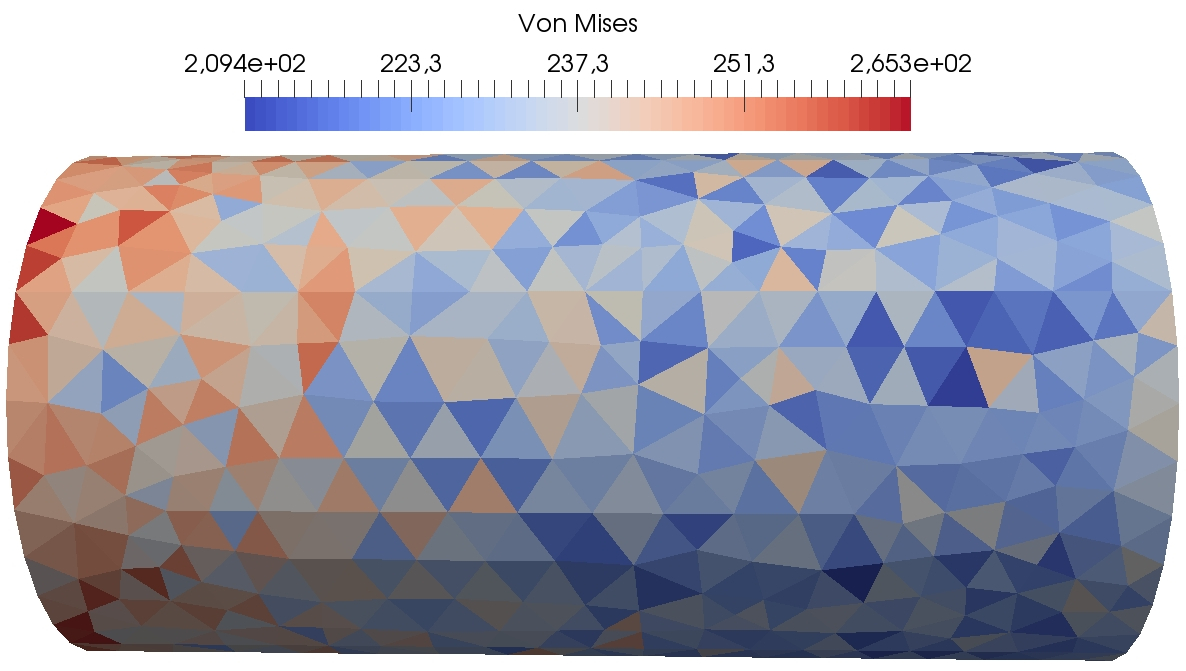}
}
\subfloat{
\includegraphics[width=0.5\textwidth]{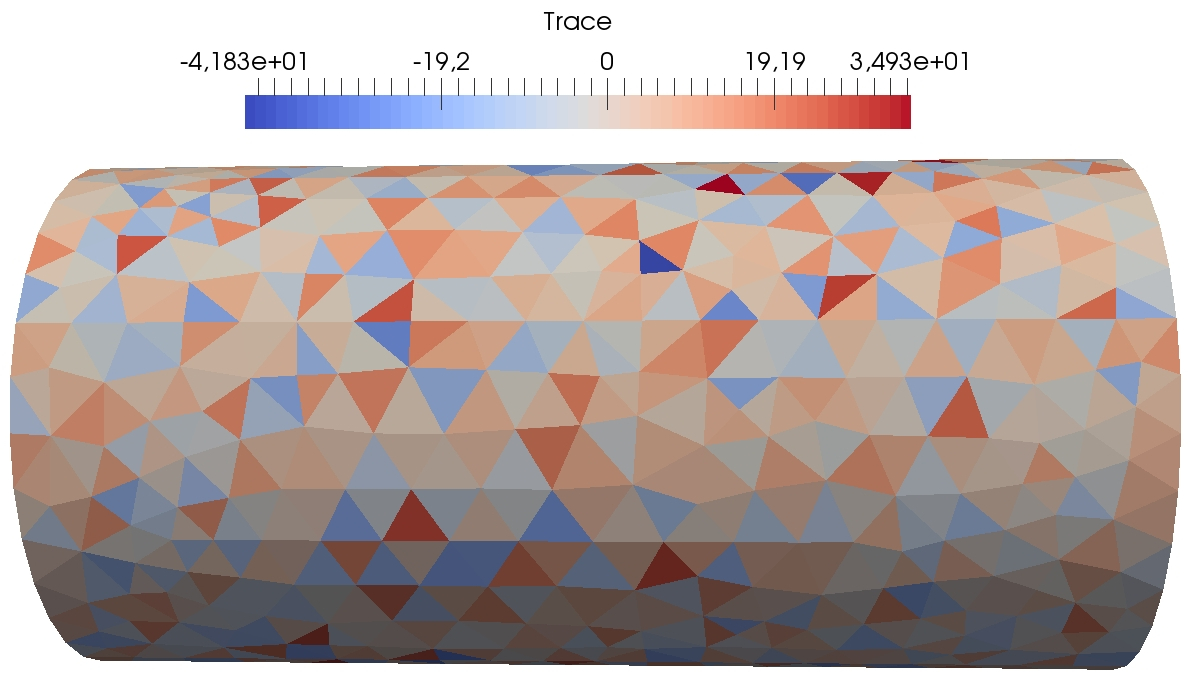}
} \\
\subfloat{
\includegraphics[width=0.5\textwidth]{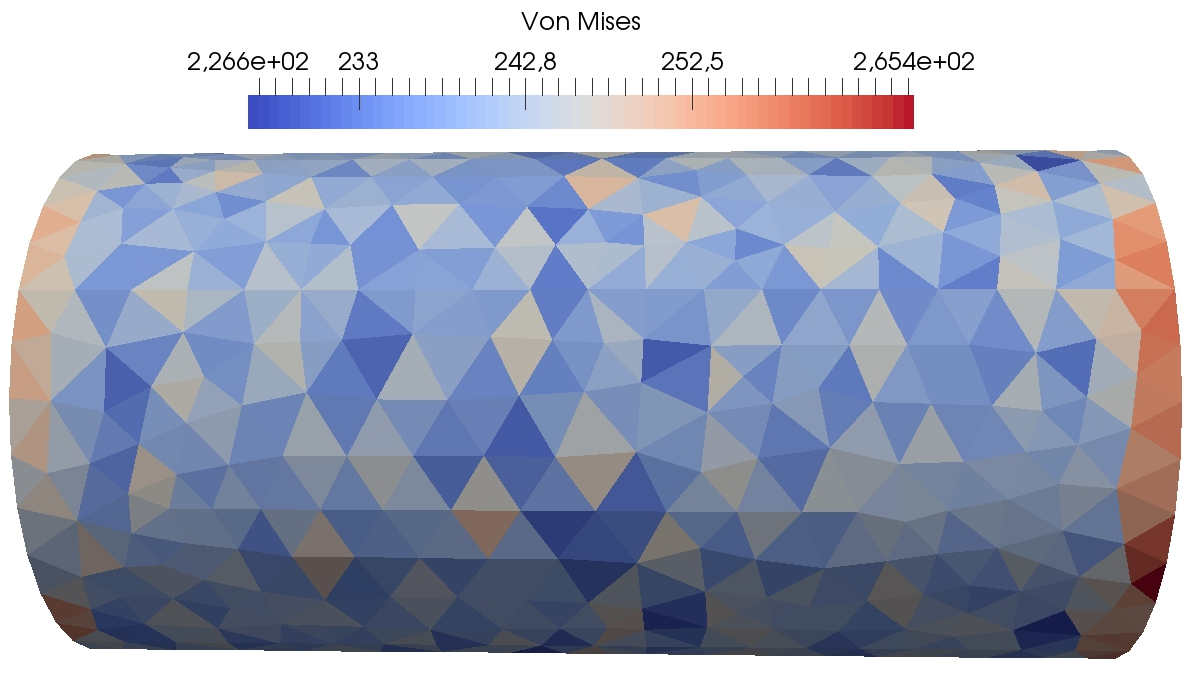}
}
\subfloat{
\includegraphics[width=0.5\textwidth]{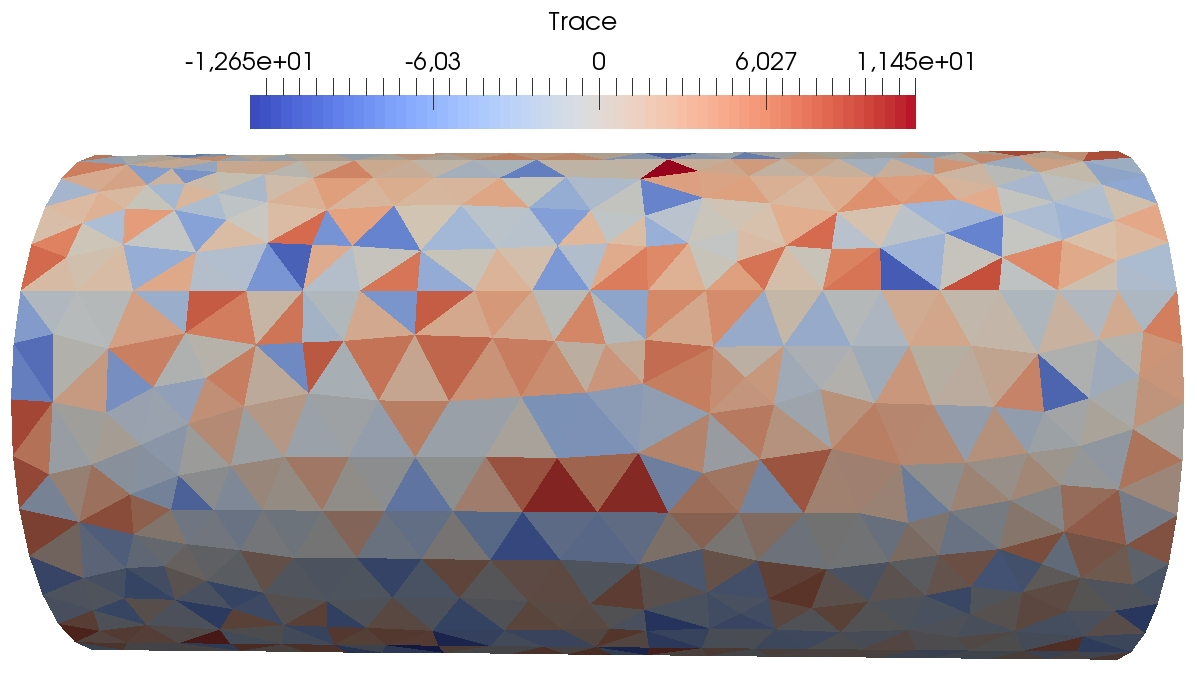}
} \\
\caption{Choice of penalty factor: Left: Von Mises equivalent stress. Right: Trace of the stress tensor. First row: DEM with $\beta = 1$. Second row: DEM with $\beta=0.1$. Third row: DEM with $\beta=0.01$. Fourth row: Penalised Crouzeix--Raviart FE with $\beta=1$ (reference).}
\label{fig:comparaison penalty parameters} 
\end{figure}

In this section, we simulate the torsion of a cylinder with various values of $\eta = \beta \mu$ obtained by varying $\beta$. The geometry is the one described in Section~\ref{sec:torsion quasi-static}.
We compare the results to first-order penalised Crouzeix--Raviart finite elements (FE) \cite{hansbo2003discontinuous}. A mesh of size $h= 29\mathrm{mm}$ is chosen for both computations. The DEM computation contains $6,332$ vectorial displacement dofs and the Crouzeix--Raviart $11,760$. The geometry and boundary conditions are similar to Figure \ref{fig:torsion plasticite parfaite}. The material is supposed to be isotropic elastic with $E=70\cdot 10^3\mathrm{Pa}$ and $\nu=0.3$.
The imposed torsion angle is $\alpha = 2.1 \cdot 10^{-2} \mathrm{rad}$. As the solution of this torsion test consists in pure shear stress, we compare in Figure \ref{fig:comparaison penalty parameters} $\mathrm{tr}(\sigma)$ and the Von Mises stress for DEM with $\beta=1$, $\beta=0.1$ and $\beta=0.01$,  and the reference penalised Crouzeix--Raviart computation with $\beta=1$. We expect the Von Mises stress to be constant on the lateral side of the cylinder and the trace of the stress tensor to be zero.

We observe that the influence of the parameter $\beta$ is rather light. Indeed, the trace of the trace tensor remains of order $10^1$ and the Von Mises stress varies between similar values. We observe that the results with the penalised Crouzeix--Raviart method vary slightly less. We believe that this is due to the higher number of dofs in this computation. As a consequence, so as to have a penalty term of the same order as the elastic term, we choose $\beta=1$ in our subsequent computations (unless explicitly mentioned).

\subsubsection{Manufactured solution}

Let us consider an isotropic two-dimensional elasticity test case in the domain
$\Omega=(-\frac12,\frac12)^2 \mathrm{m}^2$ with the Young modulus $E=70 \cdot 10^{3} \mathrm{Pa}$ and 
the Poisson ratio
$\nu=0.3$. The reference solution is $u(x,y) = \frac{a}{2} (x^2+y^2) (e_x + e_y)$ with $a=0.8\mathrm{m}^{-1}$
and $(e_x,e_y)$ is the Cartesian basis of $\Real^2$. The load term $f$, which is computed accordingly, is $f(x,y) = -a(\lambda + 3\mu) (e_x + e_y)$, where $\lambda$ and $\mu$ are the Lamé coefficients. The corresponding non-homogeneous Dirichlet boundary
condition is enforced on the whole boundary. Convergence results are reported 
in Table \ref{tab:energy errors 2} showing that the energy error converges to
first-order with the mesh size and the $L^2$-error to second-order.  
\begin{table}[!htp]
\begin{center}
   \begin{tabular}{ | c | c | c | c | c | c |}
     \hline 
     $h$ & dofs & $\Vert u - \mathfrak{R}(u_h) \Vert_{L^2(\Omega)}$ &  order & $\Vert u - u_h \Vert_{\mathrm{en}}$ & order\\ \hline
     0.0356 & $8,928$ & 5.67942e-05 & - & 4.72e-03 & - \\ \hline
     0.0177 & $35,072$ & 8.62031e-06 & 2.76 & 1.50e-03 & 1.60 \\ \hline
     0.00889 & $139,008$ &  1.80278e-06 & 2.27 & 5.87e-04 & 1.55 \\ \hline
   \end{tabular}
   \caption{Linear elasticity test case with manufactured solution: size of the mesh, number of dofs, $L^2$-error and order of convergence, energy error and order of convergence.}
   \label{tab:energy errors 2}
\end{center}
\end{table}

\section{Quasi-static elasto-plasticity}
\label{sec:quasi-static}

In this section we present the quasi-static elasto-plasticity problem,
its DEM discretization, and we perform numerical tests to assess the methodology.

\subsection{Governing equations}

The quasi-static elasto-plastic problem is a simplified formulation of~\eqref{eq:elasto-plasticity-weak}
where the inertia term in the mass bilinear form is negligible and where the time derivatives are
substituted by discrete increments. Thus we consider a sequence of loads $l^n\in V_0'$ for all
$n=1,\ldots,N$, and we consider the following sequence of nonlinear problems where $(u^n,\vesp^n,p^n) \in V_D
\times Q_0\times P$ for all $n=1,\ldots,N$:
\begin{equation}
\label{eq:elasto-plasticity-weak-qs}
\left\{ \begin{alignedat}{2}
&a(\vesp^n;u^n,\tilde{v}) = l^n(\tilde{v}),
&\quad&\forall \tilde{v}\in V_0,\\
&(\vesp^n,p^n,\mathbb{C}_{\mathrm{ep}}^n,\sigma^n) = \texttt{PLAS\_IMP}(\vesp^{n-1},p^{n-1},\veps^{n-1},\veps^n),
&\quad&\text{in $\Omega$},
\end{alignedat} \right.
\end{equation}
where $\veps^{n-1} := \veps(u^{n-1})$, $\veps^n:=\veps(u^n)$, and where variables with a superscript 
${}^{n-1}$ come from the solution of the quasi-static problem~\eqref{eq:elasto-plasticity-weak-qs} 
at the previous load increment or from a prescribed
initial condition if $n=1$. Given a quadruple $(\vesp\uold,p\uold,\veps\uold,\veps\unew)$, the procedure
\texttt{PLAS\_IMP} returns a quadruple $(\vesp\unew,p\unew,\mathbb{C}_{\mathrm{ep}}\unew,\sigma\unew)$ such that
\begin{equation}
\label{eq:PLAS_IMP}
\left\{ \begin{aligned}
&\lambda\unew \ge 0, \quad \varphi(\sigma\unew,p\unew) \leq 0, \quad
\lambda\unew\varphi(\sigma\unew,p\unew)=0,\\
&\delta p := p\unew-p\uold = \lambda\unew, \quad 
\delta\vesp:=\vesp\unew-\vesp\uold = \lambda\unew \frac{\partial \varphi}{\partial \sigma}(\sigma\unew).
\end{aligned} \right.
\end{equation}
Moreover $\sigma\unew=\mathbb{C} : (\veps\unew-\vesp\unew)$ is the new stress tensor, and $\mathbb{C}_{\mathrm{ep}}\unew$
is the consistent elastoplastic modulus \cite{simo1985consistent} such that
\begin{equation}
\sigma\unew = \sigma\uold +\mathbb{C}_{\mathrm{ep}}\unew : \delta \veps,
\quad  \delta \veps:=\veps\unew-\veps\uold,
\quad \sigma\uold := \mathbb{C} : (\veps\uold-\vesp\uold).
\end{equation}
The consistent elastoplastic modulus is instrumental to solve~\eqref{eq:elasto-plasticity-weak-qs}
iteratively by using an implicit radial return mapping technique (close to Newton--Raphson iterations) \cite{son1977elastic}:
Starting from $k=0$, we solve iteratively
the linear problem in $u^{n,k}\in V_D$ such that
\begin{equation}
\label{eq:elasto-plasticity-weak-qs-k}
\left\{ \begin{alignedat}{2}
&(\mathbb{C}_{\mathrm{ep}}^{n,k} : \veps(u^{n,k+1}),\veps(\tilde{v}))_Q = r^{n,k}(\tilde{v}) := l^n(\tilde{v}) - (\sigma^{n,k},\veps(\tilde{v}))_Q,
&\quad&\forall \tilde{v}\in V_0,\\
&(\vesp^{n,k},p^{n,k},\mathbb{C}_{\mathrm{ep}}^{n,k},\sigma^{n,k}) = \texttt{PLAS\_IMP}(\vesp^{n-1},p^{n-1},\veps^{n-1},\veps^{n,k}),
&\quad&\text{in $\Omega$},
\end{alignedat} \right.
\end{equation}
where the state for $k=0$ comes from the previous loading step or the initial condition.
Convergence of the iterative process in $k$ is reached when the norm of the residual $r^{n,k}$ is small enough. 

\subsection{DEM space discretization}

Using the DEM space discretization, the sequence of quasi-static problems~\eqref{eq:elasto-plasticity-weak-qs}
amounts to seeking a discrete triple $(u_h^n,\vespC^n,p_{\mathcal{C}}^n) \in V_{hD}\times Q_h\times P_h$ for all
$n=1,\ldots,N$, solving the following nonlinear problem:
\begin{equation}
\label{eq:elasto-plasticity-weak-qs-DEM}
\left\{ \begin{alignedat}{2}
&a_h(\vespC^n;u_h^n,\tilde{v}_h) = l_h^n(\tilde{v}_h),
&\quad&\forall \tilde{v}_h\in V_{h0},\\
&(\vespc^n,p_{c}^n,\mathbb{C}_{ep,c}^n,\sigma_c^n) = \texttt{PLAS\_IMP}(\vespc^{n-1},p_c^{n-1},\vespc^{n-1},\veps_c^n),
&\quad&\forall c\in\mathcal{C},
\end{alignedat} \right.
\end{equation} 
where $l_h^n$ represents a suitable discretization of the load linear form $l^n$. 
Using the radial return mapping technique as in~\eqref{eq:elasto-plasticity-weak-qs-k} and starting from $k=0$, we solve iteratively
the linear problem in $u_h^{n,k}\in V_{hD}$ such that
\begin{equation}
\label{eq:elasto-plasticity-weak-qs-k-DEM}
\left\{ \begin{alignedat}{2}
&\sum_{c\in\mathcal{C}} |c| (\mathbb{C}_{ep,c}^{n,k} : \veps_c(\mathcal{R}(u_h^{n,k+1}))) : \veps_c(\mathcal{R}(\tilde{v}_h)) 
= r_{\mathcal{C}}^{n,k}(\tilde{v}_h),
&\quad&\forall \tilde{v}_h\in V_{h0},\\
&(\vespc^{n,k},p_c^{n,k},\mathbb{C}_{ep,c}^{n,k},\sigma_c^{n,k}) = \texttt{PLAS\_IMP}(\vespc^{n-1},p_c^{n-1},\veps_c^{n-1},\veps_c^{n,k}),
&\quad&\forall c\in\mathcal{C},
\end{alignedat} \right.
\end{equation}
with the residual $r_{\mathcal{C}}^{n,k}(\tilde{v}_h) := l^n(\tilde{v}_h) - \sum_{c\in\mathcal{C}} |c| 
\sigma_c^{n,k} \cdot \veps_c(\mathcal{R}(\tilde{v}_h))$, and where the discrete state for $k=0$ 
comes from the 
previous loading step or by interpolating the values of the initial condition 
at the cell barycentres and the 
boundary vertices. Convergence of the iterative process in $k$ is reached when the 
norm of the residual $r_{\mathcal{C}}^{n,k}$ 
is small enough (we use a scaled Euclidean norm).

\subsection{Numerical tests}
\label{sec:quasi-static numerical tests}

This section contains two three-dimensional tests, a beam in quasi-static flexion and a beam in quasi-static torsion,
and a two-dimensional test case on the swelling of an infinite cylinder with internal pressure. 

\subsubsection{Beam in quasi-static traction}
\begin{figure} [!htp]
\begin{center}
\begin{tikzpicture}[scale=1]
\coordinate (a) at (-3,-0.5);
\coordinate (b) at (-3,0.5);
\coordinate (c) at (3,0.5);
\coordinate (d) at (3,-0.5);
\coordinate (ap) at (-2.5,0);
\coordinate (bp) at (-2.5,1.);
\coordinate (cp) at (3.5,1.);
\coordinate (dp) at (3.5,0);

\draw (a) -- (b) -- (c) -- (d) -- cycle;
\draw[dashed] (a) -- (ap);
\draw (b) -- (bp);
\draw[dashed] (bp) -- (ap);
\draw (c) -- (cp);
\draw (d) -- (dp);
\draw (cp) -- (dp);
\draw (bp) -- (cp);
\draw[dashed] (bp) -- (cp);
\draw[dashed] (ap) -- (dp);
\draw [<->] (-3,-1.2) -- (3,-1.2);
\draw (0,-0.9) node{$L$}; 
\draw (4.6,0.4) node {$u=u_D(t)$};
\draw (0,0.8) node{$\sigma \cdot n = 0$};
\draw(-3,0.2) node[left] {$u\cdot n=0$};
\end{tikzpicture}
\caption{Beam in quasi-static traction: problem setup.}
\label{fig:traction simple}
\end{center}
\end{figure}
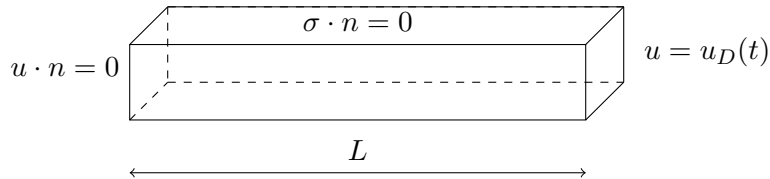

\begin{figure}[htp]
\begin{center}
\includegraphics[width=0.5\textwidth]{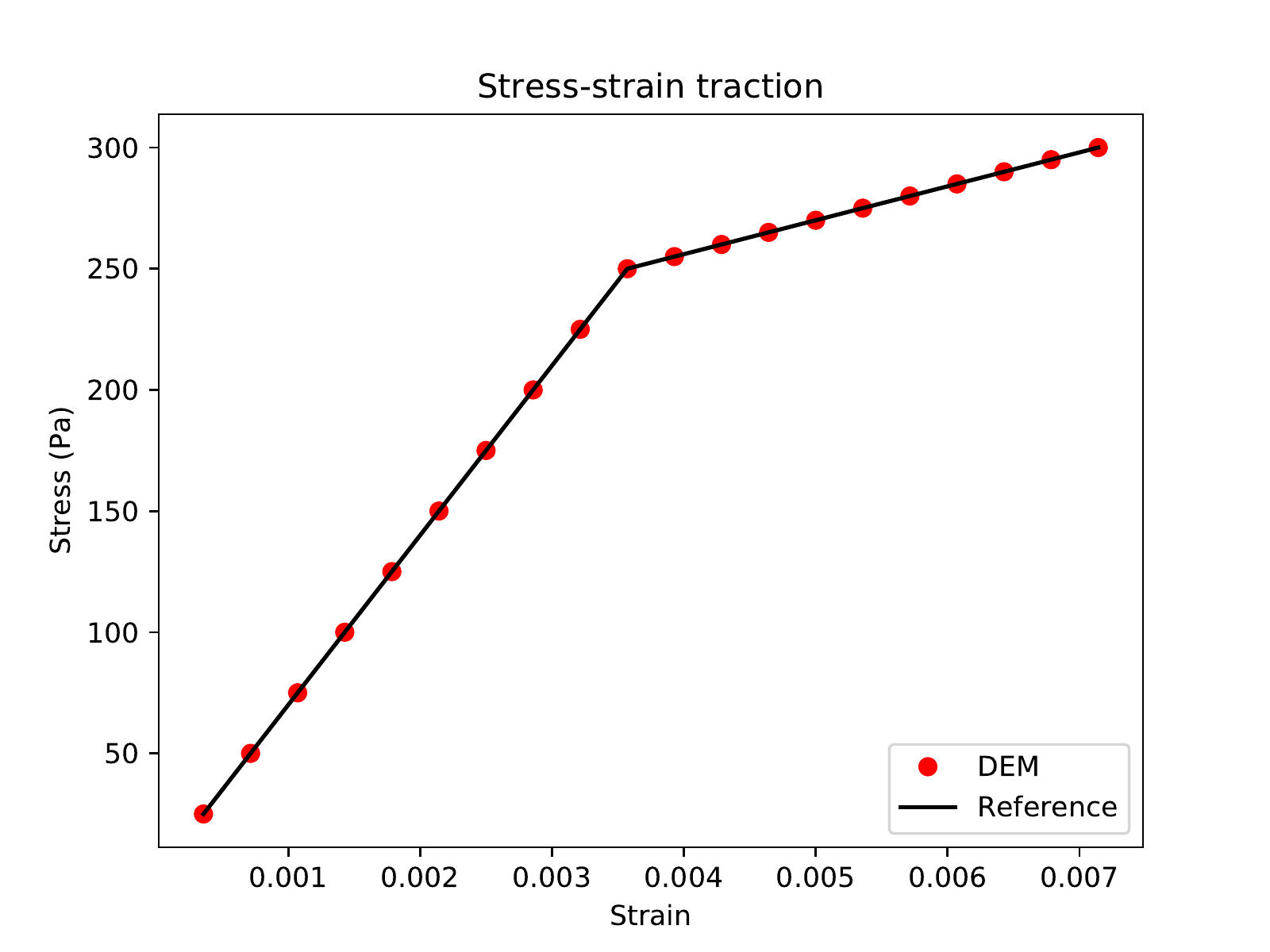}
\caption{Beam in quasi-static traction: stress-strain response curve for the analytical solution and the DEM solution.}
\label{fig:convergence traction}
\end{center}
\end{figure}

A beam of square section $S = 0.016\mathrm{m}^2$ and $L=1\mathrm{m}$ is stretched by a displacement $u_D(t)$ imposed at its right extremity, whereas the normal displacement and the tangential component of the normal stress are null at the left extremity. An homogeneous Neumann condition ($\sigma\cdot n=0$) is enforced on the four remaining sides of the beam.
Figure \ref{fig:traction simple} shows a sketch of the problem setup. The Young modulus is $E = 70\cdot 10^3 \mathrm{Pa}$ and the Poisson ratio $\nu = 0.3$. The yield stress is $\sigma_0 = 250 \mathrm{Pa}$, and the material is supposed to be elasto-plastic with linear kinematic hardening. Specifically the tangent plastic modulus is set to $E_t = \frac{1}{5}E$, so that we have $R(p)=Hp$ with $H=\frac{EE_t}{E-E_t}$. 
The imposed displacement is linearly increased in $20$ loading steps from $0$ to $2\delta_y$, where $\delta_y = \frac{\sigma_0}{E}L$ is the yield displacement. For this test case the analytical solution is available.
Figure \ref{fig:convergence traction} shows the stress-strain response curve, showing perfect agreement with the analytical solution using a mesh of size $h=0.2\mathrm{m}$. Note that in this test case, the stress tensor is constant in the beam.


\subsubsection{Beam in quasi-static torsion}
\label{sec:torsion quasi-static}
\begin{figure} [!htp]
\begin{center}
\begin{tikzpicture}
\pgfmathsetmacro{\R}{0.7}
\pgfmathsetmacro{\L}{5*\R}
\coordinate (O) at (0,0);
\coordinate (Op) at (\L/6,5*\L/6);
\coordinate (Ap) at (\L/6+\R,5*\L/6);
\coordinate (Bp) at (\L/6-\R,5*\L/6);

\draw (\R,0) arc (0:360:\R);

\draw (\R,0) -- (Ap);
\draw (-\R,0) -- (Bp);
\draw (0,\R) -- (\L/6,5*\L/6+\R);
\draw[dashed] (0,-\R) -- (\L/6,5*\L/6-\R);

\draw (Ap) arc (0:180:\R);
\draw[dashed] (Bp) arc (180:360:\R);

\draw (Ap) node[right] {$u=0$};
\draw (\L/12+\R/2+0.3,5*\L/12) node[right] {$\sigma\cdot n=0$};
\draw (\R+1.,0) node{$\alpha(t)$};
\draw [->] (\R + 0.3, -0.8) arc (-45:45:1);
\draw (0,-0.4-\R) node {$u = \alpha(t)re_{\theta}$};

\draw [<->] (-\R-0.5,0) -- (\L/6-\R-0.5,5*\L/6);
\draw (-1.2,\L/2) node{$L$};

\end{tikzpicture}
\caption{Beam in quasi-static torsion: problem setup.}
\label{fig:torsion plasticite parfaite}
\end{center}
\end{figure}
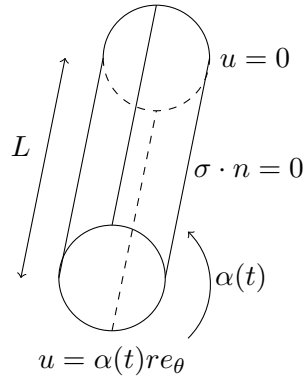

A beam of length $L = 0.2\mathrm{m}$ with a circular section of radius $R = 0.05\mathrm{m}$ is subjected to torsion at one of its extremities. The Young modulus is $E = 70\cdot 10^3 \mathrm{Pa}$ and the Poisson ratio $\nu = 0.3$. The yield stress is $\sigma_0 = 250 \mathrm{Pa}$, and the material is supposed to be perfectly plastic so that $R(p)=0$.
The beam is clamped at one of its extremities, a torsion angle $\alpha(t)$ is imposed at the other extremity, and the rest of the boundary of the beam is stress free ($\sigma\cdot n=0$). Figure \ref{fig:torsion plasticite parfaite} presents the problem setup. The torsion angle $\alpha(t)$ is increased linearly in 20 loading steps from $0$ to $\alpha_{\max{}} = 2 \alpha_{y}$, where $\alpha_y = \frac{\sigma_0 L}{\mu R \sqrt{3}}$ is the yield angle and $\mu$ is the second Lamé coefficient.  An analytical solution is available in the cylindrical frame $(e_r,e_{\theta},e_z)$: the displacement field is $u(r,z,t) = \alpha(t) r \frac{z}{L} e_{\theta}$, and the stress field is $\sigma(r,t) = \tau(r,t) \left(e_{\theta} \otimes e_z + e_z \otimes e_{\theta} \right)$, where
\[ \tau(r,t) := 
\left\{
\begin{alignedat}{2}
&\mu \alpha(t) \frac{r}{L}, \quad &  0 \leq r \leq R \frac{\alpha_y}{\alpha(t)}, \\
&\sigma_0, & R \frac{\alpha_y}{\alpha(t)} \leq r \leq R. \\
\end{alignedat}
\right.\] 

\begin{figure}[!htp]
\centering
\subfloat{
\includegraphics[width=0.5\textwidth]{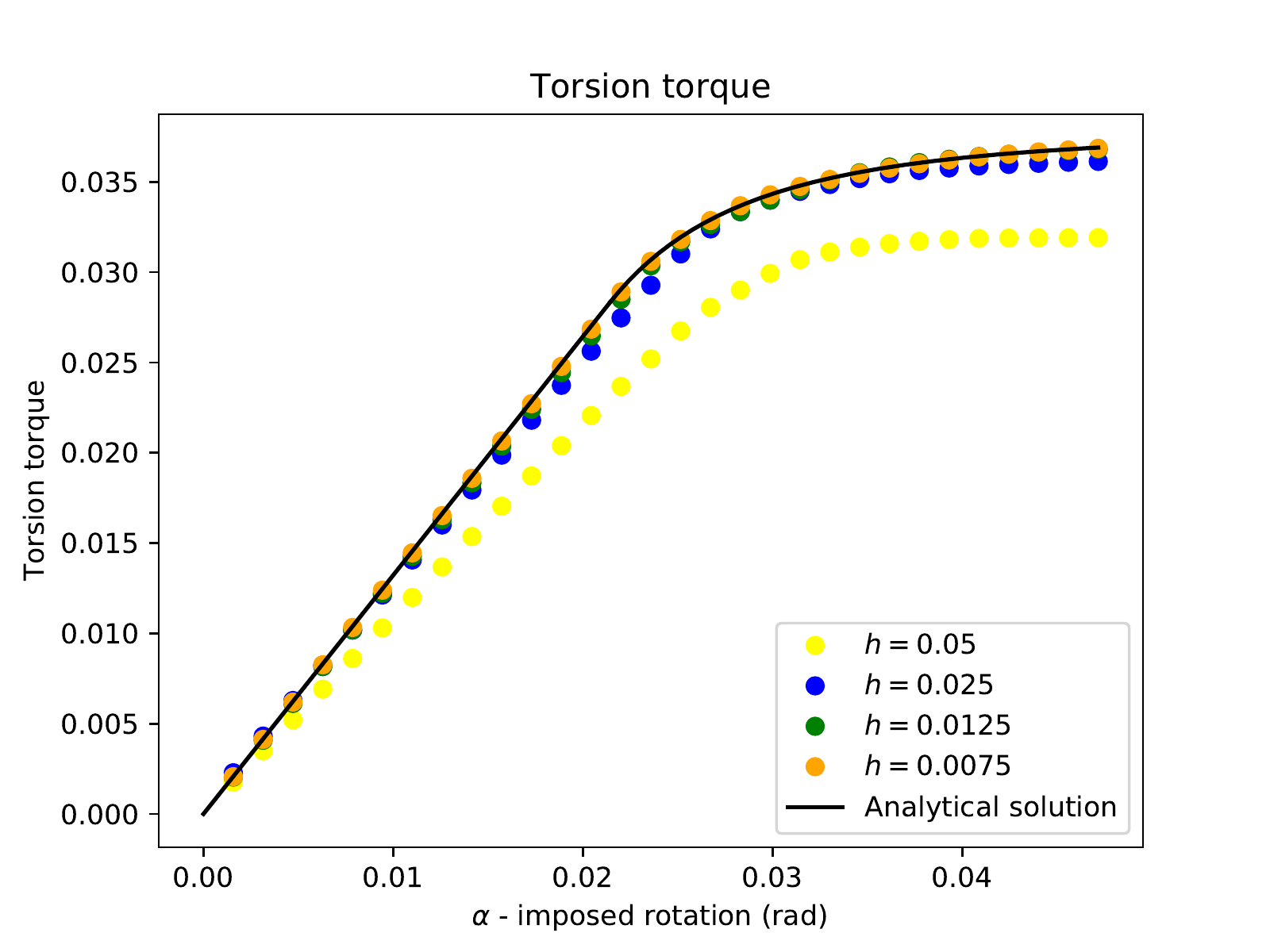}
}
\subfloat{
\includegraphics[width=0.5\textwidth]{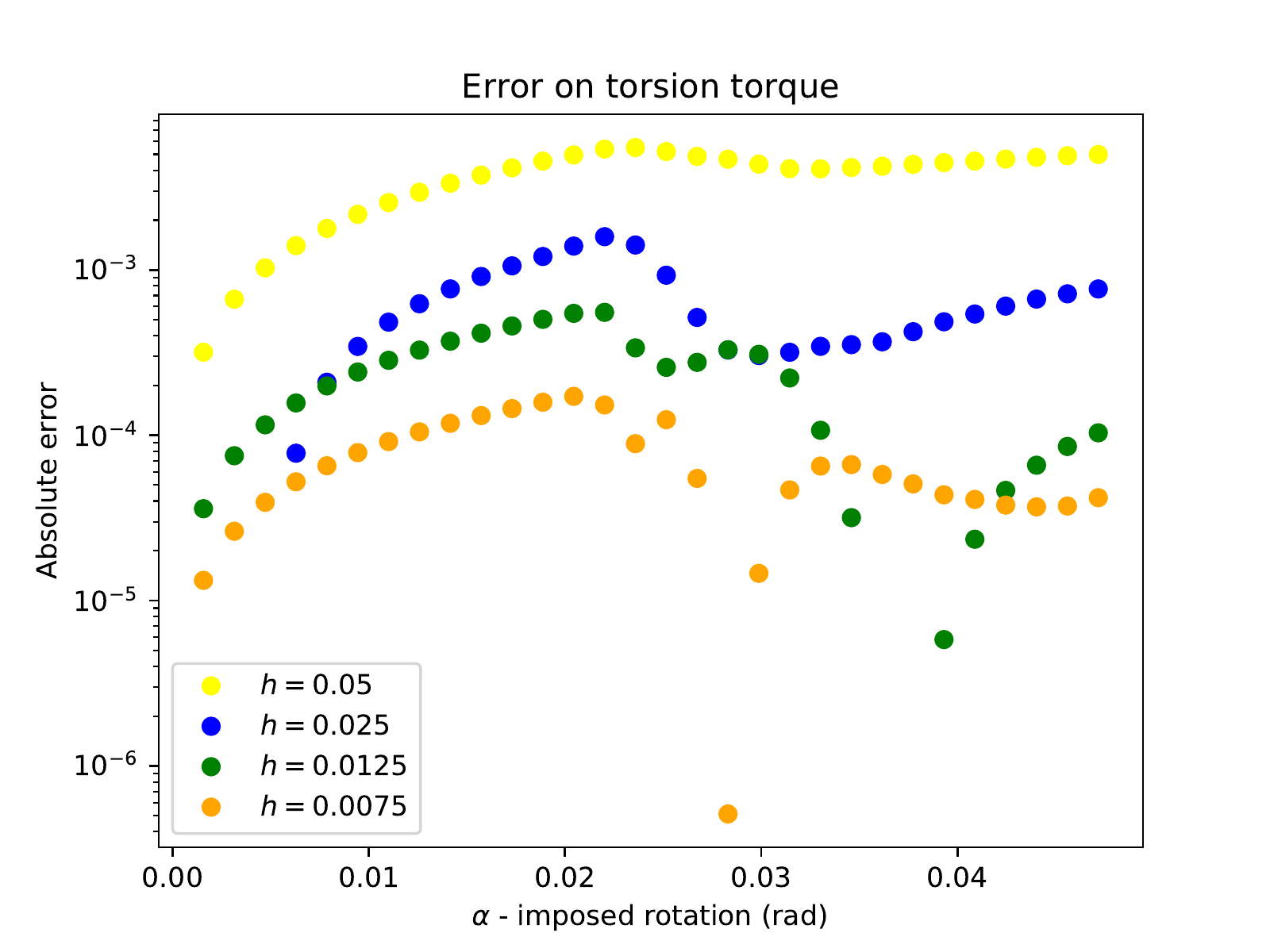}
}\\
\caption{Beam in quasi-static torsion. Left: torque-angle response curve for the analytical solution and the DEM solution on various meshes. Right: difference between the analytical solution and the DEM solution on various meshes.}
\label{fig:convergence torsion}
\end{figure}

Table \ref{tab:energy errors torsion} reports the maximum $L^2$-error on the displacement (evaluated as in Section~\ref{sec:convergence tests}) over the simulation interval and the energy error (including the energy of the penalty terms).
The errors are evaluated as described in Section~\ref{sec:convergence tests}. First-order convergence in the energy norm is observed, as expected. However, full second-order convergence in $L^2$ norm does not seem satisfied (although the convergence order is still above $1.77$). This can be due to the fact that in perfect plasticity $u \notin H^2(\Omega)$ which is a required hypothesis to obtain full second-order convergence.
\begin{table}[!htp]
\begin{center}
   \begin{tabular}{ | c | c | c | c | c | c |}
     \hline 
     $h$ & dofs & $\Vert u - \mathfrak{R}(u_h) \Vert_{L^2(\Omega)}$ &  order & $\Vert u - u_h \Vert_{\mathrm{en}}$ & order\\ \hline
     $0.05263$ & $3,753$ & 2.62e-06 & - & 5.87e-04 & -  \\ \hline 
     $0.03294$ & $12,726$ &  1.02e-06 & $2.32$ & 1.97e-04 & $2.68$\\ \hline 
     $0.02871$ & $18,996$ & 7.75e-07 & $2.04$  & 1.53e-04 & $1.89$ \\ \hline 
     $0.01965$ & $47,670$ & 4.50e-07 & $1.77$ & 8.23e-05 & $2.02$ \\ \hline 
     0.01418 & $160,146$ & 2.14e-07 & $1.84$ & 5.36e-05 & $1.06$ \\ \hline 

   \end{tabular}
   \caption{Beam in quasi-static torsion: size of the mesh, number of dofs, $L^2$-error and order of convergence, energy error and order of convergence.}
   \label{tab:energy errors torsion}
\end{center}
\end{table}
Figure \ref{fig:convergence torsion} presents the torque-angle response curve for the reference solution and the DEM solution on various meshes, showing good agreement and the convergence of the DEM predictions as the mesh is refined.

\subsubsection{Inner swelling of an infinite cylinder}

\begin{figure} [!htp]
\begin{center}
\begin{tikzpicture}
\draw (-3,-1) -- (-3,1);
\draw (2,-4) -- (0,-4);
\draw [-] (2, -4) arc (0:90:5);
\draw [-] (0, -4) arc (0:90:3);
\draw (-3,-4) node {$\times$};

\draw [<->] (-2.8,-4) -- (1.5,-2.2);
\draw (0.5,-2) node{$R_o = 1.3$};
\draw [<->] (-3,-3.8) -- (-2,-1.2);
\draw (-2.9,-1.4) node[below]{$R_i = 1$};

\draw (1.5,0) node {$\sigma\cdot n = 0$};
\draw (-4,0) node {$u \cdot n = 0$};
\draw (1,-4.5) node {$u \cdot n = 0$};
\draw (-1.5,-2.5) node {$ \sigma\cdot n = \varpi n$};

\end{tikzpicture}
\caption{Inner swelling of an infinite cylinder: problem setup.}
\label{fig:cylindre gonflé}
\end{center}
\end{figure}
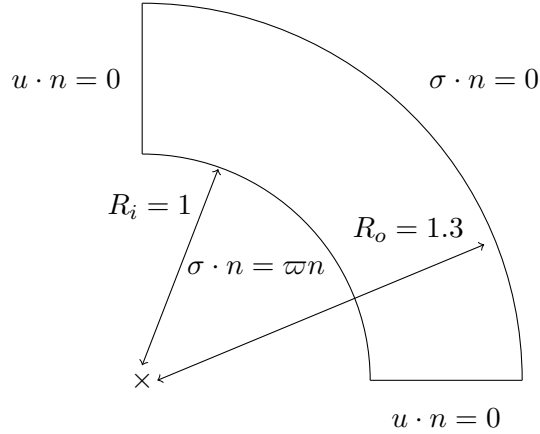

This test case consists in the inner swelling of an infinite cylinder. 
The inner radius is $R_i=1\mathrm{m}$ and the outer radius is $R_o=1.3\mathrm{m}$.
Owing to the symmetries, the computation is carried out on a quarter of a planar section of the cylinder with a plane strain formulation. A sketch of the problem setup is presented in Figure \ref{fig:cylindre gonflé}. On the lateral sides of the quarter of cylinder, a null normal displacement and a null tangential component of the normal stress are enforced. The outer side of the cylinder is stress free ($\sigma\cdot n=0$), and the inner pressure $\varpi$ imposed on the inner side is linearly increased from $0$ to $p_{\max{}} = \frac{2}{\sqrt{3}} \sigma_0 \mathrm{ln}\left(\frac{R_o}{R_i} \right)$, where $\sigma_0 = 250 \mathrm{N.m}^{-2}$ 
is the initial yield stress. The Young modulus and the tangent plastic modulus are set to $E = 70\cdot 10^3 \mathrm{Pa}$ and $E_t = \frac{1}{100}E$.

Table \ref{tab:energy errors 3} reports the $L^2$-error on the displacement (evaluated as in Section~\ref{sec:convergence tests}) and on the cumulated plastic strain. The reference solution is computed on the finest mesh and is based on $P^2$-Lagrange finite elements (FE) using the implementation available in \cite{bleyer2018numericaltours}.
The results in Table \ref{tab:energy errors 3} show that the method converges at order two in the $L^2$-norm and at order one in the energy norm.

\begin{table}[!htp]
\begin{center}
   \begin{tabular}{ | c | c | c | c | c | c |}
     \hline 
     $h$ & dofs & $\Vert u - \mathfrak{R}(u_h) \Vert_{L^2(\Omega)}$ &  order & $\Vert u - u_h \Vert_{\mathrm{en}}$ & order\\ \hline
     0.07735 & 992 & 5.15e-4 & - & 3.05e-03 & - \\ \hline
     0.04217 & 3412 & 1.69e-4 & $1.80$ & 8.38e-04 & $2.09$ \\ \hline
     0.02879 & 7588 &  7.26e-5 & $2.12$ & 3.76e-04 & $2.01$ \\ \hline
     0.02172 & 13380 & 3.85e-5 & $2.23$ & 2.10e-04 & $2.05$ \\ \hline
   \end{tabular}
   \caption{Inner swelling of an infinite cylinder: size of the mesh, number of dofs, $L^2$-error and order of convergence, energy error and order of convergence.}
   \label{tab:energy errors 3}
\end{center}
\end{table}

\section{Space-time fully discrete elasto-plasticity}
\label{sec:MEMM}

In this section we consider the dynamic elasto-plasticity equations from Section~\ref{sec:governing}.
The time discretization is performed by means of an explicit, pseudo-energy conserving, time-integration 
scheme recently introduced in \cite{MARAZZATO2019906}.
The space discretization is
achieved by means of the DEM scheme discussed in Sections~\ref{sec:dofs} and \ref{sec:quasi-static}. The main difference with Section \ref{sec:quasi-static numerical tests} is that no iterative procedure is necessary in this section.
Three-dimensional test cases are 
presented to assess the proposed methodology.

\subsection{Time semi-discretization of dynamic elasto-plasticity}

For simplicity we consider in this section only the time semi-discretization 
of the dynamic elasto-plasticity 
equations~\eqref{eq:elasto-plasticity-weak}. We deal with the space-time fully 
discrete setting in the next section.  
The time-integration scheme \cite{MARAZZATO2019906} is a two-step method of order two which ensures a 
discrete pseudo-energy conservation, if the integration of forces is exact, 
even for nonlinear systems. 
Symmetric Gaussian quadratures of the forces can be used in practice as long as they are of 
order at least two. 
The time interval $(0,T)$ is discretized using the time nodes $0 = t_0 < \ldots 
< t_n < \ldots < t_N = T$,
and for simplicity we consider a constant time step $\Delta t$. We define the half-time nodes 
$t_{n+\frac12} = \frac12(t_n + t_{n+1})$ for all $n=0,\ldots,N$. The time step is limited by a CFL
condition which we specify in the fully discrete setting in the next section. 

The key idea in the scheme \cite{MARAZZATO2019906} is to approximate the displacement field at the time nodes by means of functions
$u^n$, for all $n=0,\ldots,N$, with $u^0$ specified by the initial condition on the displacement, and the velocity
field at the half-time nodes by means of functions $v^{n+\frac12}$, for all $n=0,\ldots,N$, with $v^{\frac12}$ specified
by the initial condition on the velocity. For all $n=0,\ldots,N$, given $u^n$ and $v^{n+\frac12}$, the scheme 
performs two substeps: \textup{(i)} A time-dependent displacement field is predicted on the time interval
$[t_n,t_{n+1}]$ using the free-flight expression $\tilde u(t)=u^n+(t-t_n)v^{n+\frac12}$ for all $t\in [t_n,t_{n+1}]$;
\textup{(ii)} The velocity field $v^{n+\frac32}$ is predicted by means of a quadrature on the time-integration of the 
forces in the time interval $[t_n,t_{n+1}]$. Let $\{t_{n,k}\}_{k\in\mathcal{K}}$ and
$\{\omega_{n,k}\}_{k\in\mathcal{K}}$ be the nodes and the weights for the quadrature in the time interval
$[t_n,t_{n+1}]$. We then set 
\begin{equation} 
\left\{ \begin{alignedat}{2}
&u^{n,k} = u^n + (t_{n,k}-t_n)v^{n+\frac12},&\quad&\forall k\in\mathcal{K}, \\
&(\vesp^{n,k},p^{n,k}) = \texttt{PLAS\_EXP}(\vesp^{n,k-1},p^{n,k-1},\veps^{n,k}),&\quad&\forall k\in\mathcal{K}, \\
&\frac12 m(v^{n+\frac32}-v^{n-\frac12},\tilde{v}) = \sum_{k\in\mathcal{K}} \omega_{n,k} \Big( 
l(t_{n,k},\tilde{v}) - a(\vesp^{n,k};u^{n,k},\tilde{v})\Big),&\quad&\forall \tilde{v}\in V_0, \\
\end{alignedat} \right.
\end{equation}
where $\veps^{n,k}:=\veps(u^{n,k})$ is known from the free-flight displacement prediction
and where the state for the first Gauss node $k=1$ comes from the previous time step or the initial condition.
Given a triple $(\vesp\uold,p\uold,\veps\unew)$, the procedure \texttt{PLAS\_EXP} returns a pair
$(\vesp\unew,p\unew)$ such that, letting $\sigma\unew:=\mathbb{C} : (\veps\unew-\vesp\unew)$, we have
\begin{equation} \label{eq:PLAS_EXP}
\left\{ \begin{aligned}
&\lambda\unew \ge 0, \quad \varphi(\sigma\unew,p\unew) \leq 0, 
\quad \lambda\unew\varphi(\sigma\unew,p\unew)=0,\\
&\delta p := p\unew-p\uold =\lambda\unew, \quad 
\delta\vesp:=\vesp\unew-\vesp\uold = \lambda\unew \frac{\partial \varphi}{\partial \sigma}(\mathbb{C}:(\veps\unew-\vesp\uold)).
\end{aligned} \right.
\end{equation}
The main difference with respect to the procedure \texttt{PLAS\_IMP} described
in~\eqref{eq:PLAS_IMP} is on the increment of the tensor of remanent plastic strain.

\subsection{Space-time fully discrete scheme}

Full space-time discretization is achieved by combining the time-integration scheme \cite{MARAZZATO2019906}
described in the previous section with the DEM space discretization scheme from 
Section~\ref{sec:dofs}.
For all $n=1,\ldots,N$, we compute a discrete displacement field $u_h^n\in V_{hD}$ and a discrete velocity field $v_h^{n+\frac12} \in W_{hD}(t_{n+\frac12})$ (recall
that these spaces depend on $n$ if the prescribed Dirichlet condition on the displacement
is time-dependent).
Moreover, we compute
a (trace-free) tensor of remanent plastic strain $\vespc^{n,k}$ and a scalar cumulated 
plastic deformation $p_c^{n,k}$ for every mesh cell $c\in\mathcal{C}$ and every Gauss time-node
$k\in\mathcal{K}$. We set $\vespC^{n,k}:=(\vespc^{n,k})_{c\in\mathcal{C}}$ and
$p_{\mathcal{C}}^{n,k}:=(p_c^{n,k})_{c\in\mathcal{C}}$.
The fully discrete scheme reads as follows: For all $n=1,\ldots,N$, given
$u_h^n$, $v_h^{n-\frac12}$ and $v_h^{n+\frac12}$, compute $\{u_h^{n,k}\}_{k\in\mathcal{K}}$,
$v_h^{n+\frac32}$, $\{\vespC^{n,k}\}_{k\in\mathcal{K}}$, and 
$\{p_{\mathcal{C}}\}_{k\in\mathcal{K}}$ such that
\begin{equation} 
\left\{ \begin{alignedat}{2}
&u_h^{n,k} = u_h^n + (t_{n,k}-t_n)v_h^{n+\frac12},&\quad&\forall k\in\mathcal{K}, \\
&(\vespc^{n,k},p_c^{n,k}) = \texttt{PLAS\_EXP}(\vespc^{n,k-1},p_c^{n,k-1},\veps_c^{n,k}),&\quad&\forall k\in\mathcal{K},\ \forall c\in\mathcal{C}, \\
&\frac12 m_h(v_h^{n+\frac32}-v_h^{n-\frac12},\tilde{v}_h) = \sum_{k\in\mathcal{K}} \omega_{n,k} \Big( 
l_h(t_{n,k},\tilde{v}_h) - a_h(\vespC^{n,k};u_h^{n,k},\tilde{v}_h)\Big),&\quad&\forall \tilde{v}_h\in V_{h0}, \\
\end{alignedat} \right.
\end{equation}
where $\veps_c^{n,k}:=\veps_c(\mathcal{R}(u_h^{n,k}))$. Moreover, 
$m_h$ and $l_h$ are, respectively, the discrete mass bilinear form and the 
discrete load linear form. For the first Gauss node $k=1$, the first two arguments
in \texttt{PLAS\_EXP} come from the previous time step or the initial condition. 
The initial displacement $u_h^0$ and the initial velocity $v_h^{\frac12}$ are evaluated by
using the values of the prescribed initial displacement $u_0$ and the prescribed 
initial velocity $v_0$ at the cell barycentres and the boundary vertices. 

The time step is restricted by the following CFL stability condition:
\begin{equation}
\label{eq:CFL stability MEMM}
\Delta t <  2\sqrt{\frac{\mu_{\min{}}}{\lambda_{\max{}}}},
\end{equation}
where $\mu_{\min{}}$ is the smallest entry of the diagonal mass matrix associated with the discrete mass bilinear form $m_h(\cdot,\cdot)$ and $\lambda_{\max{}}$ is the largest 
eigenvalue of the stiffness matrix associated the discrete stiffness bilinear form 
$a_h(0;\cdot,\cdot)$ (i.e., this maximal eigenvalue is computed in the worst-case 
scenario when there is no plasticity). The CFL condition~\eqref{eq:CFL stability MEMM}
guarantees the stability of the time-integration scheme
in the linear case \cite{MARAZZATO2019906}, i.e.,
when there is no plasticity. We expect that this
condition is still reasonable in the nonlinear case since plasticity does not
increase the stiffness of the material. 

Finally, let us write the discrete equivalent of the energy conservation property \eqref{eq:energy_balance}.
Define the discrete elastic energy at $t_{n+1}$ as
\[ E_{\mathrm{elas},h}^{n+1}:=\frac12 \sum_{c\in\mathcal{C}} |c| \sigma_{c}^{n+1} : 
\mathbb{C}^{-1} : \sigma_{c}^{n+1}, \]
 with $\sigma_{c}^{n+1}:=\mathbb{C}:(\veps_c(\mathcal{R}(u_h^{n+1}))-\vespc^{n+1})$, for all $c \in \mathcal{C}$. The discrete kinetic energy at $t_{n+1}$ is defined by 
 \[ E_{\mathrm{kin},h}^{n+1}:=\frac12 m_h(v_h^{n+3/2},v_h^{n+1/2}), \]
 the discrete plastic dissipation at $t_{n+1}$ is defined by
\[ E_{\mathrm{plas},h}^{n+1} := \sum_{c\in\mathcal{C}} |c| \left( \sigma_0 p_c^{n+1} +\omega_p(p_c^{n+1}) \right), \]
and finally the work of external loads at $t_{n+1}$ is defined by
\[ E_{\mathrm{ext},h}^{n+1} := \sum_{m=0}^{n} l_h(t_{m+1}; v_h^{m+1/2}) \Delta t. \]
Then assuming a homogeneous Dirichlet condition for simplicity, we have the following energy balance equation:
\begin{equation}
\label{eq:discrete conservation energy}
E_{\mathrm{elas},h}^{n+1} + E_{\mathrm{kin},h}^{n+1} + E_{\mathrm{plas},h}^{n+1}
= E_{\mathrm{ext},h}^{n+1} + E_{\mathrm{elas},h}^{0} + E_{\mathrm{kin},h}^{0} + \mathcal{O}(\Delta t^2),
\end{equation}
where the term $\mathcal{O}(\Delta t^2)$ results from the use of quadratures to compute the integral of forces (recall that $E_{\mathrm{plas},h}^{0} = 0$ in our setting).
The interested reader is referred to \cite{MARAZZATO2019906} for further details on the effect of quadratures on the conservation properties of the integration method.

\subsection{Numerical tests}
\label{sec:dynamic numerical tests}

This section contains two three-dimensional test cases: a beam in dynamic flexion 
and a beam in dynamic torsion.
We use the midpoint quadrature for the integration of the forces in each time step
within the time-integration scheme. We refer
the reader to \cite{MARAZZATO2019906} for a study of the influence of the quadrature on the 
scheme accuracy for various nonlinear problems with Hamiltonian dynamics. 

Although the material parameters are indicated below using physical units, they are best interpreted in terms of characteristic times. For instance, considering a one-dimensional domain of length $L$, the characteristic time of the numerical experiments is $T_{\mathrm{ref}} := L \sqrt{\frac{\rho}{E}}$. The CFL condition \eqref{eq:CFL stability MEMM} restricts the time-step as $\Delta t < 2 \sqrt{\frac{\mu_{\min{}}}{\lambda_{\max{}}}}$. Thus, one has
\[ \frac{\Delta t}{T_{\mathrm{ref}}} < \frac{2}{L} \sqrt{ \frac{\mu_{\min{}} E}{\lambda_{\max{}} \rho} }. \]
Since the ratio $\frac{E}{\lambda_{\max}}$ is independent of $E$, the same conclusion holds for $\frac{\Delta t}{T_\mathrm{ref}}$. Therefore we will report this time ratio in the computations.

\subsubsection{Beam in dynamic flexion}
\label{sec:dynamic flexion}

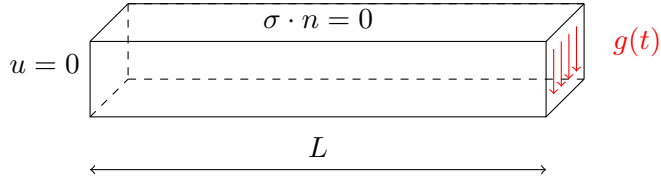
\begin{figure} [!htp]
\begin{center}
\begin{tikzpicture}
\coordinate (a) at (-3,-0.5);
\coordinate (b) at (-3,0.5);
\coordinate (c) at (3,0.5);
\coordinate (d) at (3,-0.5);
\coordinate (ap) at (-2.5,0);
\coordinate (bp) at (-2.5,1.);
\coordinate (cp) at (3.5,1.);
\coordinate (dp) at (3.5,0);

\draw (a) -- (b) -- (c) -- (d) -- cycle;
\draw[dashed] (a) -- (ap);
\draw (b) -- (bp);
\draw[dashed] (bp) -- (ap);
\draw (c) -- (cp);
\draw (d) -- (dp);
\draw (cp) -- (dp);
\draw (bp) -- (cp);
\draw[dashed] (bp) -- (cp);
\draw[dashed] (ap) -- (dp);
\draw [<->] (-3,-1.2) -- (3,-1.2);
\draw (0,-0.9) node{$L$}; 
\draw (4.2,0.5) node[red] {$g(t)$};
\draw [<-,red] (3.4, 0.1) -- (3.4, 0.7);
\draw [<-,red] (3.3, -0.) -- (3.3, 0.6);
\draw [<-,red] (3.2, -0.1) -- (3.2, 0.5);
\draw [<-,red] (3.1, -0.2) -- (3.1, 0.4);
\draw (0,0.8) node{$\sigma \cdot n = 0$};
\draw(-3,0.2) node[left] {$u=0$};
\end{tikzpicture}
\caption{Beam in dynamic flexion: problem setup.}
\label{fig:dynamical test case sketch}
\end{center}
\end{figure}

This test case consists in computing the oscillations of an elastic and linearly isotropic plastic beam of length $L=1\mathrm{m}$ with a rectangular section of $0.04 \times 0.1 \mathrm{m}^2$. The simulation time is $T=2.5\mathrm{s}$. The beam is clamped at one end, it is loaded by a uniform vertical traction $g(t)$ at the other end, and the four remaining lateral faces are stress free ($\sigma\cdot n=0$). The load term $g(t)$ is defined as
\begin{equation}
g(t) :=
\begin{cases}
-\frac{5t}{4}e_x& \text{for $0 \leq t \leq \frac45$}, \\
0&\text{for $\frac45 \leq t \leq T$}.
\end{cases}
\end{equation}
Figure \ref{fig:dynamical test case sketch} displays the problem setup.
The material parameters are $E = 10^3 \mathrm{Pa}$ for the Young modulus, $\nu = 0.3$ for the Poisson ratio, $\rho = 1\mathrm{kg {\cdot} m^{-3}}$ for the density, $\sigma_0 = 25\mathrm{Pa}$ for the yield stress, and $E_t = \frac{1}{100}E$ for the tangent plastic modulus. 
The present three-dimensional implementation used as a starting point \cite{bleyer2018numericaltours}, where $P^1$-Lagrange FE and an implicit time-integration scheme are considered for a purely elastic material.

The proposed DEM is compared to penalised Crouzeix--Raviart FE \cite{hansbo2003discontinuous}. This method is chosen as reference since it is known to be robust in the incompressible limit as $\nu \to 0.5$.
The penalty parameter is chosen as $\eta=\beta\mu$ with $\beta=0.5$. The reference computation is performed using $14,376$ vector-valued dofs and a time-step $\Delta t = 20 \mathrm{\mu s}$.
Two computations are performed with the proposed DEM. The first uses a coarse mesh containing $4,668$ vector-valued dofs and a time-step $\Delta t = 1.4\mathrm{\mu s}$, which is stable for the explicit integration. The second uses a fine mesh containg $13,302$ vector-valued dofs and a time-step $\Delta t = 1.1\mathrm{\mu s}$, also stable for the explicit integration.
Thus one has: $4.4 \cdot 10^{-7} \leq \frac{\Delta t}{T_{\mathrm{ref}} } \leq 8.0 \cdot 10^{-6}$.
The penalty parameter for both computations is similar to the reference computation with $\beta=0.5$.
As already mentioned, we used a midpoint quadrature of the forces. Higher-order symmetric quadratures have been found to give overlapping results with respect to the midpoint quadrature. In all the computations, the time-discretization error is actually smaller than the space-discretization error, but larger time-steps cannot be considered owing to the CFL restriction. 

\begin{figure}[htp]
\centering
\subfloat{
\includegraphics[width=0.5\textwidth]{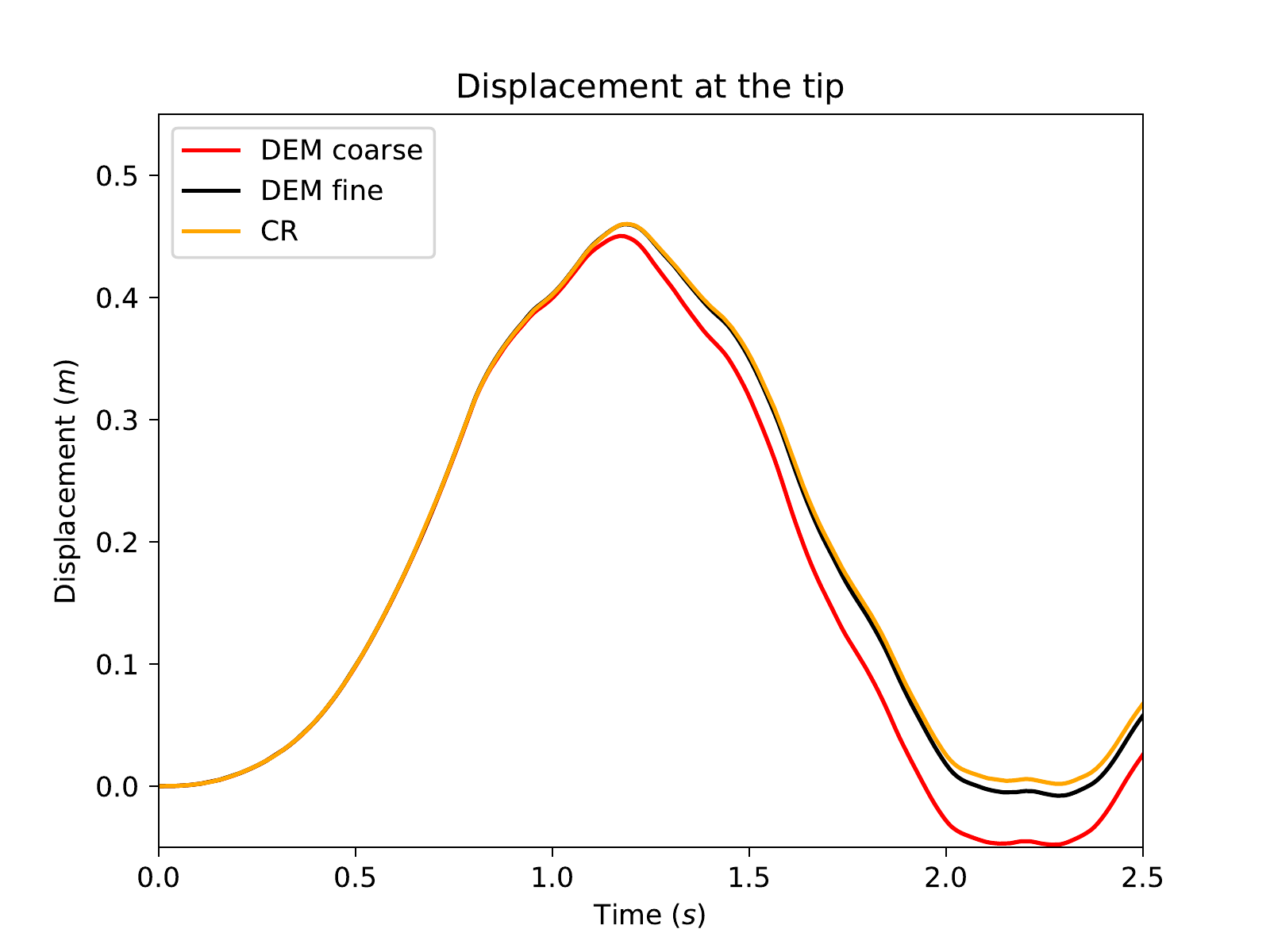}
}
\subfloat{
\includegraphics[width=0.5\textwidth]{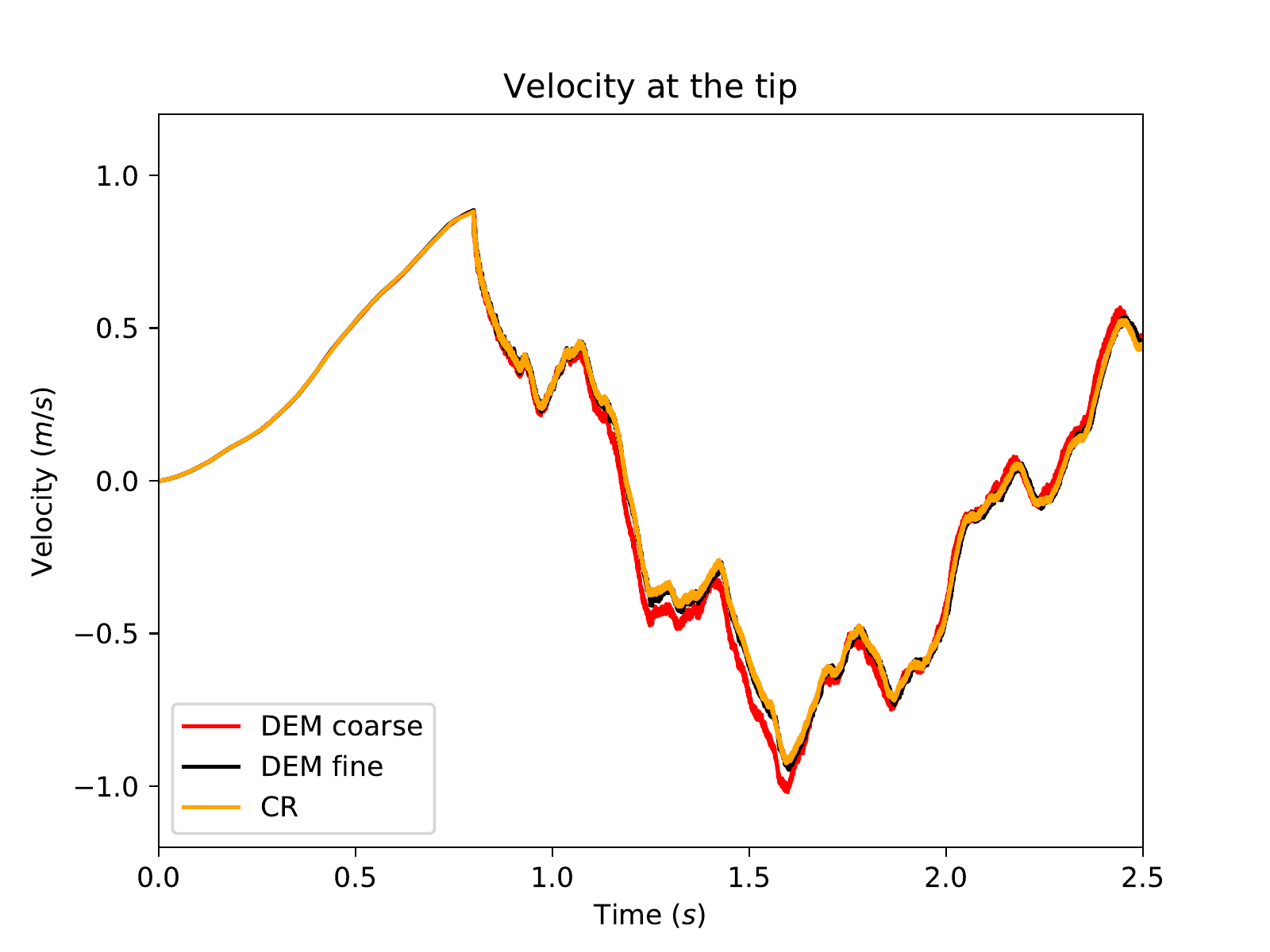}
}
\caption{Beam in dynamic flexion: comparison between the proposed scheme (DEM) on a coarse and a fine mesh and the reference solution (CR). Left: Displacement at the loaded tip of the beam. Right: Velocity at the same point.}
\label{fig:comparaison dynamique coarse}
\end{figure}
\begin{figure}[htp]
\centering

\subfloat{
\includegraphics[width=0.5\textwidth]{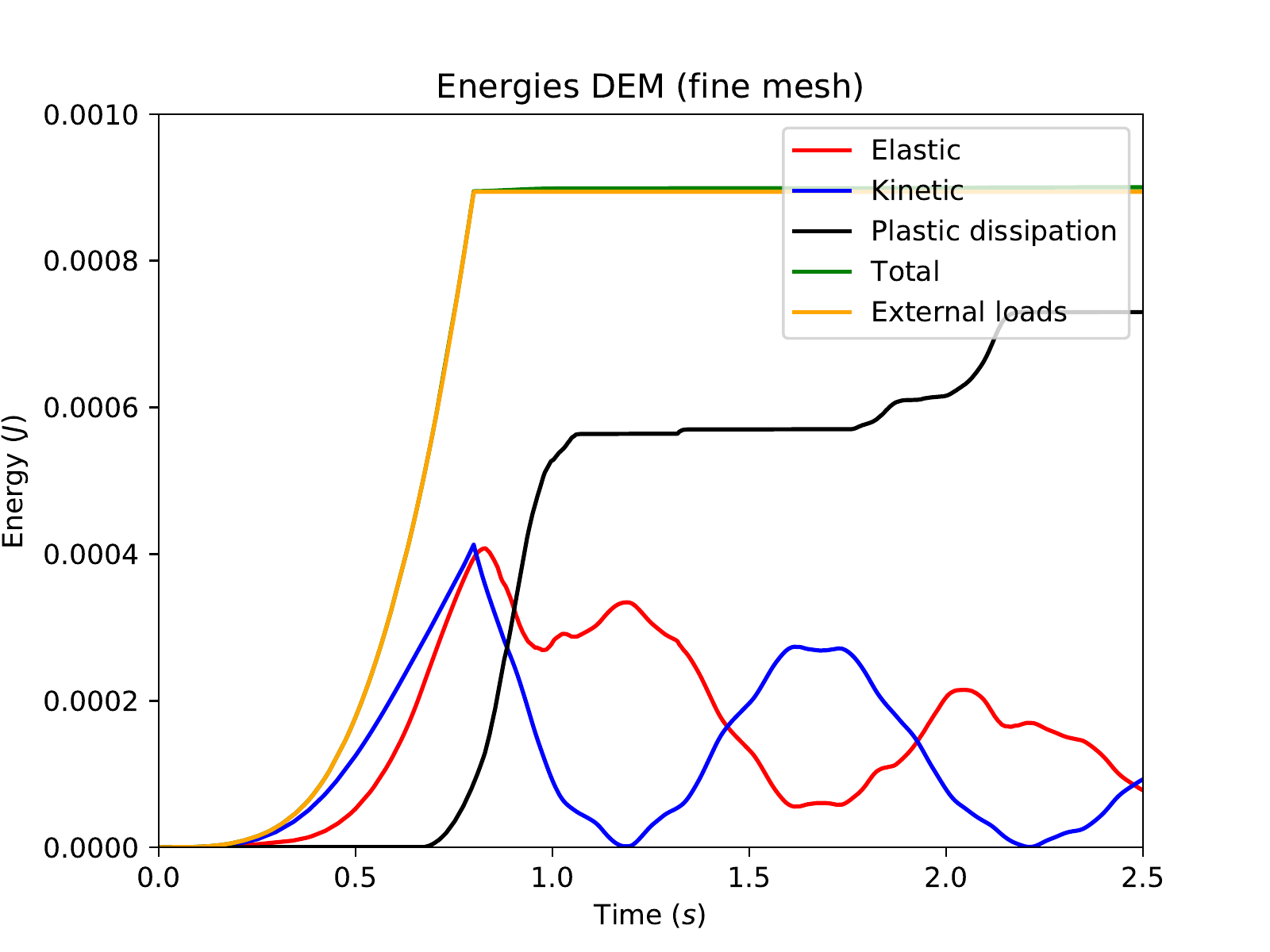}
}
\subfloat{
\includegraphics[width=0.5\textwidth]{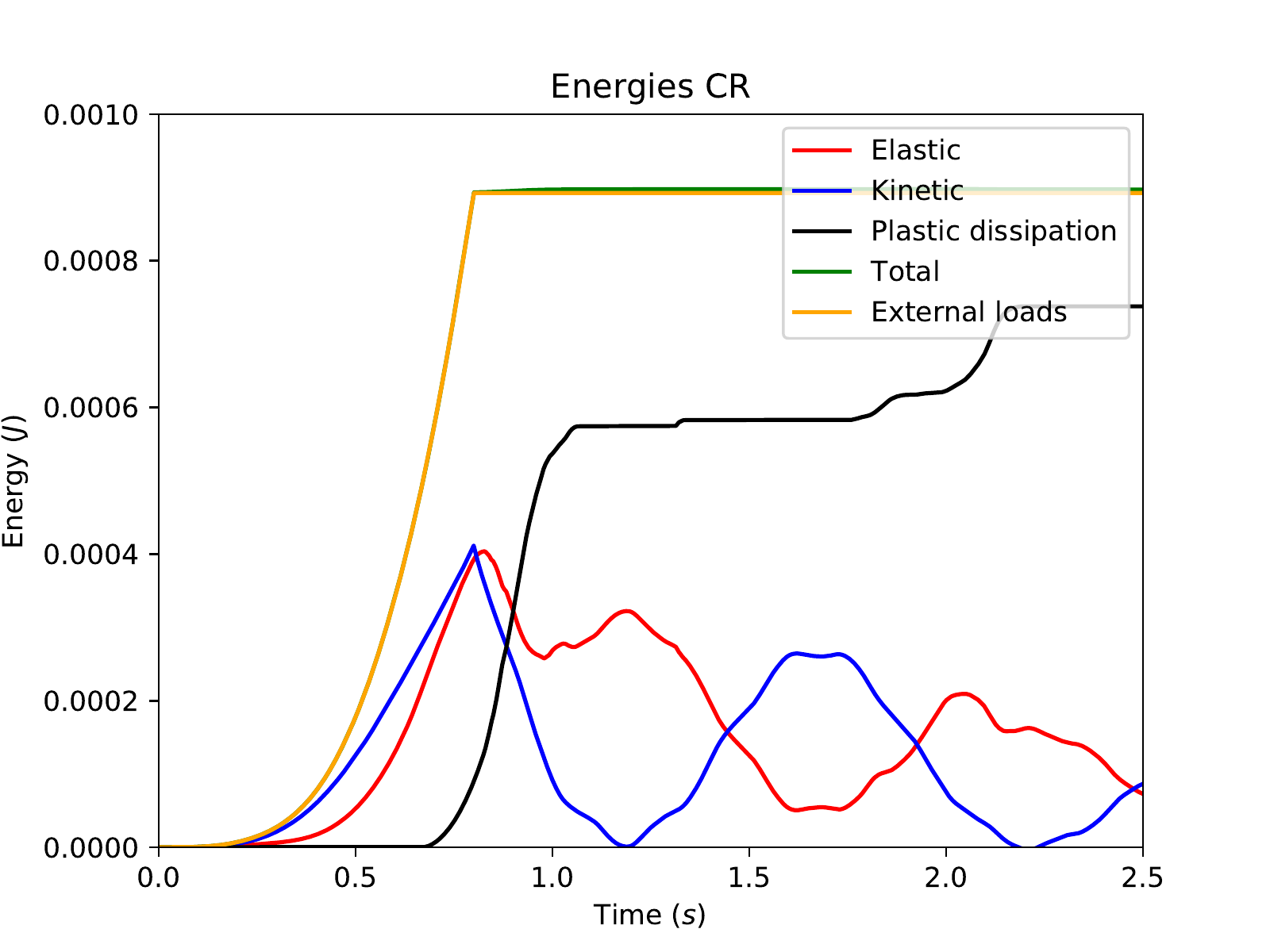}
}
\caption{Beam in dynamic flexion: energies during the simulation. Left: DEM (fine mesh). Right: reference solution (CR).}
\label{fig:comparaison dynamique energies flexion}
\end{figure}

The displacement and the velocity at the center of the loaded tip of the beam are compared in Figure \ref{fig:comparaison dynamique coarse}.
We notice the excellent agreement between the DEM prediction on the fine mesh and the reference computation.
Figure \ref{fig:comparaison dynamique energies flexion} shows the balance of energies for the reference computation and the fine DEM computation. 
One can first notice that the total energy for both DEM and Crouzeix--Raviart space semi-discretizations is well preserved by the time-integrator \cite{MARAZZATO2019906} since the total mechanical energy (kinetic energy, elastic energy and plastic dissipation) and the work of the external load are perfectly balanced at all times. We also notice that the amount of plastic dissipation is rather significant at the end of the simulation. 
\begin{figure}[htp]
\centering
\subfloat{
\includegraphics[width=0.33\textwidth]{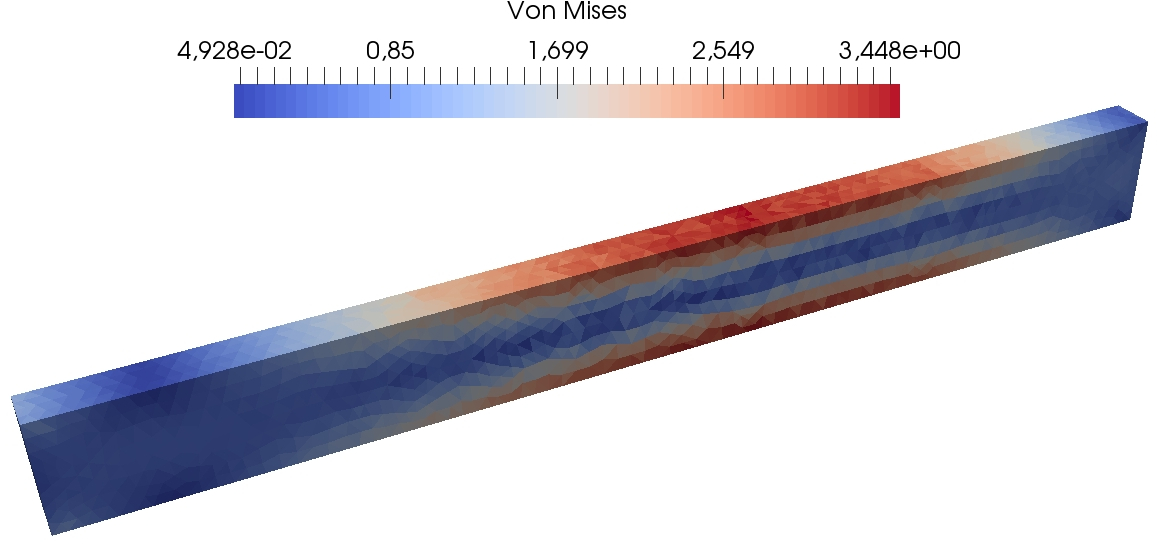}
}
\subfloat{
\includegraphics[width=0.33\textwidth]{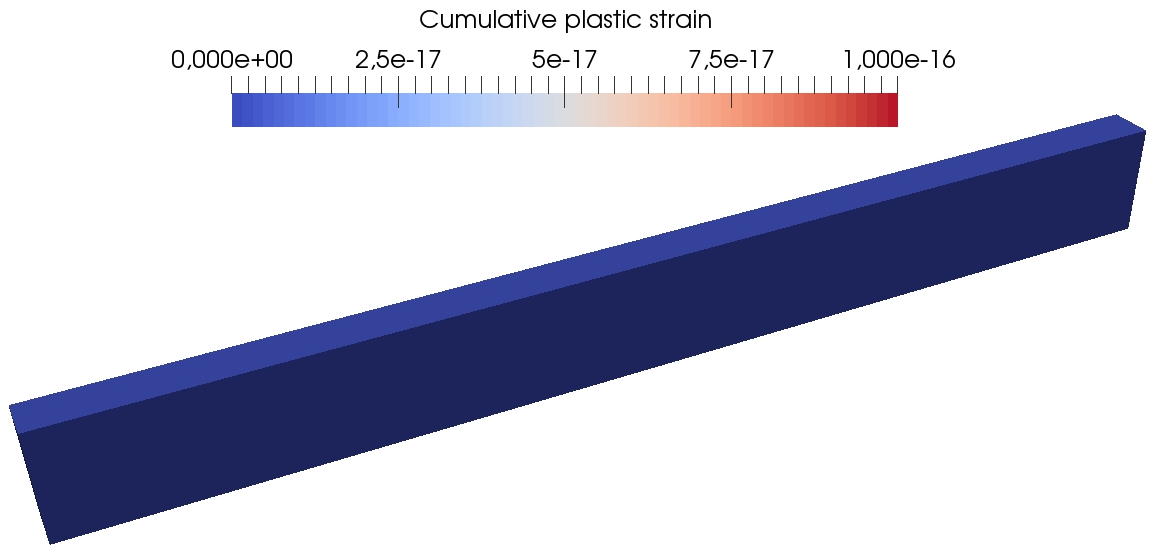}
}
\subfloat{
\includegraphics[width=0.33\textwidth]{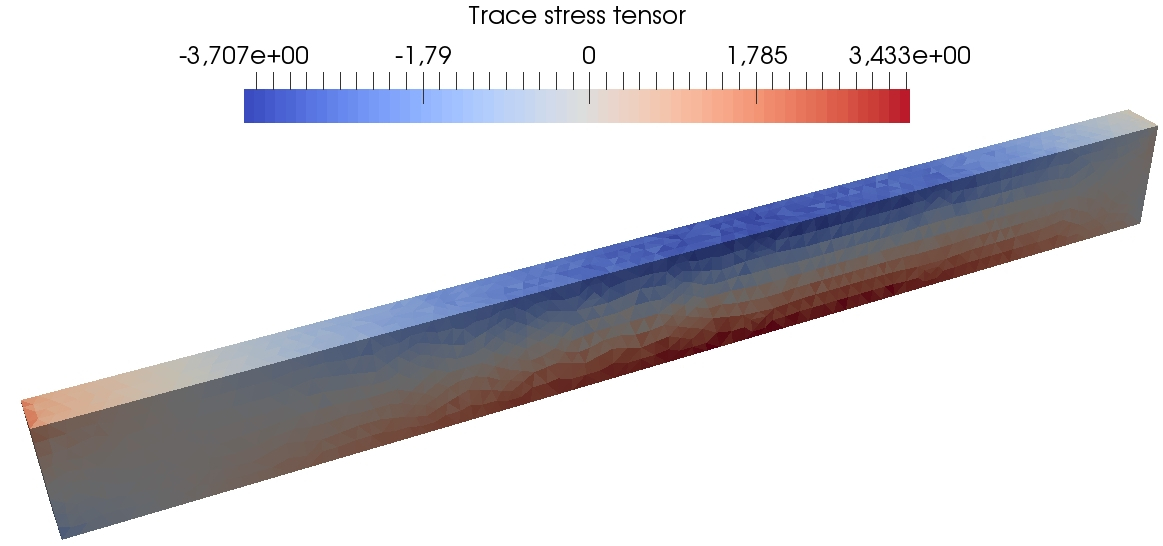}
}\\

\subfloat{
\includegraphics[width=0.33\textwidth]{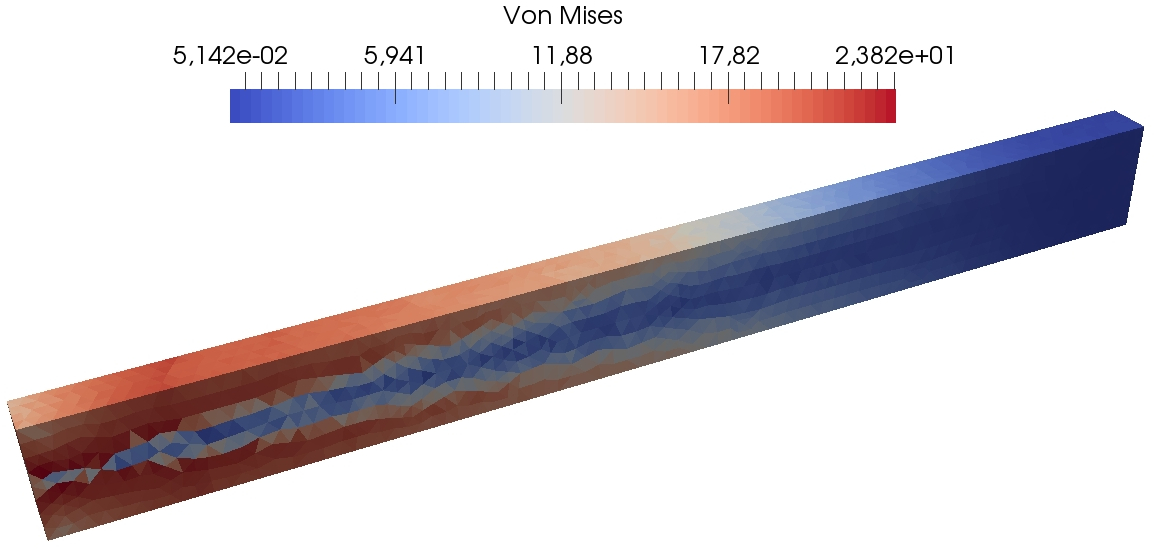}
}
\subfloat{
\includegraphics[width=0.33\textwidth]{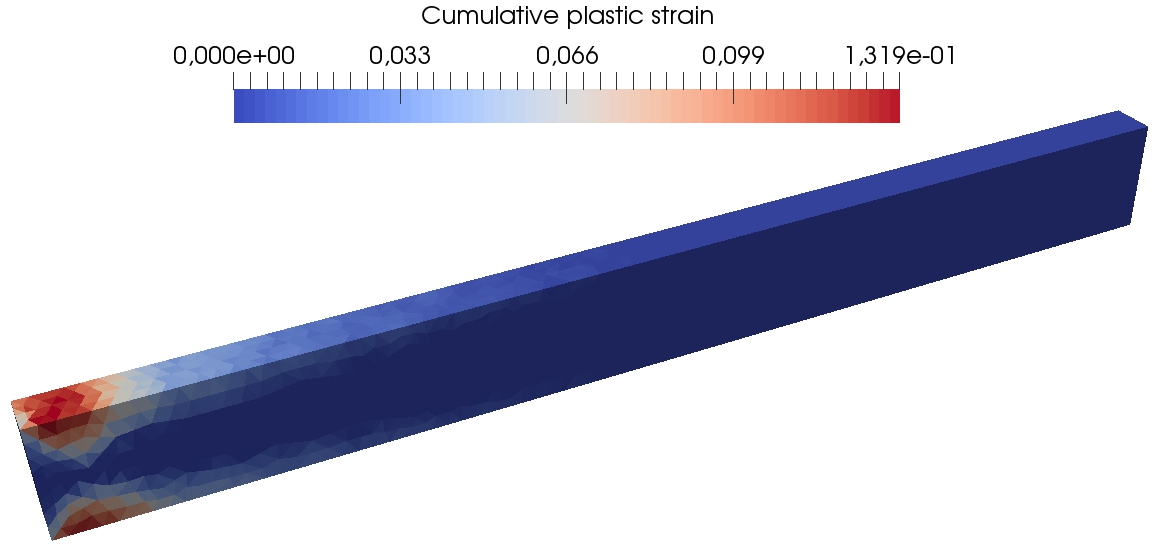}
}
\subfloat{
\includegraphics[width=0.33\textwidth]{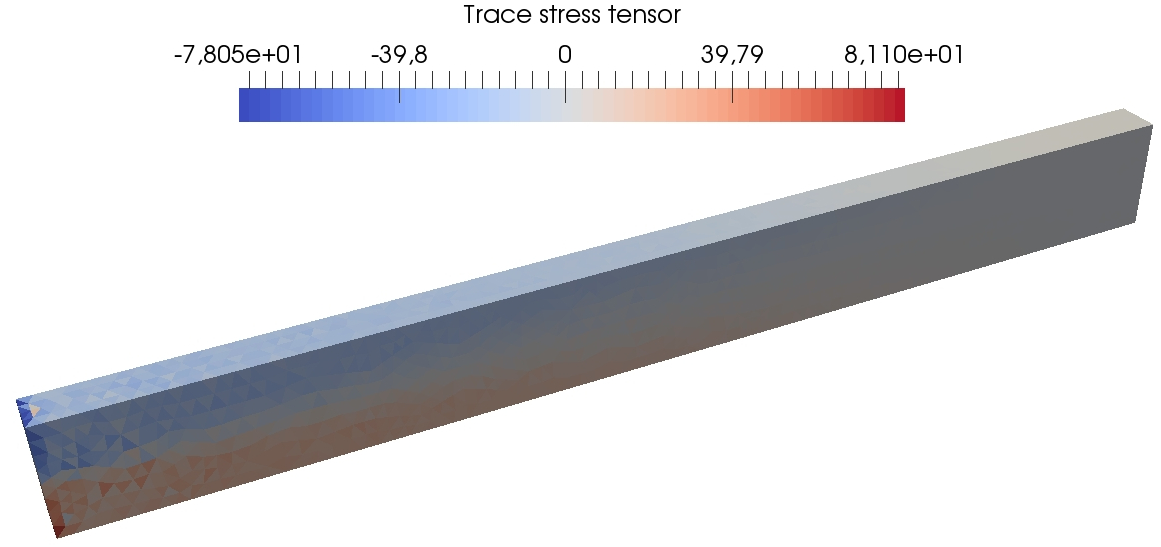}
}\\

\subfloat{
\includegraphics[width=0.33\textwidth]{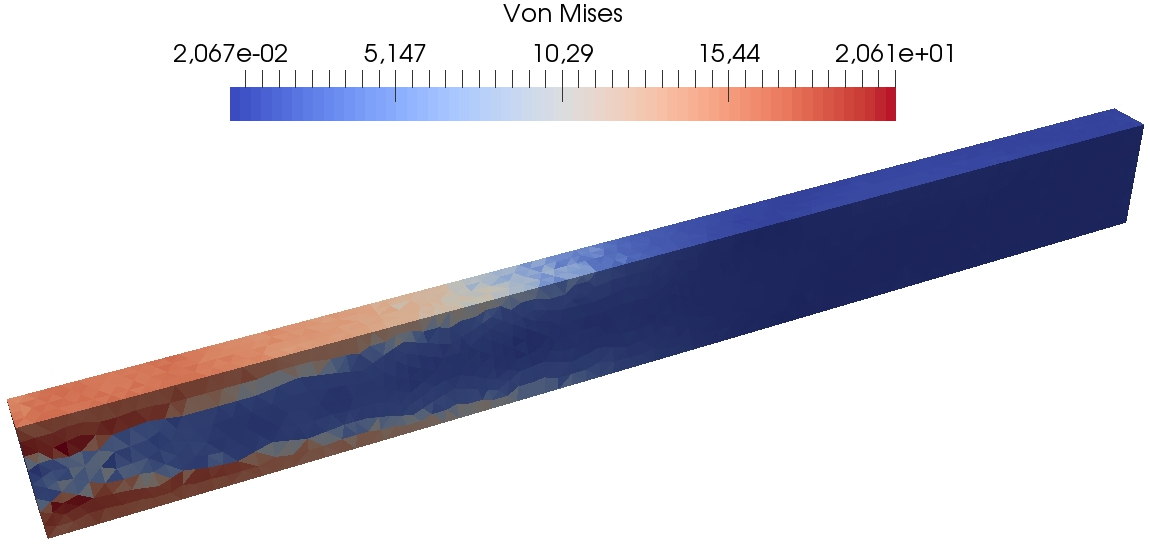}
}
\subfloat{
\includegraphics[width=0.33\textwidth]{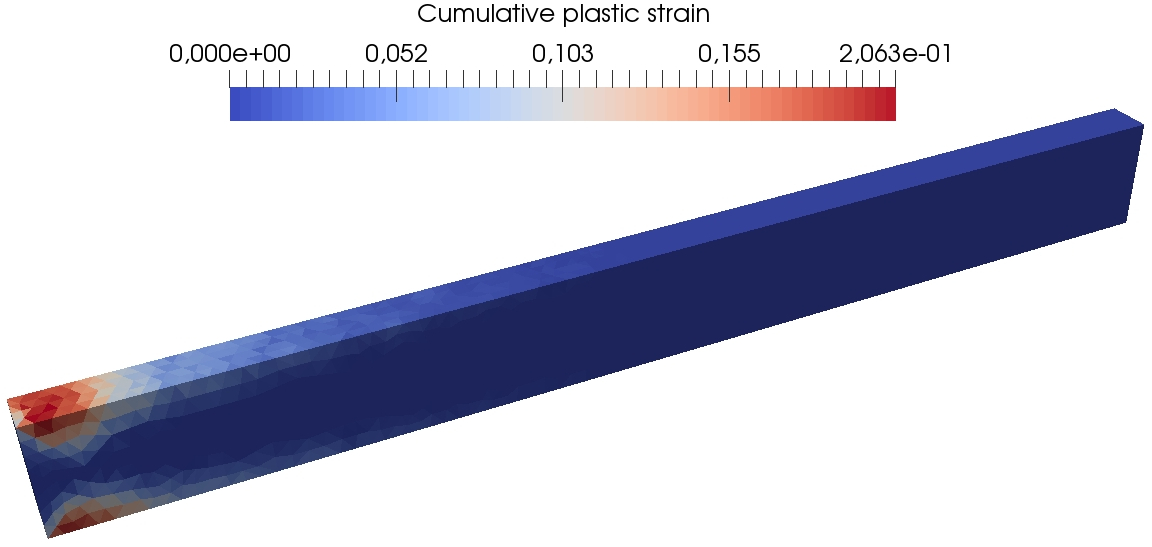}
}
\subfloat{
\includegraphics[width=0.33\textwidth]{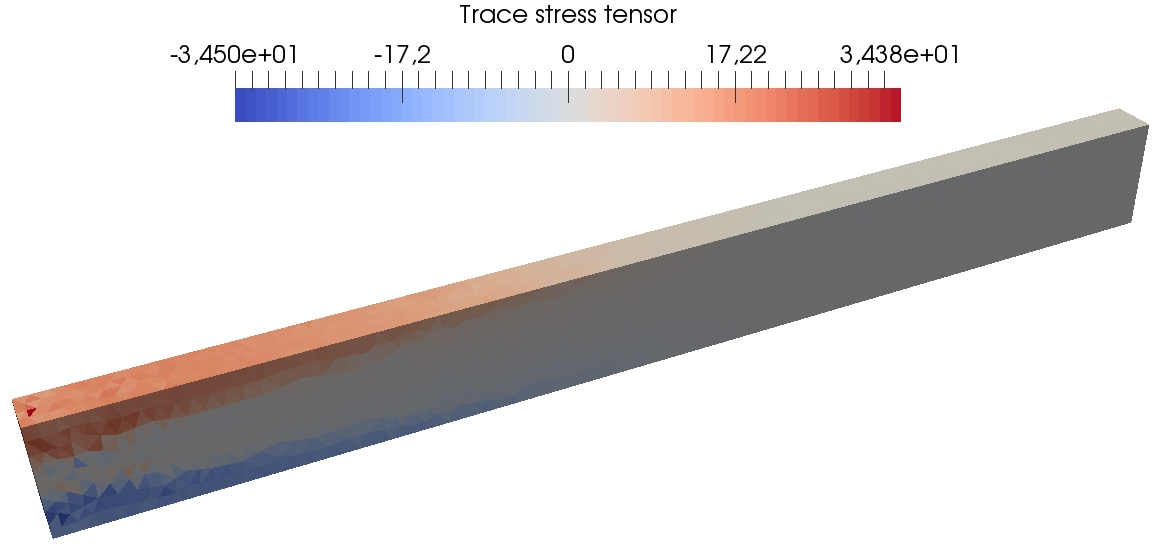}
}\\

\caption{Beam in dynamic flexion: DEM on the fine mesh. Von Mises equivalent stress (left column), $p$ (middle column) and $\mathrm{tr}(\sigma)$ (right column) at $t=\frac{1}{10}T$ (top line), $t=\frac12 T$ (middle line) and $t=T$ (bottom line).}
\label{fig:results computation flexion}
\end{figure}
Figure \ref{fig:results computation flexion} presents some further results of the DEM
computations on the fine mesh so as to visualize at three different times during the
simulation the spatial localization of the von Mises equivalent stress, the cumulated plastic strain, and the trace of the stress tensor. One can see that the plastic strain is concentrated close to the clamped tip of the beam, where the material undergoes the greatest stresses. The method does not exhibit any locking due to plastic incompressibility as indicated by the smooth behavior of the trace of the stress tensor.

\subsubsection{Beam in dynamic torsion}

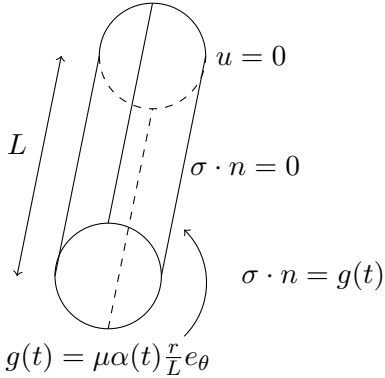
\begin{figure} [!htp]
\begin{center}
\begin{tikzpicture}
\pgfmathsetmacro{\R}{0.7}
\pgfmathsetmacro{\L}{5*\R}
\coordinate (O) at (0,0);
\coordinate (Op) at (\L/6,5*\L/6);
\coordinate (Ap) at (\L/6+\R,5*\L/6);
\coordinate (Bp) at (\L/6-\R,5*\L/6);

\draw (\R,0) arc (0:360:\R);

\draw (\R,0) -- (Ap);
\draw (-\R,0) -- (Bp);
\draw (0,\R) -- (\L/6,5*\L/6+\R);
\draw[dashed] (0,-\R) -- (\L/6,5*\L/6-\R);

\draw (Ap) arc (0:180:\R);
\draw[dashed] (Bp) arc (180:360:\R);

\draw (Ap) node[right] {$u=0$};
\draw (\L/12+\R/2+0.3,5*\L/12) node[right] {$\sigma\cdot n=0$};
\draw (\R+2.,0) node{$\sigma\cdot n = g(t)$};
\draw [->] (\R + 0.3, -0.8) arc (-45:45:1);
\draw (0,-0.4-\R) node {$g(t) = \mu \alpha(t) \frac{r}{L}e_{\theta}$}; %

\draw [<->] (-\R-0.5,0) -- (\L/6-\R-0.5,5*\L/6);
\draw (-1.2,\L/2) node{$L$};
\end{tikzpicture}
\caption{Beam in dynamic torsion: problem setup.}
\label{fig:torsion dynamic}
\end{center}
\end{figure}

The setting is similar to the quasi-static torsion test case presented in Section \ref{sec:torsion quasi-static}. The two differences are the material parameters and the plastic law which are similar to Section \ref{sec:dynamic flexion}, and the boundary conditions on one side of the beam. Figure \ref{fig:torsion dynamic} shows the problem setup. On one of its extremities the beam is clamped, and on the other extremity the following normal stress is imposed:
\begin{equation}
g(t) = \mu \alpha(t) \frac{r}{L} e_{\theta},
\end{equation}
where $r$ and $e_{\theta}$ are defined in Section \ref{sec:torsion quasi-static}.
The angle $\alpha(t)$ is increased from $0$ at $t=0$ to $5\alpha_y$ at $t=T=0.5\mathrm{s}$, 
where $\alpha_y$ is the yield angle defined in Section \ref{sec:torsion quasi-static}.
The plastic parameters are the same as those in Section \ref{sec:dynamic flexion}.

Three different space discretizations are considered for this test case: $P^1$-Lagrange FE, penalised Crouzeix--Raviart FE and the proposed DEM. The Crouzeix--Raviart computations are used as reference since the method is known to be robust with respect to the incompressible limit. The Lagrange FE computations are used to illustrate their inability to deal with large plastic (deviatoric) strains. The goal of this test case is to show the ability of the proposed DEM to deal with deviatoric plasticity. The computations are not performed on the same meshes but rather with meshes leading to a similar number of dofs so as to give comparable results for DEM and Lagrange FE, whereas the meshes used with penalised Crouzeix--Raviart FE lead to twice as many dofs since they are employed to obtain a reference solution. The number of dofs and the time-steps used in the computations are presented in Table \ref{tab:dofs and time-steps torsion}.
\begin{table}[!htp]
\begin{center}
   \begin{tabular}{ | c | c | c | c | c | c | c |}
     \hline
	 Method & \multicolumn{2}{|c|}{DEM} & \multicolumn{2}{|c|}{Lagrange FE} & \multicolumn{2}{|c|}{penalised CR FE} \\ \hline
	  Computation & coarse & fine & coarse & fine & fine & very fine \\ \hline
      Vectorial dofs & $6,978$ & $14,438$ & $6,584$ & $12,853$ & $13,052$ & $27,711$ \\ \hline
     $\Delta t \text{ (}\mu s)$ & $4.1$ & $1.3$ & $3.9$ & $2.6$ & $2.3$ & $0.42$ \\ \hline
   \end{tabular}
   \caption{Beam in dynamic torsion: number of vectorial dofs and time-step for all the computations.}
   \label{tab:dofs and time-steps torsion}
\end{center}
\end{table}
One thus has: $8.4 \cdot 10^{-7} \leq \frac{\Delta t}{T_{\mathrm{ref}} } \leq 8.2 \cdot 10^{-6}$.
All the reported time-steps are compatible with the CFL restriction. Also, for all computations, we have verified that the time discretization error is negligible with respect to the space discretization error.
The time-integration scheme \cite{MARAZZATO2019906} is used with a midpoint quadrature of the forces for the coarse computations and with a Gauss-Legendre quadrature of the forces of order 5 for the fine computations. For details on the effect of the quadrature rule, see \cite{MARAZZATO2019906}. 

\begin{figure}[htp]
\centering
\subfloat{
\includegraphics[width=0.5\textwidth]{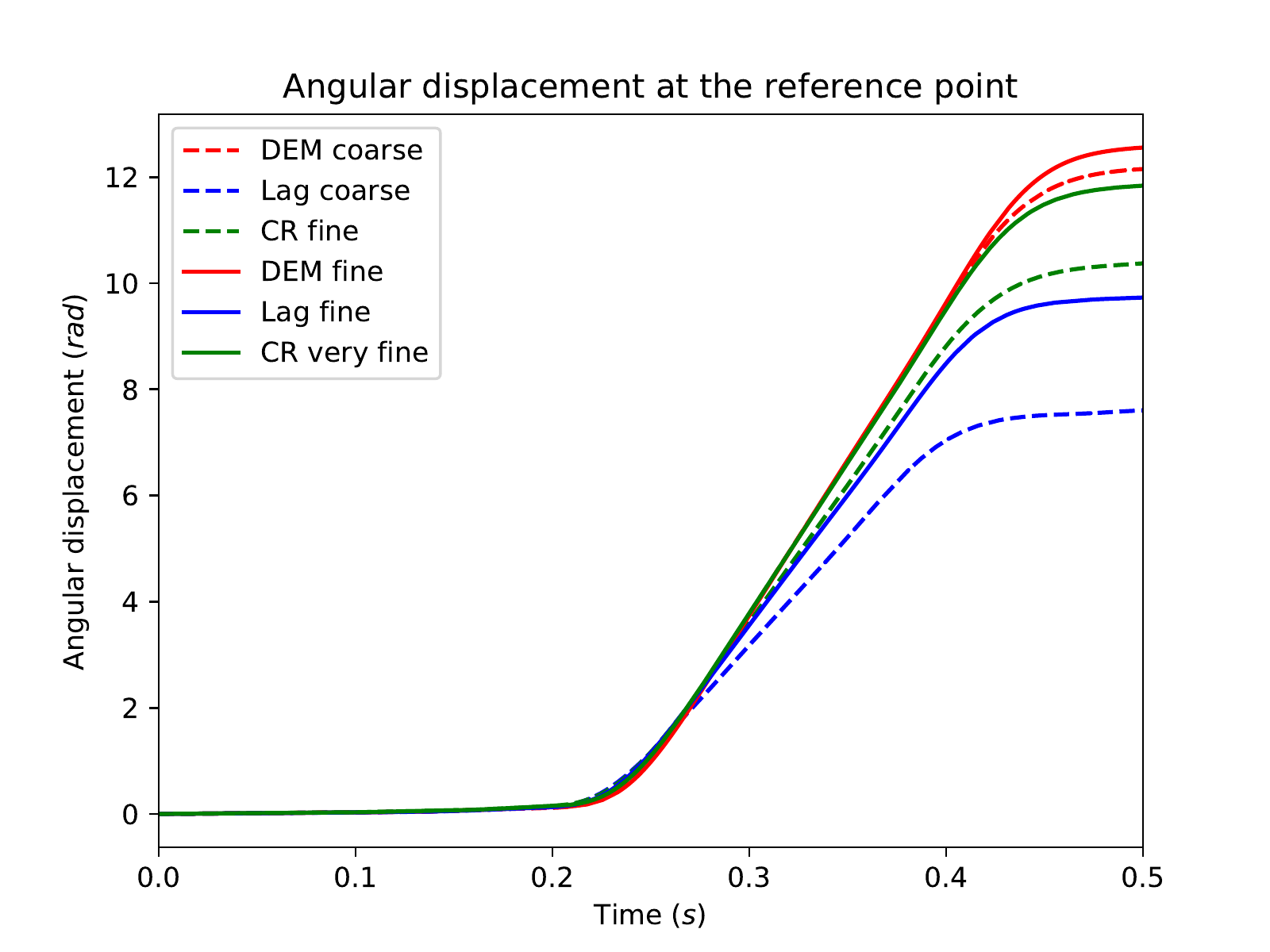}
}
\subfloat{
\includegraphics[width=0.5\textwidth]{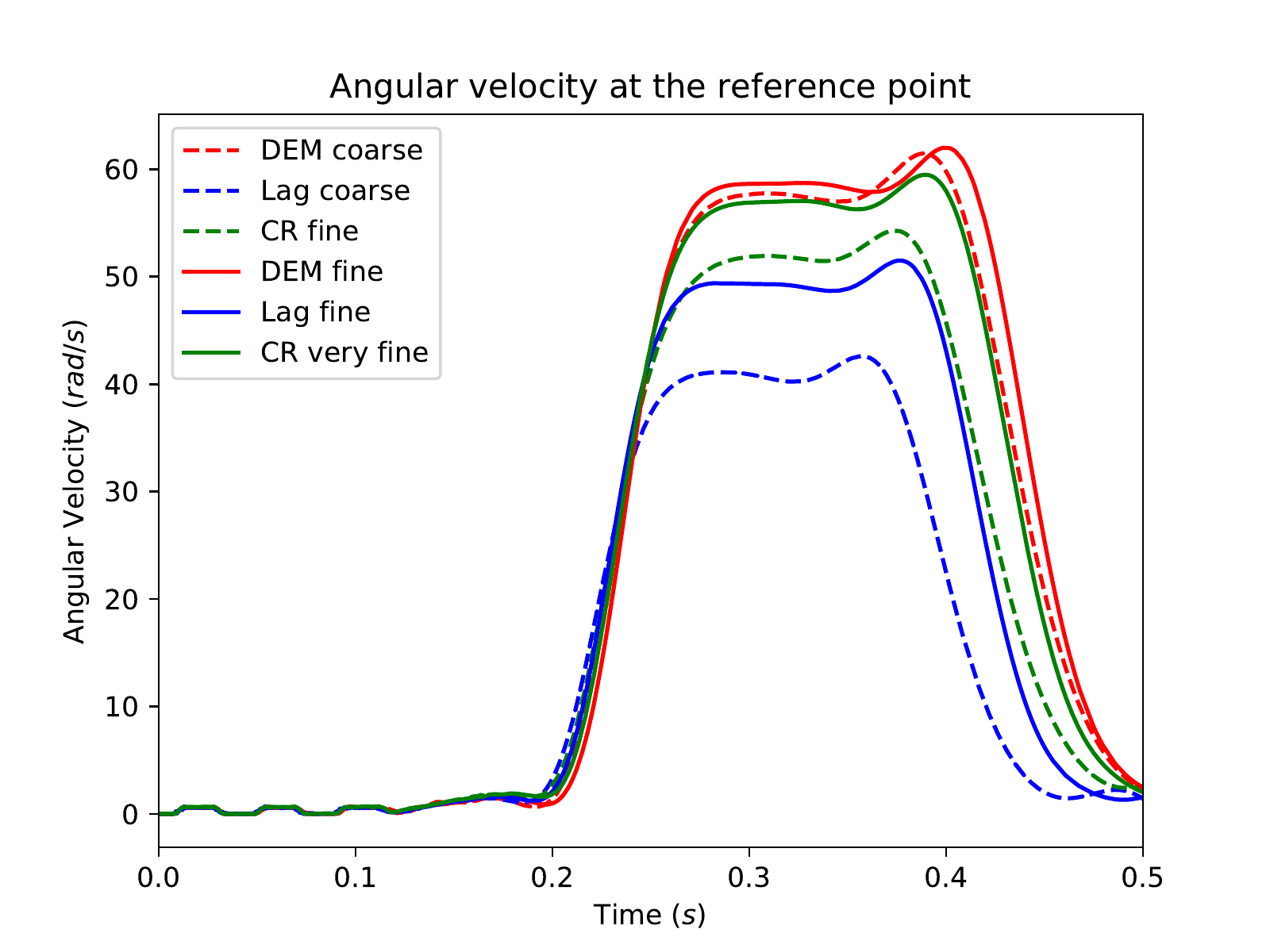}
}\\
\caption{Beam in dynamic torsion: comparison between DEM and FEM. Left: Displacement at the chosen point. Right: Velocity at the same point.}
\label{fig:comparaison dynamique 3}
\end{figure}

\begin{figure}[htp]
\centering
\subfloat{
\includegraphics[width=0.5\textwidth]{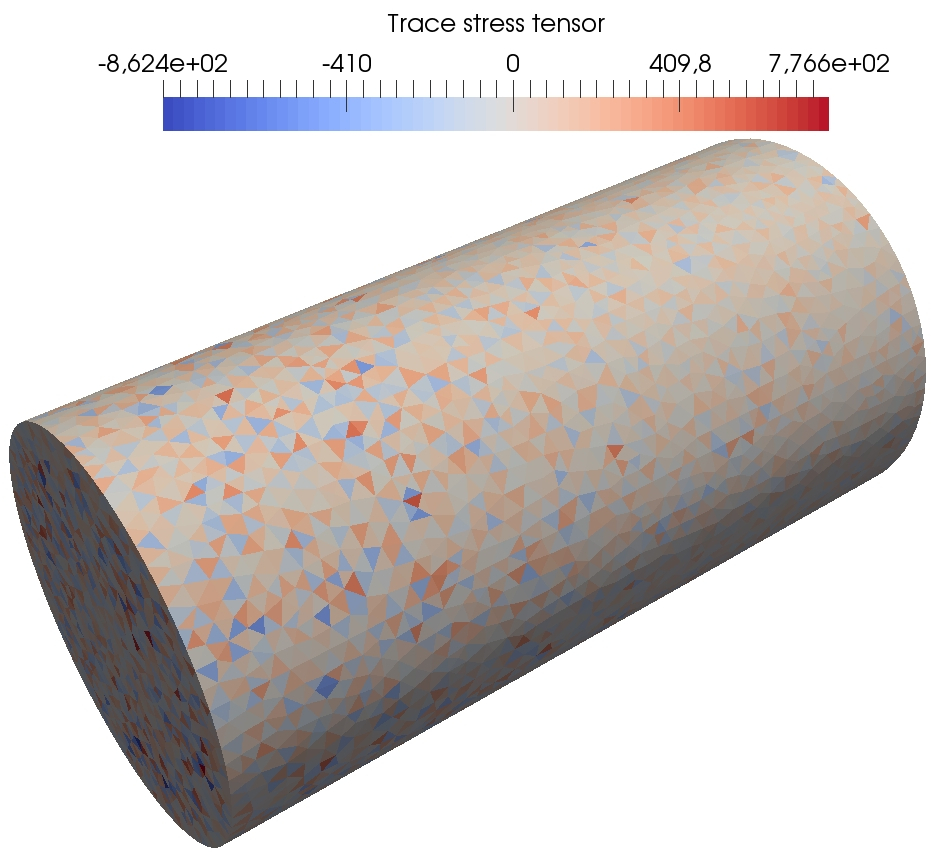} 
}
\subfloat{
\includegraphics[width=0.5\textwidth]{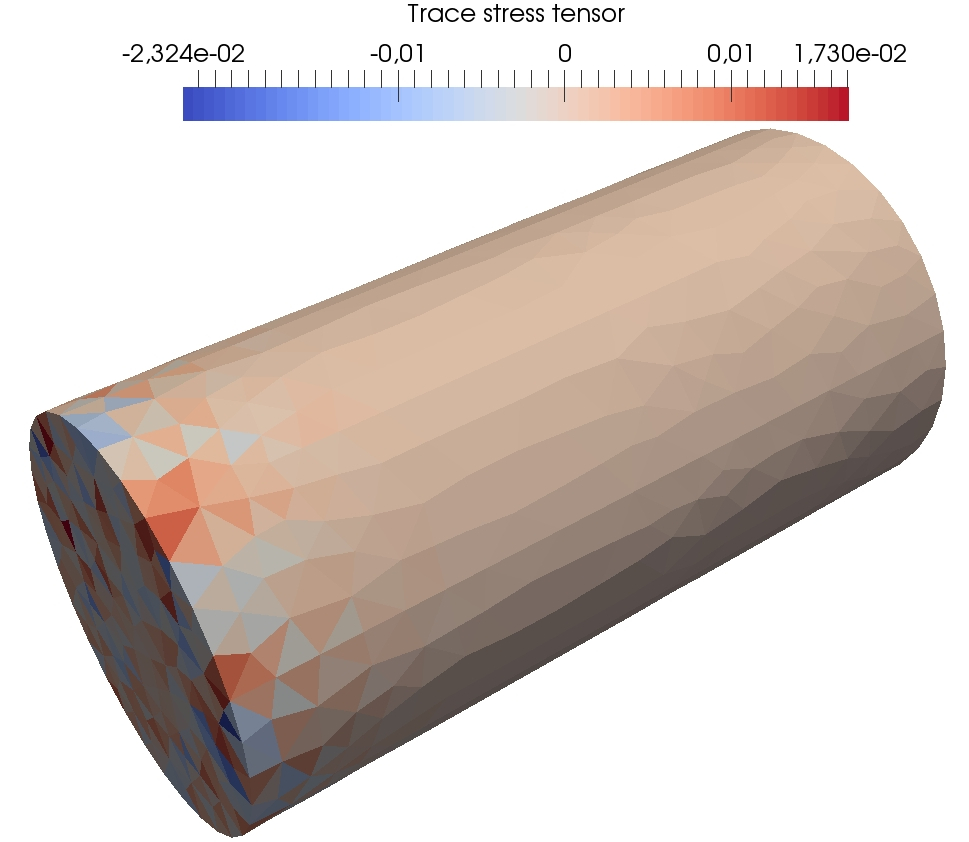}
}\\
\caption{Beam in dynamic torsion: $\mathrm{tr}(\sigma)$ at $t=\frac12 T$. Left: Lagrange FE (fine mesh). Right: Penalised Crouzeix--Raviart FE (fine mesh). }
\label{fig:comparaison torsion plasticite}
\end{figure}

\begin{figure}[htp]
\centering
\subfloat{
\includegraphics[width=0.33\textwidth]{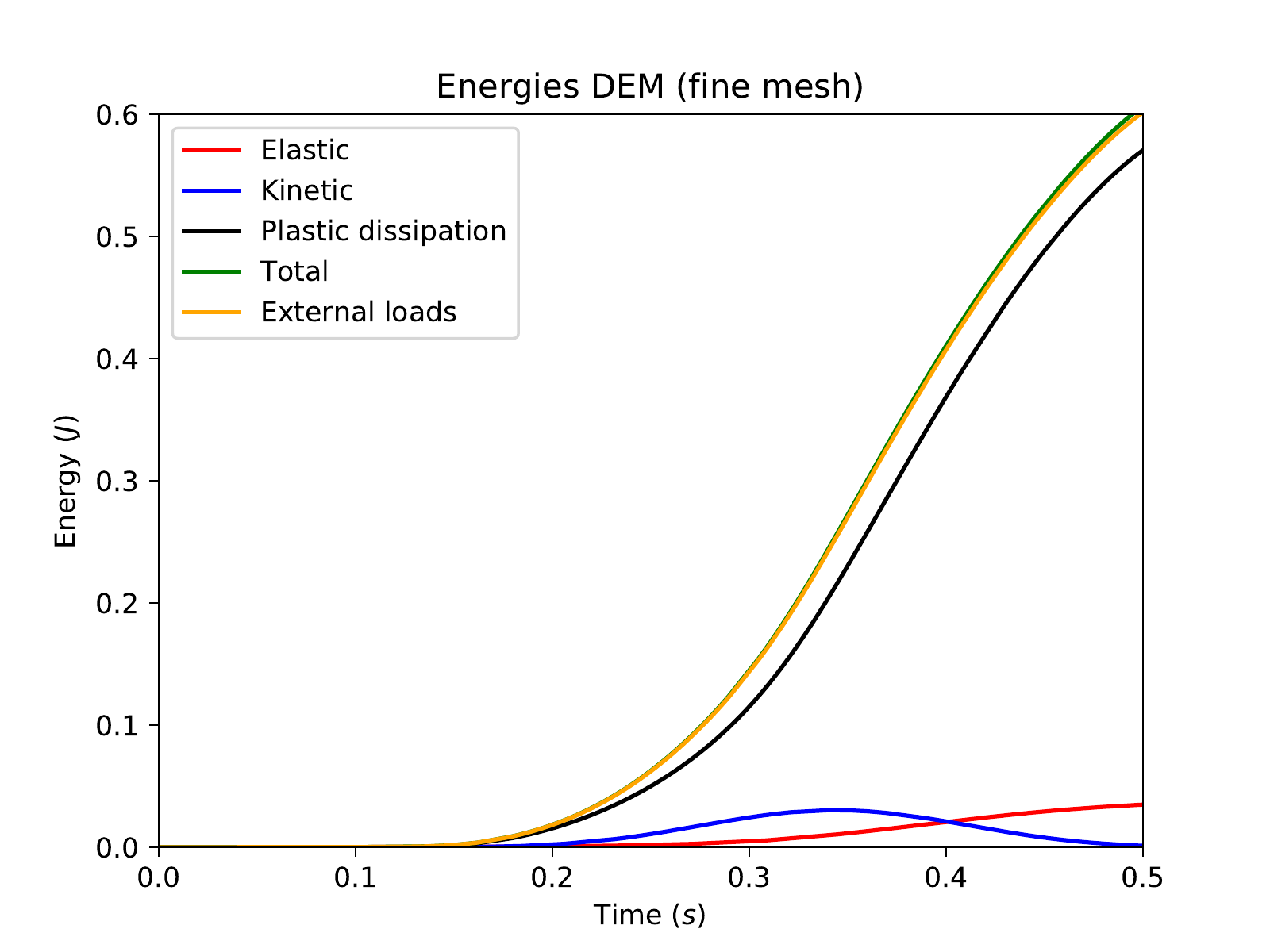}
}
\subfloat{
\includegraphics[width=0.33\textwidth]{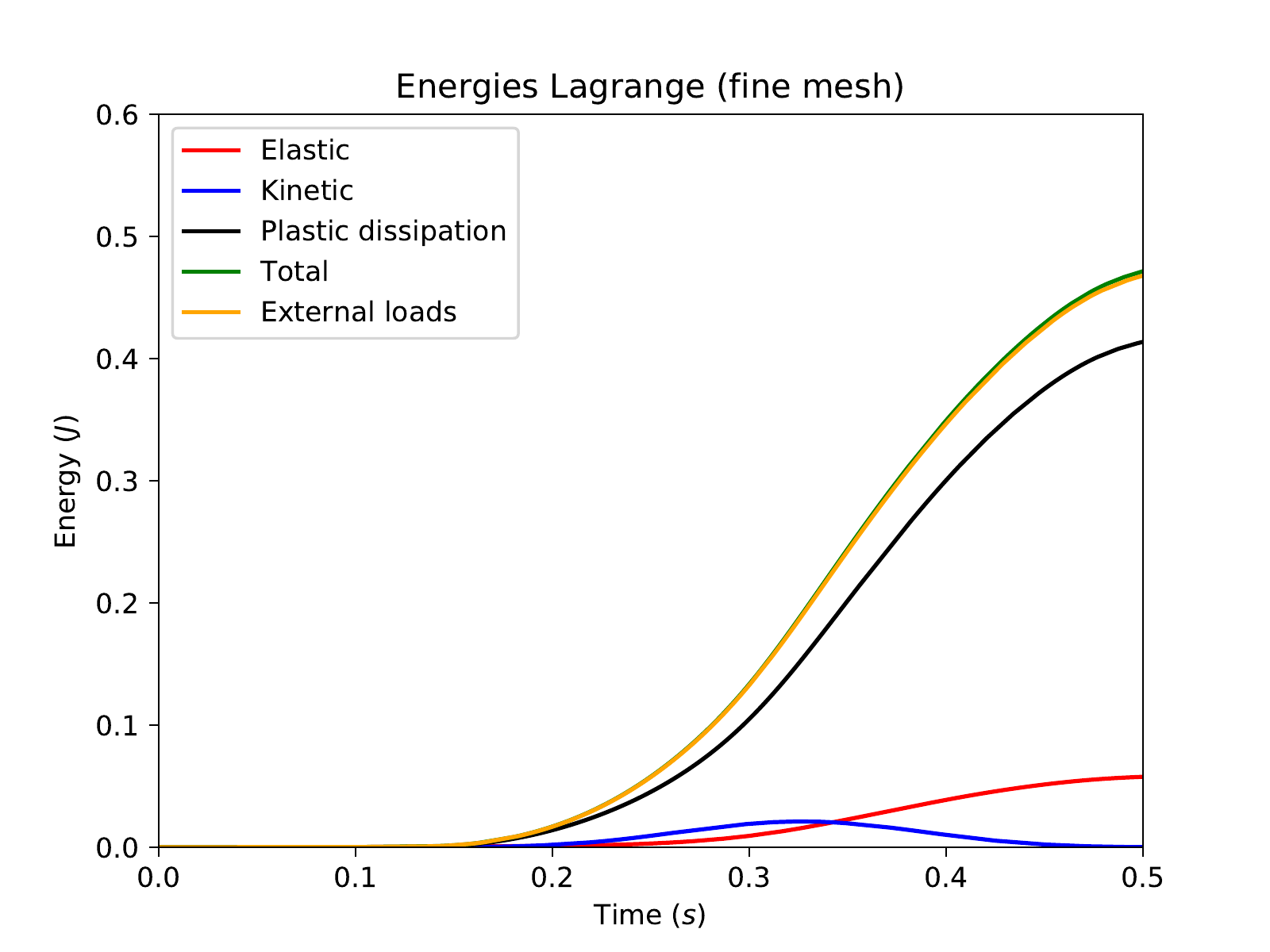}
}
\subfloat{
\includegraphics[width=0.33\textwidth]{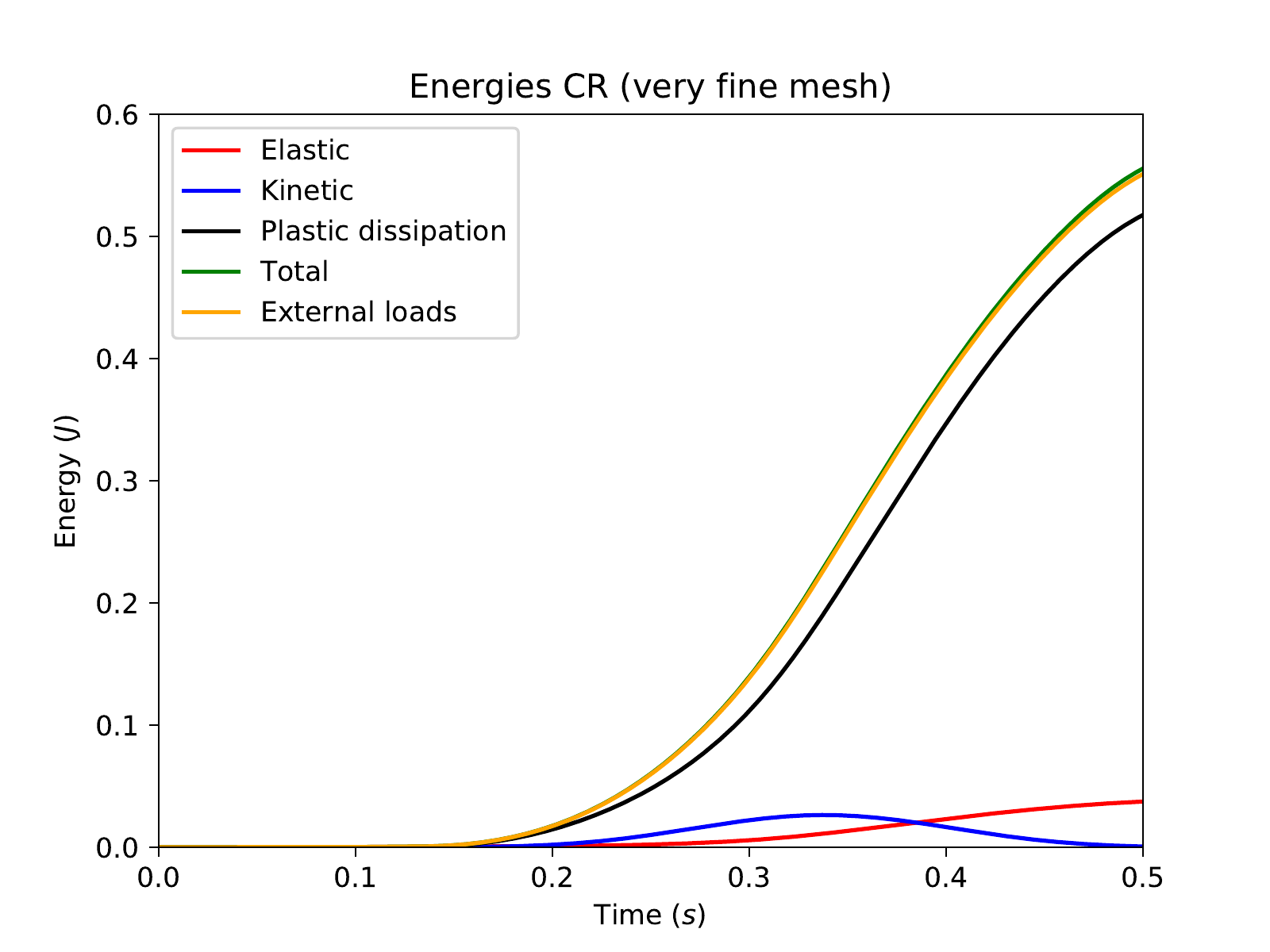}
}
\caption{Beam in dynamic torsion: energies during the simulation on the fine meshes. Left: DEM. Middle: Lagrange FE. Right: penalised Crouzeix--Raviart FE.}
\label{fig:comparaison dynamique energies torsion} 
\end{figure}

The comparison between the methods is performed by considering the displacement and the velocity of the point of coordinates $(0.9R,0,\frac16 L)$ over the simulation time $T$.
The results are reported in Figure \ref{fig:comparaison dynamique 3}. 
The angular velocity plot indicates  for times $t\le 0.2$s the presence of elastic waves with a small magnitude travelling through the beam. The wave of larger amplitude arriving afterwards is a plastic wave whose velocity is ten times smaller than the elastic waves since $E_t = \frac{E}{100}$.
We notice that the value for the simulation time is too long for the simulation to remain 
physically relevant within the small strain assumption owing to the large value reached by the angular displacement of the reference point. However this setting allows us to reach substantial amounts of plastic dissipation and thereby to probe the robustness of the space semi-discretization methods with respect to incompressibility. Recall that the remanent plastic strain tensor is trace-free, so that the stress tensor is nearly deviatoric in the entire beam at the end of the simulation. Such a situation is challenging for the $P^1$-Lagrange FEM since this method is  known to lock in the incompressible limit.
To highlight the volumetric locking incurred by Lagrange FE,  Figure~\ref{fig:comparaison torsion plasticite} displays at the time $t=\frac12 T$ 
the trace of the stress tensor predicted
by Lagrange FE (fine mesh) and penalised Crouzeix--Raviart FE (coarse mesh), for a similar number of dofs. 
For Lagrange FE, significant oscillations are visible in the whole beam (the amplitude of these 
oscillations is about ten times the maximal value of the von Mises equivalent stress). Also, the amplitudes of the oscilliations of the trace tensor are about four times larger than for penalised Crouzeix--Raviart FE.
Figure~\ref{fig:comparaison dynamique energies torsion} reports the energies on the fine meshes for the DEM, Lagrange FE and penalised Crouzeix--Raviart FE. First, we notice as in the previous test case the prefect balance of the work of external loads with the different components of the mechanical energy.
The orders of magnitude of the energies and plastic dissipations are similar for the three methods. However, a significant discrepancy in the plastic dissipation can be observed for Lagrange FE with respect to penalised Crouzeix--Raviart FE and DEM which both give a plastic dissipation similar to the reference computation.

\begin{figure}[htp]
\centering
\subfloat{
\includegraphics[width=0.33\textwidth]{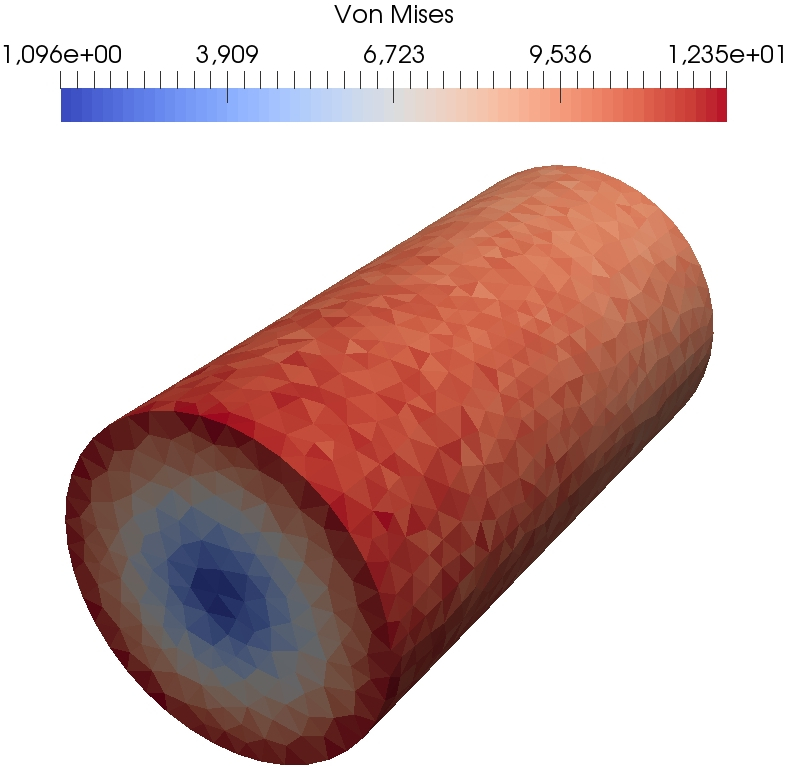}
}
\subfloat{
\includegraphics[width=0.33\textwidth]{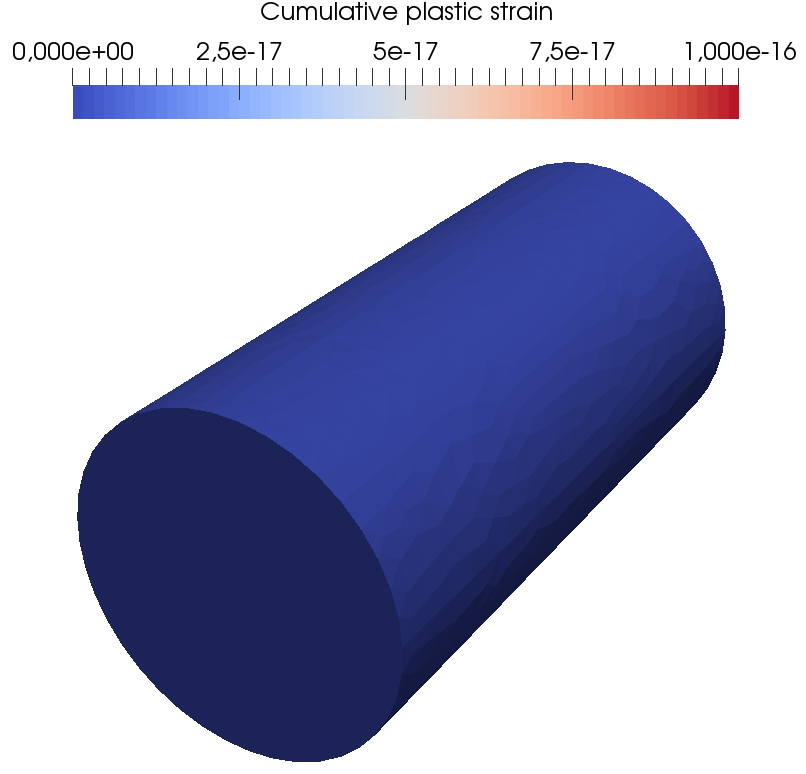}
}
\subfloat{
\includegraphics[width=0.33\textwidth]{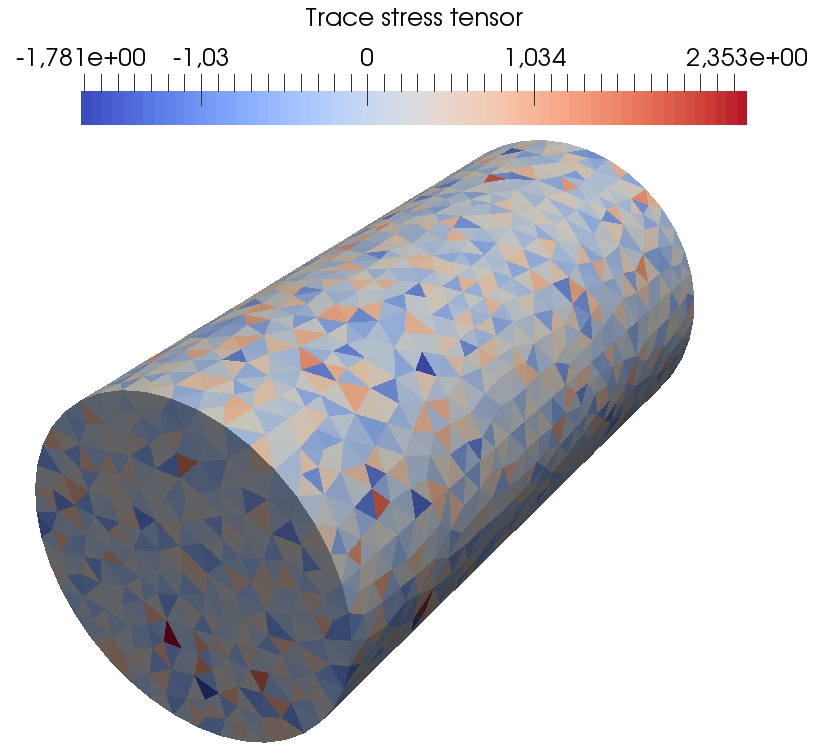}
}\\

\subfloat{
\includegraphics[width=0.33\textwidth]{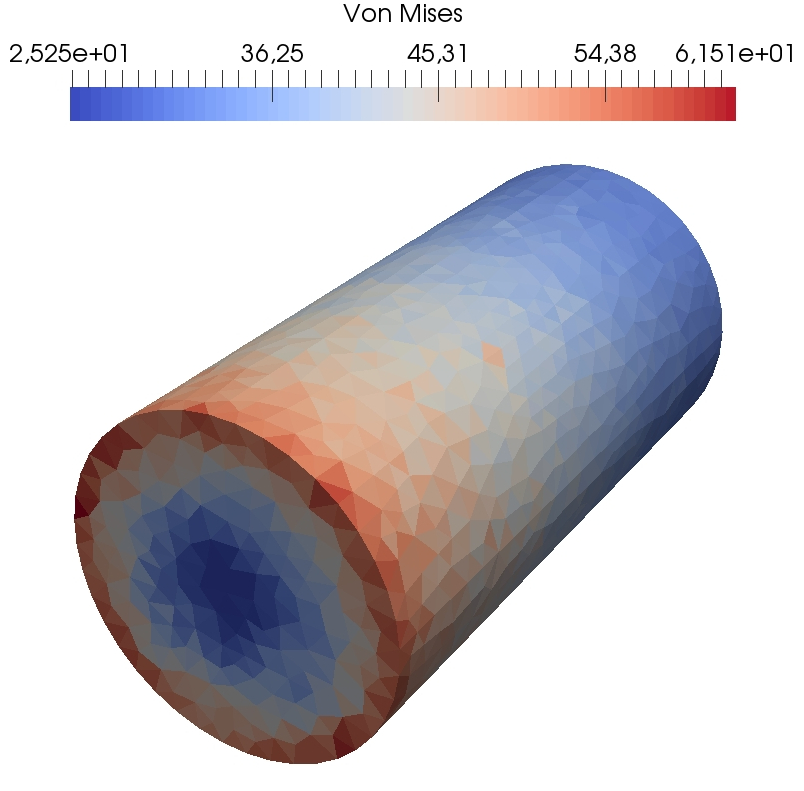}
}
\subfloat{
\includegraphics[width=0.33\textwidth]{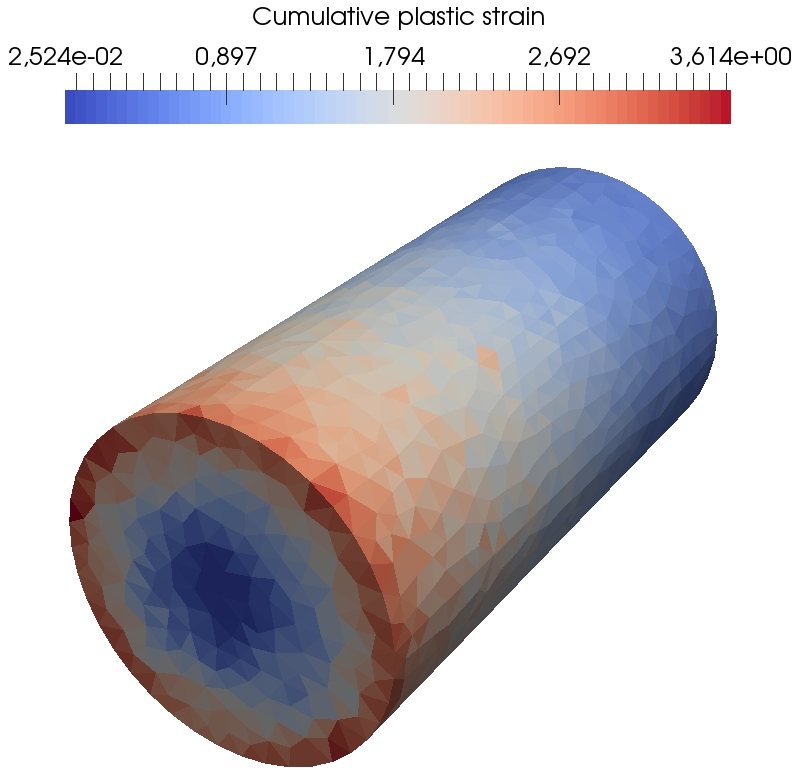}
}
\subfloat{
\includegraphics[width=0.33\textwidth]{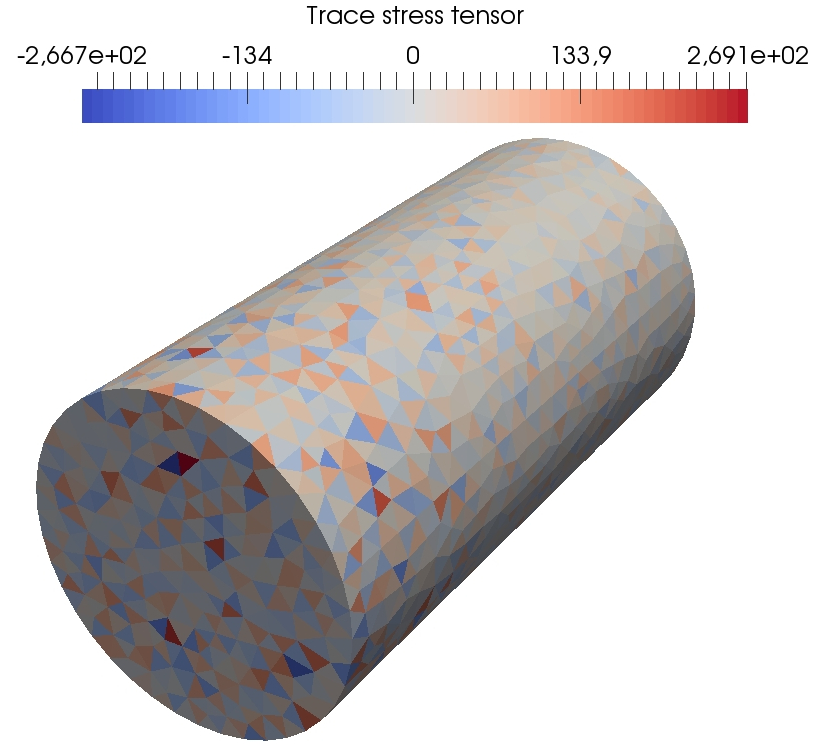}
}\\

\subfloat{
\includegraphics[width=0.33\textwidth]{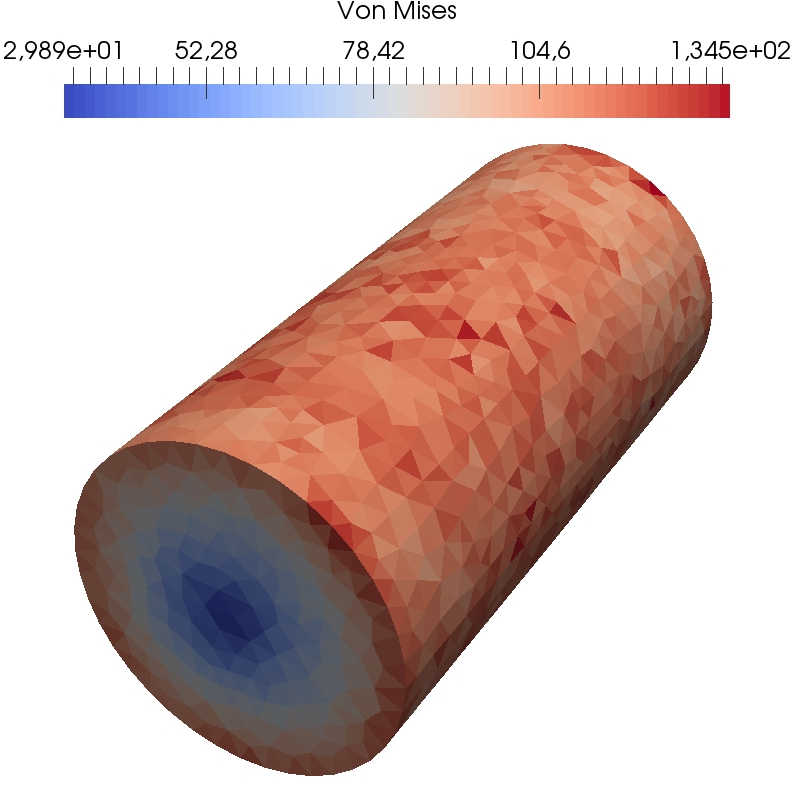}
}
\subfloat{
\includegraphics[width=0.33\textwidth]{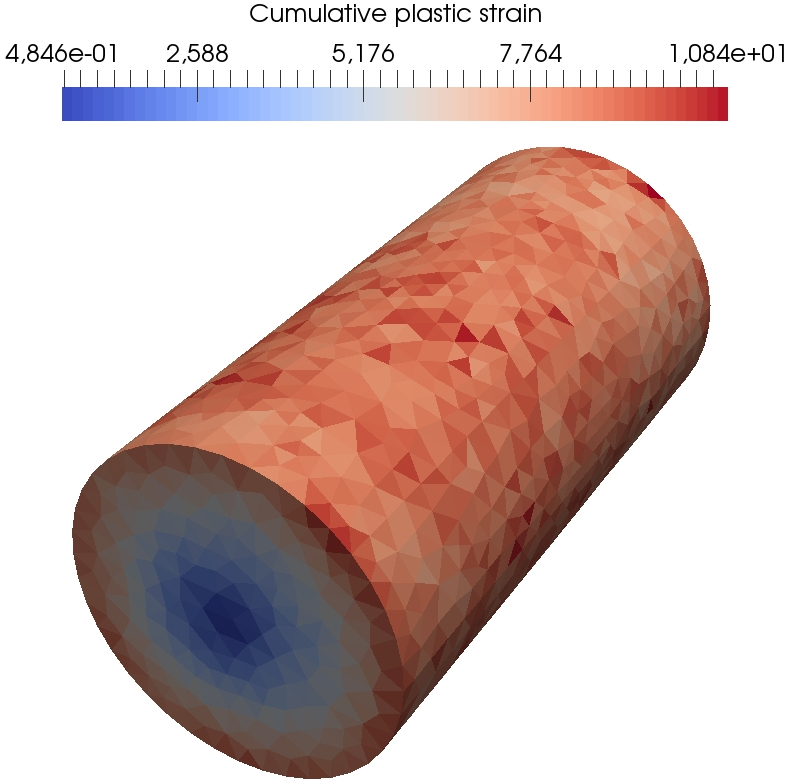}
}
\subfloat{
\includegraphics[width=0.33\textwidth]{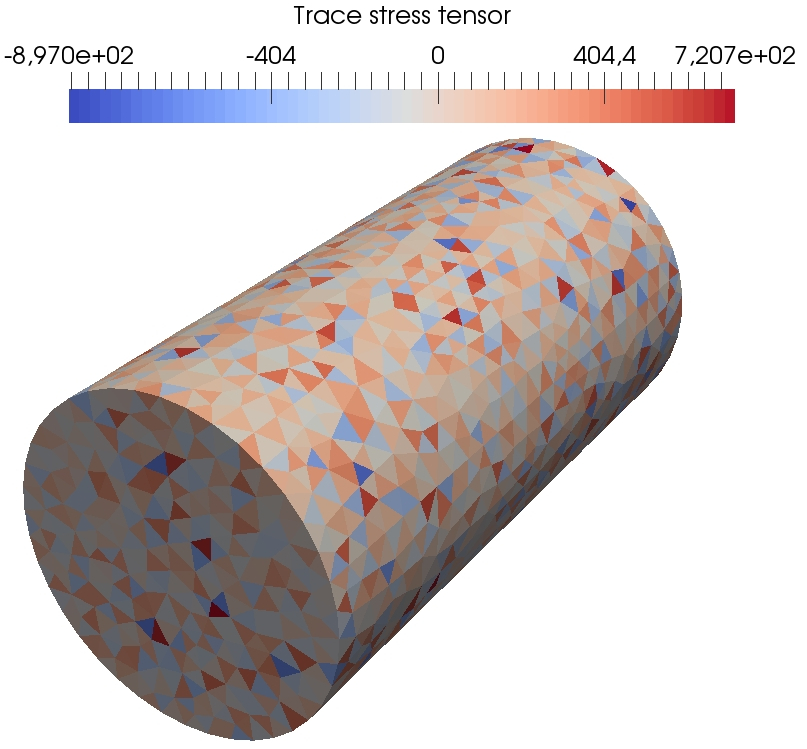}
}\\

\caption{Beam in dynamic torsion: DEM on the fine mesh. Von Mises equivalent stress (left column), $p$ (middle column) and $\mathrm{tr}(\sigma)$ (right column) at $t=\frac{1}{10}T$ (top line), $t=\frac12 T$ (middle line) and $t=T$ (bottom line).}
\label{fig:results computation torsion}
\end{figure}

Figure~\ref{fig:results computation torsion} presents some further results of the DEM
computations on the fine mesh so as to visualize at three different times during the
simulation the spatial localization of the von Mises equivalent stress, the cumulated plastic strain, and the trace of the stress tensor.
We can see in the first row that the von Mises stress is nonzero in the entire beam and thus the elastic waves have already travelled through the entire beam whereas the cumulative plastic strain is still zero and thus no plastic flow has occurred. In the second row, we can see that a plastic wave has started to propagate from one end of the beam. In the last row, we see that the plastic wave has reached the other side of the beam at the end of the simulation.
Also, regarding the robustness with respect to $\nu \to 0.5$, the magnitude of the oscillations of the trace of the stress tensor is similar to penalised Crouzeix--Raviart FE. Indeed for the DEM, the extreme values of the trace of the stress tensor are $-3 \cdot 10^2 \mathrm{Pa}$ and $3 \cdot 10^2 \mathrm{Pa}$ at $t=\frac{T}{2}$ and $-9 \cdot 10^2 \mathrm{Pa}$ and $7 \cdot 10^2 \mathrm{Pa}$ at $t=T$ whereas for penalised Crouzeix--Raviart FE the extermes values are $-2 \cdot 10^2 \mathrm{Pa}$ and $2 \cdot 10^2 \mathrm{Pa}$ at $t=\frac{T}{2}$ and $-6 \cdot 10^2 \mathrm{Pa}$ and $6 \cdot 10^2 \mathrm{Pa}$ at $t=T$.
Comparatively, for the $P^1$-Lagrange FE computations, the extreme values of the trace of the stress tensor are $-8 \cdot 10^2 \mathrm{Pa}$ and $9 \cdot 10^2 \mathrm{Pa}$ for $t=\frac{T}{2}$ and $-2 \cdot 10^3 \mathrm{Pa}$ and $2 \cdot 10^3 \mathrm{Pa}$ for $t=T$.

\section{Conclusion}
\label{sec:conclusions}

We have presented a new Discrete Element Method which is a variational discretization of a Cauchy continuum and which only requires continuum macroscopic parameters as the Young modulus and the Poisson ratio for its implementation. The displacement degrees of freedom are attached to the barycentres of the mesh cells and to the boundary vertices. The key idea is to reconstruct displacements on the mesh facets and then to use a discrete Stokes formula to devise a piecewise constant gradient and linearized strain reconstructions. A simple geometric pre-processing has been devised to ensure that for almost all the mesh facets, the reconstruction is based on an interpolation (rather than extrapolation) formula and we have shown by numerical experiments that this choice can produce highly beneficial effects in terms of the largest eigenvalue of the stiffness matrix, and thus on the time step restriction within explicit time-marching schemes. Moreover, in the case of elasto-plastic behavior, the internal variables for plasticity are piecewise-constant in the mesh cells. The scheme has been tested on quasi-static and dynamic test cases using a second-order, explicit, energy-conserving time-marching scheme. Future work can include adapting the present framework to dynamic cracking and fragmentation as well as to Cosserat continua. An extension to large strain solid dynamics could also be performed by working in the reference configuration.

\subsection*{Acknowledgements}
The authors would like to thank K. Sab (Navier, ENPC) and J.-P. Braeunig and L. Aubry (CEA) for fruitful discussions. The authors would also like to thank J. Bleyer (Navier, ENPC) for his help in dealing with the Fenics implementations. The PhD fellowship of the first author was supported by CEA.

\subsection*{Conflict of interest}
The authors declare no conflict of interest.

\bibliographystyle{plain} 
\bibliography{bib_article} 
\end{document}